\documentclass[12pt]{article}
\usepackage{amssymb}
\usepackage{mathrsfs}
\usepackage[utf8]{inputenc}
\usepackage[top=72pt,bottom=96pt,left=72pt,right=72pt]{geometry}
\def\1{\mathbf{1}}
\def\A{\mathbb{A}}
\def\at{\mathbf{at}}
\def\Aut{\mathrm{Aut}}
\def\black{\mathrm{black}}

\def\cls{\mathrm{closure}}
\def\col{\mathbf{c}}
\def\conn{\mathrm{conn}}
\def\Cross{\mathrm{Cross}}
\def\Curv{\mathrm{Curv}}
\def\Der{\mathrm{Der}}
\def\dw{\flat}
\def\Edge{\mathrm{Edge}}
\def\Ein{\mathrm{Ein}}
\def\End{\mathrm{End}}

\def\ext{\mathrm{ext}}
\def\Flag{\mathrm{Flag}}
\def\FS{\mathrm{FS}}
\def\Gr{\mathrm{Gr}}
\def\id{\mathrm{id}}
\def\IHX{\mathrm{IHX}}
\def\im{\mathrm{im}}
\def\Inv{\mathrm{Inv}}
\def\K{\mathscr{R}}
\def\KM{\mathbf{CM}}
\def\LP{\Lambda}
\def\N{\mathbb{N}}
\def\new{\mathrm{new}}
\def\OG{\mathbf{O}}
\def\old{\mathrm{old}}
\def\pf{\mathrm{pf}}
\def\pfill{\vskip6pt plus3pt minus2pt\noindent}
\def\pr{\mathrm{pr}}
\def\proof{\noindent\textbf{Proof:}\quad}
\def\qed{\ensuremath{\hfill\Box}\pfill}
\def\R{\mathbb{R}}
\def\red{\mathrm{red}}
\def\Ric{\mathrm{Ric}}
\def\S{\mathrm{Sym}}
\def\sec{\mathrm{sec}}
\def\Sec{\mathrm{Sec}}
\def\sgn{\mathrm{sgn}}
\def\so{\mathfrak{so}}
\def\span{\mathrm{span}}
\def\tr{\mathrm{tr}}
\def\U{\mathscr{U}}
\def\up{\sharp}
\def\Vert{\mathrm{Vert}}
\def\vol{\mathrm{vol}}
\def\Vol{\mathrm{Vol}}
\def\X{\mathfrak{X}}
\def\<#1,#2>{\langle\,#1,\,#2\,\rangle}
\newtheorem{Lemma}{Lemma}[section]
\newtheorem{Remark}[Lemma]{Remark}

\newtheorem{Theorem}[Lemma]{Theorem}
\newtheorem{Corollary}[Lemma]{Corollary}
\newtheorem{Definition}[Lemma]{Definition}

\newenvironment{Reference}[1]{\pfill\textbf{#1} \textit\bgroup}{\egroup\par}
\title{A Graphical Calculus for Stable Curvature Invariants}
\author{Gregor Weingart\footnote{Instituto de Matemáticas, Universidad
 Nacional Autónoma de México, Avenida Universidad s/n, Colonia
 Lomas de Chamilpa, 62210 Cuernavaca, MEXIQUE; \texttt{gw@matcuer.unam.mx}.}}
\begin{document}
\maketitle
\begin{abstract}
 \noindent
 In this article we develop a graphical calculus for stable invariants of
 Riemannian mani\-folds akin to the graphical calculus for Rozansky--Witten
 invariants for hyperkähler manifolds; based on interpreting trivalent
 graphs with colored edges as stably invariant polynomials on the space
 of algebraic curvature tensors. In this graphical calculus we describe
 explicitly the Pfaffian polynomials $(\,\pf_n\,)_{n\,\in\,\N_0}$ central
 to the Theorem of Chern--Gauß--Bonnet and the normalized moment polynomials
 $(\,\Psi^\circ_n\,)_{n\,\in\,\N_0}$ calculating the moments of sectional
 curvature considered as a random variable on the Graßmannian of planes.
 Eventually we illustrate the power of this graphical calculus by deriving
 a curvature identity for compact Einstein manifolds of dimensions greater
 than $2$ involving the Euler characteristic, the third moment of sectional
 curvature and the $L^2$--norm of the covariant derivative of the curvature
 tensor. A model implementation of this calculus for the computer algebra
 system Maxima is available \cite{w2}.
 \\[-10pt]
 \begin{center}
  \begin{tabular}{lp{268pt}}
   \textit{MSC (2020):}& 53--08, 53C25, 53E20.
   \\
   \textit{Keywords:}& Riemannian invariants, graph algebras, Einstein metrics.
  \end{tabular}
 \end{center}
\end{abstract}
\section{Introduction}
\label{intro}
 Riemannian manifolds are studied extensively in Differential Geometry from
 different points of view, among which is the classical topic of classifying
 Riemannian manifolds up to homeomorphisms, diffeomorphisms or isometries.
 Riemannian invariants provide a direct method to distinguish non--isometric
 manifolds, such invariants can be constructed for example by integrating a
 scalar valued polynomial $\psi$ in algebraic curvature tensors over the
 manifold in question. In order to have the integrand well--defined independent
 of the choice of coordinates $\psi$ needs to be a polynomial on the space
 $\Curv^-T$ of algebraic curvature tensors invariant under the orthogonal
 group $\OG(\,T,\,g\,)$. Tabulating the dimensions of the spaces of such
 invariant polynomials in dependence on their degree $n\,\in\,\N_0$ and
 the dimension $m\,\in\,\N$ of the euclidean vector space with the help
 of the computer algebra system LiE we obtain:
 $$
  \begin{array}{r|rrrrrrrrrrrrrr}
   \begin{picture}(12,12)
    \put( 1.0,11.5){\line(+1,-1){16}}
    \put(-2,- 1){$n$}
    \put( 5,+10){$m$}
   \end{picture}
   & 2 & 3 &  4 &  5 &  6 &  7 &  8 &  9 & 10 & 11 & 12 & 13 & 14 & 15
   \\
   \hline&&&&&&&&&&&&&&
   \\[-12pt]
   0 & 1 & 1 &  1 &  1 &  1 &  1 &  1 &  1 &  1 &  1 &  1 &  1 &  1 &  1
   \\[3pt]
   1 & \underline{1}
   & 1 &  1 &  1 &  1 &  1 &  1 &  1 &  1 &  1 &  1 &  1 &  1 &  1
   \\[3pt]
   2 & \underline{1} & \underline{2} & \underline{4}
   &  3 &  3 &  3 &  3 &  3 &  3 &  3 &  3 &  3 &  3 &  3
   \\[3pt]
   3 & \underline{1} & \underline{3} & \underline{9} & \underline{7}
   & \underline{8} &  8 &  8 &  8 &  8 &  8 &  8 &  8 &  8 &  8
   \\[3pt]
   4 & \underline{1} & \underline{4} & \underline{19} & \underline{20}
   & \underline{24} & \underline{25} & \underline{28}
   & 26 & 26 & 26 & 26 & 26 & 26 & 26
   \\[3pt]
   5 & \underline{1} & \underline{5} & \underline{39} & \underline{51}
   & \underline{83} & \underline{84} & \underline{101} & \underline{89}
   & \underline{90} & 90 & 90 & 90 & 90 & 90
   \\[3pt]
   6 & \underline{1} & \underline{7} & \underline{82} & \underline{150}
   & \underline{361} & \underline{359} & \underline{509} & \underline{403}
   & \underline{409} & \underline{408} & \underline{412}
   & 409 & 409 & 409
   \\[3pt]
   7 & \underline{1} & \underline{8} & \underline{151} & \underline{431}
   & \underline{1697} & \underline{1761} & \underline{3125} & \underline{2194}
   & \underline{2407} & \underline{2240} & \underline{2281} & \underline{2245}
   & \underline{2246} & 2246\hbox to0pt{\ .\hss}
  \end{array}
 $$
 Certainly the most interesting aspect of these dimensions is that they
 stabilize for fixed degree $n\,\in\,\N_0$ provided the dimension $m$ of
 the vector space $T$ is sufficiently large, more precisely the dimension
 of the space of invariant polynomials of degree $n\,\in\,\N_0$ on $\Curv^-T$
 is constant for $\dim\,T\,>\,2n$ with a last drop from $\dim\,T\,=\,2n$ to
 $\dim\,T\,=\,2n+1$ for even $n$. This phenomenon gives rise to the concept
 of stable curvature invariants.

 Stability is a somewhat vague concept, it can be made precise by introducing
 a suitably defined category $\KM$ of curvature models. Objects in this
 category are triples $(\,T,\,g,\,R\,)$ formed by an algebraic curvature
 tensor $R\,\in\,\Curv^-T$ over a euclidean vector space $(\,T,\,g\,)$,
 morphisms are the adjoints $F^*:\,T\longrightarrow\hat T$ of linear
 isometric maps $F:\,\hat T\longrightarrow T$ satisfying
 \begin{equation}\label{mkm}
  \hat R(\;F^*X,\,F^*Y;\,F^*U,\,F^*V\;)
  \;\;=\;\;
  R(\;X,\,Y;\,U,\,V\;)
 \end{equation}
 for all $X,\,Y,\,U,\,V\,\in\,T$. In turn a stable curvature invariant is a
 functor $\psi:\,\KM\longrightarrow\R$ from the curvature model category $\KM$
 to the category $\R$ of real numbers with only the identities as morphisms.
 In other words a stable curvature invariant $\psi$ associates a real number
 $\psi(R)\,\in\,\R$ to every curvature model $(\,T,\,g,\,R\,)$ regardless of
 its dimension $m\,\in\,\N_0$ with equality $\psi(R)\,=\,\psi(\hat R)$,
 whenever there exists a linear isometric map $F:\,\hat T\longrightarrow T$
 whose adjoint $F^*:\,T\longrightarrow\hat T$ satisfies the constraint
 (\ref{mkm}).

 The category $\KM$ of curvature models captures the essence of the
 stabilization phenomenon observed above: The adjoint $F^*:\,T\longrightarrow
 \hat T$ of an isometric map $F:\,\hat T\longrightarrow T$ equals the
 composition $F^*\,=\,F^{-1}\circ\pr$ of the orthogonal projection to
 the regular subspace $\im\,F\,\subseteq\,T$ followed by the isometry
 $F^{-1}:\,\im\,F\longrightarrow\hat T$. In turn the existence of a
 morphism $F^*:\,(\,T,\,g,\,R\,)\longrightarrow(\,\hat T,\,\hat g,\,\hat R\,)$
 in the category $\KM$ tells us via equation (\ref{mkm}) that the algebraic
 curvature tensor $R$ is essentially the Cartesian product $\hat R\oplus 0$
 of the algebraic curvature tensor $\hat R$ on $\im\,F\,\cong\,\hat T$ with
 the flat algebraic curvature tensor on $(\,\im\,F\,)^\perp$.

 \pfill
 The purpose of this article is to develop a graphical calculus akin to the
 calculus of Rozansky--Witten invariants \cite{nw} to describe the algebra of
 stable curvature invariants $\psi:\,\KM\longrightarrow\R$ such that the
 induced map $\Curv^-T\longrightarrow\R,\,R\longmapsto\psi(\,R\,),$ is a
 polynomial on the vector space $\Curv^-T$ of algebraic curvature tensors
 for every euclidean vector space $(\,T,\,g\,)$. Every such stable polynomial
 curvature invariant of degree $n$ is necessarily an $2n$--fold iterated sum
 over an orthonormal basis $E_1,\,\ldots,\,E_m$ for $T$ like the scalar
 curvature of degree $n\,=\,1$
 $$
  \kappa
  \;\;:=\;\;
  \sum_{\mu,\,\nu\,=\,1}^m
  R(\;E_\mu,\,E_\nu;\,E_\nu,\,E_\mu\;)
  \qquad\widehat=\qquad
  \frac14\;\;
  \raise-17pt\hbox{\begin{picture}(60,40)
   \put(15,20){\circle2}
   \put(45,20){\circle2}
   \multiput(17,20)(+2, 0){14}{\circle*1}
   \qbezier( 0,20)( 0,26)( 5,28)
   \qbezier( 5,28)(12,30)(15,21)
   \qbezier(60,20)(60,26)(55,28)
   \qbezier(55,28)(48,30)(45,21)
   \qbezier( 0,20)( 0,14)( 5,12)
   \qbezier( 5,12)(12,10)(15,19)
   \qbezier(60,20)(60,14)(55,12)
   \qbezier(55,12)(48,10)(45,19)
  \end{picture}}
 $$
 related to the Einstein--Hilbert functional or the norm square of the
 curvature tensor
 $$
  g_{\LP^2T^*\otimes\LP^2T^*}(\,R,\,R\,)
  \;\;:=\;\;
  \frac14\,\sum_{\mu,\,\nu,\,\alpha,\,\beta\,=\,1}^m
  R(\,E_\mu,\,E_\nu;\,E_\alpha,\,E_\beta\,)^2
  \qquad\widehat=\qquad
  \frac1{48}\;
  \raise-10pt\hbox{\begin{picture}(24,24)
   \put( 0, 0){\circle2}
   \put( 0,24){\circle2}
   \put(24, 0){\circle2}
   \put(24,24){\circle2}
   \multiput( 0, 2)( 0,+2){11}{\circle*1}
   \multiput(24, 2)( 0,+2){11}{\circle*1}
   \put(0.8,-0.7){\line(+1, 0){22.5}}
   \put(0.8, 0.7){\line(+1, 0){22.5}}
   \put(0.8,23.3){\line(+1, 0){22.5}}
   \put(0.8,24.7){\line(+1, 0){22.5}}
  \end{picture}}
 $$
 of degree $n\,=\,2$. In order to encode the contraction pattern of these
 $2n$--fold iterated sums over orthonormal bases in graphs we consider
 trivalent graphs, possibly with loops or multiple edges between vertices,
 and with edges colored red and black, rendered in this article as dotted
 and solid lines respectively, such that every vertex is adjacent to exactly
 one red edge. Every black edge corresponds to one sum over an orthonormal
 basis, while every red edge connects two different vertices and corresponds
 to a copy of the sectional curvature tensor $\Sec\,\in\,\Curv^+T$ associated
 to the algebraic curvature tensor $R\,\in\,\Curv^-T$. A concise formulation
 of this construction can be found in Definition \ref{slci}.
 
 Having outlined the general idea of how to convert a colored trivalent
 graph into a stable curvature invariant we define the graph algebra
 $\A^\bullet$ as the convolution algebra $\R\,\Gamma^\bullet$ of the graded
 monoid $\Gamma^\bullet$ of isomorphism classes of colored trivalent graphs
 under the disjoint union product. By construction this algebra comes along
 with an algebra homomorphisms
 $$
  \Inv_{(\,T,\,g\,)}:\;\;\A^\bullet\;\longrightarrow\;
  [\;\S^\bullet(\,\Curv^-T\,)^*\,]^{\OG(\,T,\,g\,)},\qquad
  [\,\gamma\,]\;\longmapsto\;[\,\gamma\,]\ ,
 $$
 for every euclidean vector space $T$. It turns out that this algebra
 homomorphism factorizes over the quotient $\overline\A^\bullet$ of the
 algebra $\A^\bullet$ of colored trivalent graphs by the ideal of
 IHX--relations
 $$
  \raise-17pt
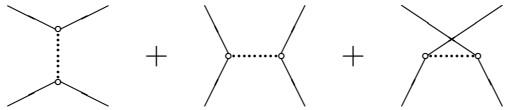
\hbox{\begin{picture}(40,40)
   \put(20,10){\circle2}
   \put(20,30){\circle2}
   \multiput(20,12)( 0,+2){9}{\circle*1}
   \put( 1, 1){\line(+2,+1){18}}
   \put(39, 1){\line(-2,+1){18}}
   \put( 1,39){\line(+2,-1){18}}
   \put(39,39){\line(-2,-1){18}}
  \end{picture}}
  \;\;\;+\;\;\;
  \raise-17pt\hbox{\begin{picture}(40,40)
   \put(10,20){\circle2}
   \put(30,20){\circle2}
   \multiput(12,20)(+2, 0){9}{\circle*1}
   \put( 1, 1){\line(+1,+2){9}}
   \put(39, 1){\line(-1,+2){9}}
   \put( 1,39){\line(+1,-2){9}}
   \put(39,39){\line(-1,-2){9}}
  \end{picture}}
  \;\;\;+\;\;\;
  \raise-17pt\hbox{\begin{picture}(40,40)
   \put(10,20){\circle2}
   \put(30,20){\circle2}
   \multiput(12,20)(+2, 0){9}{\circle*1}
   \put( 1, 1){\line(+1,+2){9}}
   \put(39, 1){\line(-1,+2){9}}
   \put( 1,39){\line(+3,-2){28}}
   \put(39,39){\line(-3,-2){28}}
  \end{picture}}
 $$
 on red edges, which arose historically in the study of invariants of knots;
 in the present context they reflect the first Bianchi identity. The reduced
 algebra $\overline\A^\bullet$ of colored trivalent graphs modulo the
 IHX--relations equals the algebra of stable curvature invariants:

 \begin{Theorem}[Stable Algebra Isomorphism]
 \hfill\label{iso}\break
  The algebra homomorphism from the reduced graph algebra $\overline
  \A^\bullet$ to the algebra of invariant polynomials on the space of
  algebraic curvature tensors over a euclidean vector space $T$
  $$
   \overline\Inv_{(\,T,\,g\,)}:\;\;\overline\A^\bullet\;\longrightarrow
   [\;\S^\bullet(\,\Curv^-T\,)^*\;]^{\OG(\,T,\,g\,)}
  $$
  induces isomorphisms $\overline\A^n\stackrel\cong\longrightarrow
  [\,\S^n(\,\Curv^-T\,)^*\,]^{\OG(\,T,\,g\,)}$ in all degrees
  $n\,<\,\frac12\dim\,T$.
 \end{Theorem}

 \pfill
 Unluckily we will not even outline the proof of this theorem in this
 article, mainly because it involves quite a lot of the representation
 theory of orthogonal and symplectic Lie algebras. However we intend to
 be more specific about Theorem \ref{iso} in a future publications, focussing
 for the time being on the description of the resulting graphical calculus.
 In case of doubt the reader may easily verify Theorem \ref{iso} directly
 in small degrees: Using the general relation between the Hilbert function
 of a free commutative associative graded algebra $\overline\A^\bullet$
 and the numbers $m^{\overline\A}_d\,\in\,\N_0$ of generators of
 $\overline\A^\bullet$ necessary in degree $d\,\in\,\N$
 \begin{equation}\label{pseries}
  \sum_{n\,>\,0}\Big(\;\sum_{d\,|\,n}d\,m^{\overline\A}_d\;\Big)\,t^n
  \;\;=\;\; 
  t\frac d{dt}\;
  \ln\,\Big(\;\sum_{n\,\geq\,0}(\,\dim\;\overline\A^n\,)\,t^n\;\Big)
 \end{equation}
 we can deduce from the table of dimensions presented above that the algebra
 $\overline\A^\bullet$ of stable curvature invariants should be a free
 commutative graded algebra with $1,\,2,\,5,\,15,\,54,\,270$ and $1639$
 generators respectively of degrees $1,\,2,\,3,\,4,\,5,\,6$ and $7$; numbers
 of generators which coincide up to degree $4$ with the number of generators
 of the reduced graph algebra $\overline\A^\bullet$ calculated in Section
 \ref{graphs} as a direct consequence of Corollary \ref{gga}.

 \pfill
 Although characteristic numbers of compact manifolds like the Pontryagin
 numbers are polynomial curvature invariants, it turns out that they are
 not stable polynomial curvature invariants in the sense of this article,
 not to the least so, because the corresponding polynomial $\psi\,\in\,
 [\,\S^\bullet(\,\Curv^-T\,)^*\,]^{\mathbf{SO}(\,T,\,g\,)}$ is invariant
 under the special orthogonal subgroup $\mathbf{SO}(\,T,\,g\,)$, but not
 under the full orthogonal group $\OG(\,T,\,g\,)$. The sole exception to
 this rule is the Euler characteristic: We can write the Theorem of
 Chern--Gauß--Bonnet for every compact, not necessarily oriented Riemannian
 manifold $M$ of even dimension $m$ in the form
 $$
  \chi(\,M\,)
  \;\;=\;\;
  \frac1{(\,2\pi\,)^{\frac m2}}\int_M\pf_{\frac m2}(\,R\,)\,|\,\vol_g\,|
 $$
 with a sequence $(\,\pf_n\,)_{n\,\in\,\N_0}$ of elements of degree $n$ in
 algebra $\A^\bullet$ of colored trivalent graphs. This sequence of Pfaffian
 polynomials will be studied in more detail in Section \ref{special} together
 with the sequence $(\,\Psi^\circ_n\,)_{n\,\in\,\N_0}$ of normalized moment
 polynomials, which calculate the moments of the sectional curvature considered
 as a random variable on the Graßmannian of planes. Both sequences have
 strikingly similar expansions in the power series completion
 \begin{eqnarray}
  \label{pfn}
  \sum_{n\,\geq\,0}\pf_n
  &=&
  \exp\left(\;\sum_{[\,\gamma\,]\,\in\,\Gamma^\bullet_\conn}
  \frac{(-1)^{e(\gamma)}}{6^{n(\gamma)}}\;
  \frac{2^{g(\gamma)}}{\#\overline\Aut\;\gamma}\;
  [\,\gamma\,]\;\right)
  \\[2pt]
  \label{psin}
  \sum_{n\,\geq\,0}\Psi^\circ_n
  &=&
  \exp\left(\;
   \sum_{{\scriptstyle[\,\gamma\,]\,\in\,\Gamma^\bullet_\conn}\atop
    {\scriptstyle\gamma_\black\mathrm{\;even\;cycles}}}
  (-1)^{n(\gamma)}\frac{2^{g(\gamma)}}{\#\overline\Aut\;\gamma}\;
  [\,\gamma\,]\;\right)
 \end{eqnarray}
 of the graph algebra $\A^\bullet$ derived in Lemmas \ref{efp} and \ref{efpsi}
 respectively, where $n(\gamma)\,:=\,\frac12\,\#\,\Vert\,\gamma$ is just the
 degree of the graph $\gamma$, whereas $e(\gamma)$ and $g(\gamma)$ denote
 the numbers of cycles of the bivalent black subgraph $\gamma_\black$ of
 even length and of length greater than $2$ respectively.

 \pfill
 In order to illustrate the power of the graphical calculus developed in this
 article we eventually consider in Section \ref{standard} the drastic algebraic
 simplifications in the reduced algebra $\overline\A^\bullet$ of colored
 trivalent graphs brought about by assuming the Riemannian manifold $M$ to
 be an Einstein manifold and use these simplification to prove the following
 curvature identity:

 \begin{Reference}{%
  Theorem \ref{ineq} (Cubic Curvature Identity for Einstein Manifolds)}
 \hfill\break
  For every compact connected Einstein manifold $M$ of dimension $m\,\geq\,3$
  with scalar curvature $\kappa\,\in\,\R$ the following identity of integrated
  stable curvature invariants of degree $3$ holds true:
  \begin{eqnarray*}
   \lefteqn{\int_M\pf_3(\,R\,)\,|\,\vol_g\,|
    \;-\;\frac1{40}\,\int_M\Psi^\circ_3(\,R\,)\,|\,\vol_g\,|
    \;+\;\frac2{15}\,|\!|\,\nabla R\,|\!|^2_{T^*\otimes\LP^2T^*\otimes\LP^2T^*}}
   \qquad\qquad
   &&
   \\
   &=&
   \kappa^3\,\frac{m^2-18m+40}{60\,m^2}\,\Vol(\,M,\,g\,)\;+\;
   \kappa\,\frac{3m-104}{30\,m}\,|\!|\,R\,|\!|^2_{\LP^2T^*\otimes\LP^2T^*}
  \end{eqnarray*}
 \end{Reference}

 \noindent
 In Section \ref{tensors} we provide a leisurely introduction to algebraic
 and sectional curvature tensors. Section \ref{graphs} is certainly the
 central section of this article and details the construction of stable
 curvature invariants from graphs and the construction of the graph algebras
 $\A^\bullet$ and $\overline\A^\bullet$. Graphs are evaluated combinatorially
 on algebraic curvature tensors of constant sectional curvature in Section
 \ref{values}. In Section \ref{special} we define the Pfaffian and normalized
 moment polynomials and present the combinatorial arguments behind the
 expansions (\ref{pfn}) and (\ref{psin}) of their generating series. Last
 but not least we establish Theorem \ref{ineq} in Section \ref{standard}.
 A model implementation for the computer algebra system Maxima of the
 functionalities of the graphical calculus presented below can be found
 under the link \cite{w2}.
\section*{Acknowledgements}
 Research into the construction of stable curvature invariants and the
 concurrent implementation of the graphs algebras $\A^\bullet$ and
 $\overline\A^\bullet$ in Maxima leading to numerous new insights
 was done while the author enjoyed a very pleasant research stay at
 the Max Planck Institut für Mathematik in Bonn, whose hospitality
 and generosity is gratefully acknowledged.
\section{Algebraic and Sectional Curvature Tensors}
\label{tensors}
 Algebraic curvature tensors are discussed in detail in many textbooks on
 differential geometry, a good reference is for example \cite{besse}. In
 this introductory section we will focus on the less well--known symmetric
 counterparts of algebraic curvature tensors, the sectional curvature tensors,
 with the aim to establish the equivalence between both ways to describe
 curvature. Sectional curvature tensors are however easier to deal with, a
 statement illustrated by the derivation of the polarization formula for
 algebraic curvature tensors given in this section. In the development of
 a graphical calculus for stable curvature invariants this simplicity of
 sectional compared to algebraic curvature tensors will be a crucial advantage.

 \pfill
 A euclidean vector space will be for the purpose of this article a finite
 dimensional vector space $T$ over $\R$ endowed with a non--degenerate,
 not necessarily positive definite symmetric bilinear form $g:\,T\times T
 \longrightarrow\R$ called its scalar product. An isometric map between
 euclidean vector spaces $T$ and $\hat T$ is a linear map $F:\,T\longrightarrow
 \hat T$ with the characteristic property that $\hat g(FX,FY)\,=\,g(X,Y)$ for
 all $X,\,Y\,\in\,T$, isometries are invertible isometric maps. Declaring for
 example the mutually inverse musical isomorphisms $\dw:\,T\longrightarrow T^*,
 \,X\longmapsto g(X,\,\cdot\,),$ and $\up\,:=\,\dw^{-1}$ to be isometries
 defines a scalar product $g^{-1}:\,T^*\times T^*\longrightarrow\R$ such
 that $g^{-1}(\alpha,\beta)\,:=\,\alpha(\beta^\up)$. The scalar products
 $g$ and $g^{-1}$ extend to symmetric and exterior powers of both $T$ and
 $T^*$ by using Gram's permanent or determinant respectively, for example:
 $$
  g_{\LP^kT^*}(\;\alpha_1\,\wedge\,\ldots\,\wedge\,\alpha_k,\;
  \beta_1\,\wedge\,\ldots\,\wedge\,\beta_k\;)
  \;\;:=\;\;
  \det\;\Big\{\;g^{-1}(\,\alpha_\mu,\beta_\nu\,)\;
  \Big\}_{\mu,\,\nu\,=\,1,\,\ldots,\,k}\ .
 $$
 Coming back to our main topic we define an algebraic curvature
 tensor on a euclidean vector space $T$ to be a quadrilinear form
 $R:\,T\times T\times T\times T\longrightarrow\R$, which is skew
 symmetric in the first $R(X,Y;U,V)\,=\,-\,R(Y,X;U,V)$ and second
 pair $R(X,Y;U,V)\,=\,-\,R(X,Y;V,U)$ of its arguments $X,\,Y,\,Z,
 \,U,\,V\,\in\,T$ and satisfies the so called first Bianchi identity:
 \begin{equation}\label{1BI}
  R(\;X,\,Y;\;Z,\,V\;)\;+\;R(\;Y,\,Z;\;X,\,V\;)\;+\;R(\;Z,\,X;\;Y,\,V\;)
  \;\;=\;\;
  0\ .
 \end{equation}
 Due to the presence of the musical isomorphisms $\dw$ and $\up$ algebraic
 curvature tensors can be interpreted alternatively as trilinear products
 $R:\,T\times T\times T\longrightarrow T,\,(\,X,\,Y;\,U\,)\longmapsto
 R_{X,\,Y}U$ on the vector space $T$ by setting $R_{X,\,Y}U\,:=\,R(X,Y;U,
 \,\cdot\,)^\up$ or even as $2$--forms with values in the Lie subalgebra
 $\so(\,T,\,g\,)\,\subseteq\,\End\,T$ of skew symmetric endomorphisms of
 $T$:
 $$
  R:\;\;T\;\times\;T\;\longrightarrow\;\so(\,T,\,g\,),\qquad
  (\,X,\,Y\,)\;\longmapsto\;R_{X,\,Y}\ .
 $$
 The latter interpretation is particularly interesting due to the vector space
 isomorphism
 $$
  \LP^2T\;\stackrel\cong\longrightarrow\;\so(\,T,\,g\,),
  \qquad X\,\wedge\,Y\;\longmapsto\;\Big(\;U\;\longmapsto\;
  g(\,X,\,U\,)\,Y\;-\;g(\,Y,\,U\,)X\;\Big)
 $$
 characterized completely by the identity $g_{\LP^2T}(\,\X,\,
 U\wedge V\,)\,=\,g(\,\X\,U,\,V\,)$ for every bivector
 $\X\,\in\,\LP^2T$ and all $U,\,V\,\in\,T$. In accordance with
 this isomorphism an algebraic curvature tensor $R$ can be interpreted
 as a bivector valued $2$--form on $T$, namely the $2$--form
 \begin{equation}\label{2v}
  R
  \;\;=\;\;
  \frac14\,\sum_{\mu,\,\nu,\,\alpha,\,\beta\,=\,1}^m
  R(\,E_\mu,\,E_\nu;\,E_\alpha,\,E_\beta\,)\;
  dE_\mu\wedge dE_\nu\;\otimes\;dE_\alpha^\up\wedge dE_\beta^\up\ ,
 \end{equation}
 where the sum is over an arbitrary pair of dual bases $\{\,E_\mu\,\}$ and
 $\{\,dE_\mu\,\}$ for the euclidean vector space $T$ and its dual $T^*$. Of
 course the latter sum simplifies somewhat for an orthonormal basis due to
 the equations $E_\alpha^\dw\,=\,\pm dE_\alpha$ and $dE_\alpha^\up\,=\,\pm
 E_\alpha$ characterizing orthonormal bases in general. Last but not least
 we want to recall the definitions of the Ricci tensor
 \begin{equation}\label{ric}
  \Ric(\;X,\,Y\;)
  \;\;:=\;\;
  \tr\Big(\;\;U\;\longmapsto\;R_{U,\,X}Y\;\;\Big)
  \;\;\stackrel!=\;\;
  \sum_{\mu\,=\,1}^mR(\;E_\mu,\,X;\,Y,\,dE_\mu^\up\;)
 \end{equation}
 and the scalar curvature $\kappa\,\in\,\R$ associated to an algebraic
 curvature tensor $R$:
 \begin{equation}\label{scal}
  \kappa
  \;\;:=\;\;
  \sum_{\nu\,=\,1}^m\Ric(\;E_\nu,\,dE_\nu^\up\;)
  \;\;=\;\;
  \sum_{\mu,\,\nu\,=\,1}^mR(\;E_\mu,\,E_\nu;\,dE_\nu^\up,\,dE_\mu^\up\;)\ .
 \end{equation}
 Sectional curvature tensors on a euclidean vector space $T$ are defined in
 complete analogy to algebraic curvature tensors as quadrilinear forms
 $S:\,T\times T\times T\times T\longrightarrow\R$, which are symmetric in
 the first and second pair $S(X,Y;U,V)\,=\,S(Y,X;U,V)\,=\,S(X,Y;V,U)$ of
 arguments and satisfy $S(X,X;X,V)\,=\,0$ for all $X,\,V\,\in\,T$; the
 latter becomes
 \begin{equation}\label{1BI'}
  S(\;X,\,Y;\;Z,\,V\;)\;+\;S(\;Y,\,Z;\;X,\,V\;)\;+\;S(\;Z,\,X;\;Y,\,V\;)
  \;\;=\;\;
  0
 \end{equation}
 upon polarization in $X,\,Y,\,Z\,\in\,T$. Algebraic curvature tensors
 are well--known to satisfy the symmetry in pairs $R(X,Y;U,V)\,=\,R(U,V;X,Y)$
 identically in their arguments:
 \begin{eqnarray*}
  2\,R(\,X,Y;U,V\,)
  &\!=\!&
  +\,R(\,X,Y;U,V\,)\,-\,R(\,X,Y;V,U\,)
  \\
  &\!=\!&
  -\,R(\,Y,U;X,V\,)\,-\,R(\,U,X;Y,V\,)\,+\,R(\,Y,V;X,U\,)\,+\,R(\,V,X;Y,U\,)
  \\
  &\!=\!&
  +\,R(\,Y,U;V,X\,)\,+\,R(\,V,Y;U,X\,)\,-\,R(\,X,U;V,Y\,)\,-\,R(\,V,X;U,Y\,)
  \\
  &\!=\!&
  -\,R(\,U,V;Y,X\,)\,+\,R(\,U,V;X,Y\,)
  \\
  &\!=\!&
  \,2\,R(\,U,V;X,Y\,)\ .
 \end{eqnarray*}
 Mutatis mutandis the reader may easily verify that sectional curvature
 tensors are symmetric in pairs $S(X,Y;U,V)\,=\,S(U,V;X,Y)$ as well.
 However it is more elegant to use the unpolarized first Bianchi identity
 $S(X,X;X,V)\,=\,0$ directly to obtain the equality 
 \begin{eqnarray*}
  0
  &=&
  \left.\frac d{dt}\right|_0S(\,X+tU,X+tU;\,X+tU,U\,)
  \,-\,\left.\frac d{dt}\right|_0S(\,U+tX,U+tX;\,U+tX,X\,)
  \\
  &=&
  2\,S(\,X,U;\,X,U\,)\;+\;S(\,X,X;\,U,U\,)
  \;-\;2\,S(\,U,X;\,U,X\,)\;-\;S(\,U,U;\,X,X\,)
  \\[2pt]
  &=&
  S(\,X,X;\,U,U\,)\;-\;S(\,U,U;\,X,X\,)
 \end{eqnarray*}
 establishing the restricted symmetry $S(X,X;U,U)\,=\,S(U,U;X,X)$ for all
 $X,\,U\,\in\,T$. The general symmetry in pairs $S(X,Y;U,V)\,=\,S(U,V;X,Y)$
 follows directly from its restricted version by means of the following
 polarization formula for sectional curvature tensors
 \begin{eqnarray}
  \nonumber\lefteqn{16\,S(\;X,Y;\,U,V\;)}
  &&
  \\
  \label{pol}
  &=&
  +\;S(\;X+Y,X+Y;\,U+V,U+V\;)\;-\;S(\;X-Y,X-Y;\,U+V,U+V\;)
  \\
  &&
  \nonumber
  -\;S(\;X+Y,X+Y;\,U-V,U-V\;)\;+\;S(\;X-Y,X-Y;\,U-V,U-V\;)
 \end{eqnarray}
 valid for all $X,\,Y,\,U,\,V\,\in\,T$, which is simply an iteration of the
 binomial polarization formula $4\,a(X,Y)\,=\,a(X+Y,X+Y)\,-\,a(X-Y,X-Y)$ for
 symmetric bilinear forms $a\,\in\,\S^2T^*$. By definition the sets of
 algebraic or sectional curvature tensors are natural subspaces $\Curv^-T$
 or $\Curv^+T$ respectively of the vector space $\bigotimes^4T^*$ of
 quadrilinear forms:
 
 \begin{Lemma}[Equivalence of Algebraic and Sectional Curvature]
 \hfill\label{eqasc}\break
  For every euclidean vector space $T$ the vector spaces $\Curv^-T$ and
  $\Curv^+T$ of algebraic and sectional curvature tensors on $T$ are
  naturally isomorphic via the mutually inverse isomorphisms $\Phi^+:
  \,\Curv^-T\longrightarrow\Curv^+T$ and $\Phi^-:\,\Curv^+T\longrightarrow
  \Curv^-T$ defined by:
  \begin{eqnarray*}
   (\,\Phi^+R\,)(\;X,\,Y;\,U,\,V\;)
   &:=&
   -\;2\;\Big(\;R(\;X,\,U;\,Y,\,V\;)\;+\;R(\;X,\,V;\,Y,\,U\;)\;\Big)
   \\
   (\,\Phi^-S\,)(\;X,\,Y;\,U,\,V\;)
   &:=&
   -\;\frac16\;\Big(\;S(\;X,\,U;\,Y,\,V\;)\;-\;S(\;X,\,V;\,Y,\,U\;)\;\Big)\ .
  \end{eqnarray*}
  The scalar factors in the definitions of $\Phi^+$ and $\Phi^-$ are a matter
  of taste up to their product being equal to $\frac13$. With our choice the
  sectional curvature tensor $\Sec\,:=\,\Phi^+R$ associated to an algebraic
  curvature tensor $R\,\in\,\Curv^-T$ satisfies for all $X,\,U\,\in\,T$ the
  identity:
  $$
   \frac14\;\Sec(\;X,\,X;\,U,\,U\;)
   \;\;=\;\;
   R(\;X,\,U;\,U,\,X\;)\ .
  $$
  In particular the inequalities $\Sec(\,X,X;\,U,U\,)\,\geq\,0$ or
  $\Sec(\,X,X;\,U,U\,)\,>\,0$ for all or for all linearly independent
  arguments $X,\,U\,\in\,T$ respectively characterize the cones of
  algebraic curvature tensors of non--negative or positive sectional
  curvature.
 \end{Lemma}

 \proof
 Leaving the straightforward verification of $\Phi^+R\,\in\,\Curv^+T$ and
 $\Phi^-S\,\in\,\Curv^-T$ for all $R\,\in\,\Curv^-T$ and $S\,\in\,\Curv^+T$
 to the reader we simply expand the definitions of the linear maps $\Phi^+$
 and $\Phi^-$ to find for example for every algebraic curvature tensor
 $R\,\in\,\Curv^-T$
 \begin{eqnarray*}
  \lefteqn{(\,\Phi^-\Phi^+R\,)(\;X,\,Y;\,U,\,V\;)}
  \quad
  &&
  \\
  &=&
  -\;\frac16\;\Big(\;(\,\Phi^+R\,)(\;X,\,U;\,Y,\,V\;)
  \;-\;(\,\Phi^+R\,)(\;X,\,V;\,Y,\,U\;)\;\Big)
  \\[2pt]
  &=&
  \hphantom{+\;}\frac13\;\Big(\;R(\,X,Y;\,U,V\,)\;+\;R(\,X,V;\,U,Y\,)
  \;-\;R(\,X,Y;\,V,U\,)\;-\;R(\,X,U;\,V,Y\,)\;\Big)
  \\[2pt]
  &=&
  \hphantom{+\;}\frac23\;R(\;X,\,Y;\,U,\,V\;)
  \;-\;\frac13\;\Big(\;R(\;V,\,X;\,U,\,Y\;)\;+\;R(\,X,\,U;\,V,\,Y\;)\;\Big)
  \\[2pt]
  &=&
  \hphantom{+\;}
  \frac23\;R(\;X,\,Y;\,U,\,V\;)\;+\;\frac13\;R(\;U,\,V;\,X,\,Y\;)
  \;\;=\;\;
  R(\;X,\,Y;\,U,\,V\;)
 \end{eqnarray*}
 using the first Bianchi identity (\ref{1BI}) and the symmetry in pairs
 in the last line. The slightly simpler analogous argument establishing
 $\Phi^+\Phi^-S\,=\,S$ for every sectional curvature tensor $S\,\in\,\Curv^+T$
 is omitted to reduce redundancy.
 \qed

 \noindent
 With $\Phi^+$ and $\Phi^-$ being inverse isomorphisms we can write every
 algebraic curvature tensor in the form $R\,=\,\Phi^-\Sec$ for its associated
 sectional curvature tensor $\Sec\,:=\,\Phi^+R$ and so
 $$
  24\,R(\,X,Y;\,U,V\,)
  \;\;=\;\;
  4\,\Sec(\,X,V;\,Y,U\,)\;-\;4\,\Sec(\,X,U;\,Y,V\,)
 $$
 for all arguments $X,\,Y,\,U,\,V\,\in\,T$. In particular the polarization
 formula (\ref{pol}) for sectional curvature tensors entails a polarization
 formula for algebraic curvature tensors $R\,\in\,\Curv^-T$
 \begin{eqnarray*}
  24R(X,Y;U,V)
  &\!=\!&
  +R(X+V,Y+U;Y+U,X+V)-R(X+U,Y+V;Y+V,X+U)
  \\
  &&
  -R(X+V,Y-U;Y-U,X+V)+R(X+U,Y-V;Y-V,X+U)
  \\
  &&
  -R(X-V,Y+U;Y+U,X-V)+R(X-U,Y+V;Y+V,X-U)
  \\
  &&
  +R(X-V,Y-U;Y-U,X-V)-R(X-U,Y-V;Y-V,X-U)
 \end{eqnarray*}
 for all $X,\,Y,\,U,\,V$. Likewise the definitions (\ref{ric}) and (\ref{scal})
 of the Ricci tensor $\Ric\,\in\,\S^2T^*$ and the scalar curvature $\kappa\,
 \in\,\R$ associated to an algebraic curvature tensor $R$ can be rewritten
 $$
  \Ric(\;X,\,Y\;)
  \;\;=\;\;
  \frac14\,\sum_{\mu\,=\,1}^m
  \Sec(\,E_\mu,dE_\mu^\up;\,X,Y\,)
  \qquad\quad
  \kappa
  \;\;=\;\;
  \frac14\sum_{\mu,\,\nu\,=\,1}^m
  \Sec(\,E_\mu,dE_\mu^\up;\,E_\nu,dE_\nu^\up\,)
 $$
 in terms of the sectional curvature tensor $\Sec\,=\,\Phi^+R$. In order to
 discuss a more important consequence of Lemma \ref{eqasc} we observe that
 the vector space $\S^2T^*\otimes\S^2T^*$ of quadrilinear forms symmetric
 in their first and second pair of arguments is spanned by the tensor products
 $$
  (\,a\,\otimes\,b\,)(\;X,\,Y;\,U,\,V\;)
  \;\;:=\;\;
  a(\;X,\,Y\;)\,b(\;U,\,V\;)
 $$
 of symmetric bilinear forms $a,\,b\,\in\,\S^2T^*$. In addition we can define
 the projection
 $$
  \pr:\;\;\S^2T^*\;\otimes\;\S^2T^*\;\longrightarrow\;\Curv^+T,
  \qquad S\;\longmapsto\;\pr\,S\ ,
 $$
 from the vector space $\S^2T^*\otimes\S^2T^*$ to the vector space $\Curv^+T$
 by setting
 \begin{eqnarray}
  \nonumber
  6\,(\,\pr\;S\,)(\;X,\,Y;\,U,\,V\;)
  &:=&
  +\;2\,S(\;X,\,Y;\,U,\,V\;)\;+\;2\,S(\;U,\,V;\,X,\,Y\;)
  \\
  \label{prd}
  &&
  -\;\hphantom{2\,}S(\;X,\,U;\,Y,\,V\;)\;-\;\hphantom{2\,}S(\;X,\,V;\,Y,\,U\;)
  \\
  \nonumber
  &&
  -\;\hphantom{2\,}S(\;Y,\,U;\,X,\,V\;)\;-\;\hphantom{2\,}S(\;Y,\,V;\,X,\,U\;)
 \end{eqnarray}
 for all arguments $X,\,Y,\,U,\,V\,\in\,T$; note that $\pr\;S\,=\,S$ for
 every sectional curvature tensor $S\,\in\,\Curv^+T$ due to the first Bianchi
 identity (\ref{1BI'}), the symmetry in pairs and the compensation factor $6$.
 On the other hand the image of every quadrilinear form $S\,\in\,\S^2T^*
 \otimes\S^2T^*$ symmetric in its first and second pair of arguments is a
 well--defined sectional curvature tensor $\pr\;S\,\in\,\Curv^+T$, because
 the right hand of equation (\ref{prd}) is evidently symmetric under
 $X\,\leftrightarrow\,Y$ and $U\,\leftrightarrow\,V$ while vanishing
 for $X\,=\,Y\,=\,U$. With $\pr$ being a projection and thus surjective
 we conclude that the vector space $\Curv^-T$ of algebraic curvature
 tensors on a euclidean vector space $T$ is spanned by the Nomizu--Kulkarni
 products
 \begin{equation}\label{nk}
  \begin{array}{l}
   (\,a\,\times\,b\,)(\;X,\,Y;\,U,\,V\;)
   \\[3pt]
   \quad\;:=\;\;
   a(X,U)\,b(Y,V)
   \,-\,a(X,V)\,b(Y,U)
   \,-\,a(Y,U)\,b(X,V)
   \,+\,a(Y,V)\,b(X,U)
  \end{array}
 \end{equation}
 of symmetric bilinear forms $a,\,b\,\in\,\S^2T^*$, because we find by
 expanding definition (\ref{prd})
 \begin{eqnarray*}
  \lefteqn{-\;36\,(\,\Phi^-\circ\pr\,)(\,a\,\otimes\,b\,)
   (\;X,\,Y;\,U,\,V\;)}
  \quad
  &&
  \\
  &=&
  +\;6\,\pr(\,a\,\otimes\,b\,)(\;X,\,U;\,Y,\,V\;)
  \;-\;6\,\pr(\,a\,\otimes\,b\,)(\;X,\,V;\,Y,\,U\;)
  \\
  &=&
  +\;2\,a(X,U)\,b(Y,V)\;+\;2\,a(Y,V)\,b(X,U)
  \;-\;2\,a(X,V)\,b(Y,U)\;-\;2\,a(Y,U)\,b(X,V)
  \\
  &&
  -\;\hphantom{2\,}a(X,Y)\,b(U,V)\;-\;\hphantom{2\,}a(X,V)\,b(U,Y)
  \;+\;\hphantom{2\,}a(X,Y)\,b(V,U)\;+\;\hphantom{2\,}a(X,U)\,b(V,Y)
  \\
  &&
  -\;\hphantom{2\,}a(U,Y)\,b(X,V)\;-\;\hphantom{2\,}a(U,V)\,b(X,Y)
  \;+\;\hphantom{2\,}a(V,Y)\,b(X,U)\;+\;\hphantom{2\,}a(V,U)\,b(X,Y)
  \\
  &=&
  +\;3\,(\,a\,\times\,b\,)(\;X,\,Y;\,U,\,V\;)
 \end{eqnarray*}
 for all $X,\,Y,\,U,\,V\,\in\,T$. Put differently the Nomizu--Kulkarni
 product $a\times b\,\in\,\Curv^-T$ of two symmetric bilinear forms
 $a,\,b\,\in\,\S^2T^*$ equals the algebraic curvature tensor corresponding
 to the sectional curvature tensor $\pr(\,a\otimes b\,)\,\in\,\Curv^+T$
 up to the scalar factor $-12$:
 \begin{equation}
  \Phi^+(\,a\,\times\,b\,)
  \;\;=\;\;
  -\;12\;\pr(\,a\,\otimes\,b\,)\ .
 \end{equation}
 Taking a closer look at definition (\ref{nk}) we see that the Nomizu--Kulkarni
 product is commutative $a\,\times\,b\,=\,b\,\times\,a$, hence we may polarize
 $a\,\times\,b\,=\,\frac14\,(a+b)\,\times\,(a+b)\,-\,\frac14\,(a-b)\,\times\,
 (a-b)$ to argue that the vector space $\Curv^-T$ is actually spanned by
 Nomizu--Kulkarni squares:
 \begin{equation}\label{nkspan}
  \Curv^-T
  \;\;=\;\;
  \span_\R\{\;\;a\times a\;\;|\;\;
  a\,\in\,\S^2T^*\textrm{\ symmetric bilinear form}\;\;\}
  \;\;\subseteq\;\;
  {\textstyle\bigotimes^4}T^*\ .
 \end{equation}
 Somewhat better every algebraic curvature tensor $R\,\in\,\Curv^-T$ can be
 written as a sum of Nomizu--Kulkarni squares of symmetric bilinear forms
 $a_1,\,\ldots,\,a_r\,\in\,\S^2T^*$ of rank two:
 \begin{equation}\label{rsum}
  R
  \;\;=\;\;
  a_1\times a_1\;+\;a_2\times a_2\;+\;\ldots\;+\;a_r\times a_r\ .
 \end{equation}
 The argument relies on the following weak form of Sylvester's Theorem of
 Inertia: For every symmetric bilinear form $a\,\in\,\S^2T^*$ on a finite
 dimensional vector space $T$ over $\R$ there exists a unique tuple
 $(\,p,\,n\,)\,\in\,\N_0^2$ called the signature of $a$ such that $a$ can
 be written as a sum of signed symmetric squares of linearly independent
 forms $\alpha_1,\,\ldots,\,\alpha_p,\,\beta_1,\,\ldots,\,\beta_n\,\in\,T^*$
 $$
  a
  \;\;=\;\;
  (\;{\textstyle\frac12}\,\alpha_1^2\;+\;\ldots\;+\;{\textstyle\frac12}
  \,\alpha_p^2\;)\;-\;(\;{\textstyle\frac12}\,\beta_1^2\;+\;\ldots\;+\;
  {\textstyle\frac12}\,\beta_n^2\;)\ ;
 $$
 by convention the square $\frac12\alpha^2\,\in\,\S^2T^*$ of a linear form
 $\alpha\,\in\,T^*$ denotes the symmetric bilinear form $\alpha\otimes\alpha:
 \,T\times T\longrightarrow\R,\,(X,Y)\longmapsto\alpha(X)\alpha(Y)$. In
 particular a symmetric bilinear form $a\,\in\,\S^2T^*$ is non--degenerate
 and hence defines a scalar product on $T$, if and only if its rank $p+n\,=\,m$
 equals the dimension of $T$ so that the $\alpha_1,\,\ldots,\,\alpha_p,\,
 \beta_1,\,\ldots,\,\beta_n$ become a basis of $T^*$, namely an orthonormal
 basis with respect to the dual scalar product $a^{-1}$.

 According to equation (\ref{nkspan}) every algebraic curvature tensor
 $R\,\in\,\Curv^-T$ can be expanded into a finite sum of scaled
 Nomizu--Kulkarni squares $\lambda\,a\times a$ with $a\,\in\,\S^2T^*$
 and a scalar $\lambda\,\in\,\R\,\setminus\,\{0\}$. Replacing
 $a\,\rightsquigarrow\,\sqrt{|\lambda|}a,\;\lambda\,\rightsquigarrow
 \,\frac\lambda{|\lambda|}$ if necessary we may assume without loss of
 generality that all scalars occurring in this expansion of $R$ equal
 $\lambda\,=\,\pm1$ so that:
 $$
  R
  \;\;=\;\;
  \pm\;a_1\,\times\,a_1\;\pm\;a_2\,\times\,a_2\;\pm\;\ldots\;
  \pm\;a_r\,\times\,a_r\ .
 $$
 Sylvester's Theorem of Inertia allows us to expand each of the symmetric
 bilinear forms $a_1,\,\ldots,\,a_r\,\in\,\S^2T^*$ further into a finite
 sum of signed symmetric squares. Using the distributivity law entailed
 by the bilinearity of the Nomizu--Kulkarni product $\times$ we thus arrive
 at an expansion of the algebraic curvature tensor $R$ into a finite sum of
 terms
 $$
  \pm\;\;{\textstyle\frac12}\,\alpha^2\,\times\,{\textstyle\frac12}\,\beta^2
 $$
 with linear forms $\alpha,\,\beta\,\in\,T^*$. Recalling the conventional
 equality $\frac12\alpha^2\,=\,\alpha\otimes\alpha$ we calculate
 \begin{eqnarray*}
  (\,{\textstyle\frac12}\alpha^2\times{\textstyle\frac12}\beta^2\,)
  (\,X,Y;\,U,V\,)
  &=&
  +\;\alpha(X)\,\alpha(U)\,\beta(Y)\,\beta(V)
  \;-\;\alpha(X)\,\alpha(V)\,\beta(Y)\,\beta(U)
  \\
  &&
  -\;\alpha(Y)\,\alpha(U)\,\beta(X)\,\beta(V)
  \;+\;\alpha(Y)\,\alpha(V)\,\beta(X)\,\beta(U)
  \\
  &=&
  \bigg(\,\alpha(X)\beta(Y)\,-\,\alpha(Y)\beta(X)\,\bigg)
  \bigg(\,\alpha(U)\beta(V)\,-\,\alpha(V)\beta(U)\,\bigg)
  \\
  &=&
  (\,\alpha\wedge\beta\,\otimes\,\alpha\wedge\beta\,)(\;X,\,Y;\,U,\,V\;)
 \end{eqnarray*}
 and conclude that $\frac12\alpha^2\,\times\,\frac12\beta^2\,=\,0$ for
 linearly dependent forms $\alpha,\,\beta\,\in\,T^*$. In consequence
 $$
  \pm\;\;{\textstyle\frac12}\,\alpha^2\,\times\,{\textstyle\frac12}\,\beta^2
  \;\;=\;\;
  {\textstyle\frac18}\;
  (\;\alpha^2\;\pm\;\beta^2\;)\;\times\;(\;\alpha^2\;\pm\;\beta^2\;)
 $$
 due to the commutativity of the Nomizu--Kulkarni product and
 $\alpha^2\times\alpha^2\,=\,0\,=\,\beta^2\times\beta^2$. Needless to
 say the symmetric bilinear forms $\frac1{\sqrt8}(\alpha^2+\beta^2)$ and
 $\frac1{\sqrt8}(\alpha^2-\beta^2)$ have rank $2$ and signature $(2,0)$
 and $(1,1)$ respectively unless $\alpha$ and $\beta$ are linearly
 dependent forms.

 \begin{Corollary}[Description of Algebraic Curvature Tensors]
 \hfill\label{gact}\break
  The vector space $\Curv^-T$ of algebraic curvature tensors on a finite
  dimensional euclidean vector space $T$ over $\R$ is spanned by the
  curvature tensors $\frac12\alpha^2\,\times\,\frac12\beta^2\,=\,
  \alpha\wedge\beta\,\otimes\,\alpha\wedge\beta$:  
  $$
   \Curv^-T
   \;\;=\;\;
   \span_\R\{\;\;\alpha\,\wedge\,\beta\,\otimes\,\alpha\,\wedge\,\beta\;\;|
   \;\;\alpha,\,\beta\,\in\,T^*\textrm{\ linear forms}\;\;\}
   \;\;\subseteq\;\;
   {\textstyle\bigotimes^4}T^*\ .
  $$ 
 \end{Corollary}
\section{Graph Algebras}
\label{graphs}
 Graphs are used in many different areas of mathematics as a means to encode
 information in a form easily accessible for humans. In mathematical physics
 for examples graphs or Feynman diagrams are used to encode specific analytic
 integrals by assigning Feynman rules to the different types of edges and
 vertices comprising a graph. In this section we will use trivalent graphs
 with edges colored red and black in a similar way to encode the contraction
 scheme corresponding to a stable curvature invariant. In due course we will
 construct the algebra $\A^\bullet$ of colored trivalent graphs and its
 quotient $\overline\A^\bullet$ by the IHX--relations on red edges, which
 arose originally in knot theory, to obtain a graphical calculus for
 curvature polynomials.

 \pfill
 Recall first of all that a labelling of a finite set $S$ is a bijection
 $L:\,S\stackrel\cong\longrightarrow\{1,\ldots,n\}$ with the set of the first
 $n\,=\,\#S$ natural numbers. In turn an orientation of $S$ is an equivalence
 class $o\,=\,[\,L,\,\varepsilon\,]$ of a labelling $L$ of $S$ and a sign
 $\varepsilon\,\in\,\{\,+1,\,-1\,\}$ under the equivalence relation:
 $$
  (\,L,\,\varepsilon\,)
  \;\;\sim\;\;
  (\,\hat L,\,\hat\varepsilon\,)
  \qquad\Longleftrightarrow\qquad
  \sgn(\,L\,\circ\,\hat L^{-1}\,)
  \;\;=\;\;
  \varepsilon\,\hat\varepsilon\ .
 $$
 Because all one point sets as well as the empty set $\emptyset$ have
 unique labellings $L_{\mathrm{can}}$, these sets have two distinguished
 orientations $o_\pm\,:=\,[\,L_{\mathrm{can}},\,\pm 1\,]$. In the same vein
 every labelling $L$ of a finite set $S$ represents the orientation
 $o_L\,:=\,[\,L,\,+1\,]$. Finite sets with at least two elements have
 odd permutations and hence lack distinguished orientations, nevertheless
 we can still talk about the orientation $-o\,:=\,[\,L,\,-\varepsilon\,]$
 opposite to $+o\,=\,[\,L,\,+\varepsilon\,]$.

 It is convenient for our purposes to define a finite graph as a quadruple
 $\gamma\,:=\,(\,V,\,F;\,\theta,\,\at\,)$ consisting of finite sets $V$ and
 $F$ called the sets of vertices and flags of $\gamma$ respectively, a fix
 point free involution $\theta:\,F\longrightarrow F$ called the flag involution
 and a map $\at:\,F\longrightarrow V$ called the attaching map. By assumption
 $\theta$ is fix point free with $f\,\neq\,\theta(f)$ for all $f\,\in\,F$,
 hence all orbits of $\theta$ in the set of flags have exactly two elements,
 namely the edges of the graph $\gamma$:
 $$
  \Edge\;\gamma
  \;\;:=\;\;
  F/_{\displaystyle\langle\,\theta\,\rangle}
  \;\;=\;\;
  \{\;\;\{\,f,\,\theta(f)\,\}\;\;|\;\;f\,\in\,F\;\;\}\ .
 $$
 The analogous back references $\Vert\,\gamma\,:=\,V$ and $\Flag\,\gamma
 \,:=\,F$ reduce the need to specify the vertex and flag sets of a finite
 graph $\gamma$ by name, nevertheless $\theta$ and $\at$ will always refer
 to the flag involution and the attaching map of the finite graph in question.
 By definition every flag $f\,\in\,F$ is adjacent to the vertex $\at(f)
 \,\in\,V$, similarly an edge $e\,=\,\{f_1,f_2\}$ is adjacent to the not
 necessarily different vertices $\at(f_1)$ and $\at(f_2)$. The cardinality
 of the set $\Flag_v\gamma\,:=\,\at^{-1}(\,v\,)$ of flags adjacent to a
 vertex $v\,\in\,\Vert\,\gamma$ is called its valence $\#\,\Flag_v\gamma$,
 in a bivalent and trivalent graphs respectively all vertices are required
 to have the same valence $2$ or $3$. According to the definitions above
 finite graphs may have multiple edges between vertices and loops,
 i.e.~edges from a vertex to itself.

 In order to talk about isomorphic graphs we consider finite graphs as the
 objects in a suitable category of graphs. Morphisms $\varphi:\,\gamma
 \longrightarrow\hat\gamma$ in this category are pairs of maps between
 the sets of vertices $\varphi_\Vert:\,\Vert\,\gamma\longrightarrow\Vert\,
 \hat\gamma$ and flags $\varphi_\Flag:\,\Flag\,\gamma\longrightarrow
 \Flag\,\hat\gamma$ respectively, which intertwine the flag involutions
 $\hat\theta\,\circ\,\varphi_\Flag\,=\,\varphi_\Flag\,\circ\,\theta$ and
 the attaching maps $\widehat\at\,\circ\,\varphi_\Flag\,=\,\varphi_\Vert
 \,\circ\,\at$. In consequence the automorphism group of a finite graph
 $\gamma$ comes along with a group homomorphism $\Aut\,\gamma\longrightarrow
 S_{\Vert\,\gamma},\,\varphi\longmapsto\varphi_\Vert$, whose image and kernel
 are the groups of pure and trivial automorphisms of $\gamma$ respectively:
 \begin{equation}\label{autg}
  \begin{array}{lcl}
   \overline\Aut\;\gamma
   &:=&
   \{\;\;\varphi_\Vert\;\;|\;\;\varphi\textrm{\ automorphism of $\gamma$}\;\;\}
   \\[2pt]
   \Aut_\circ\gamma
   &:=&
   \{\;\;\varphi_\Flag\;\;|\;\;\varphi
   \textrm{\ automorphism of $\gamma$ with\ }
   \varphi_\Vert\,=\,\id_{\Vert\,\gamma}\;\;\}\ .
  \end{array}
 \end{equation}
 Remarkably the induced short exact sequence $\Aut_\circ\gamma\stackrel\subset
 \longrightarrow\Aut\,\gamma\longrightarrow\overline\Aut\,\gamma$ splits on
 the right. The underlying argument is rather important for our calculations,
 because it ensures that the group $\overline\Aut\,\gamma$ of pure
 automorphisms is effectively computable as the subgroup of $S_{\Vert\,\gamma}$
 preserving the adjacency numbers of the graph $\gamma$ for all pairs
 $v,\,w\in\,\Vert\,\gamma$ of vertices
 $$
  \overline\Aut\;\gamma
  \;\;=\;\;
  \{\;\;\sigma\,\in\,S_{\Vert\,\gamma}\;\;|\;\;
  \#\,\Flag_{\sigma(v),\,\sigma(w)}\gamma\;=\;\#\,\Flag_{v,\,w}\gamma
  \textrm{\ for all\ }v,\,w\,\in\,\Vert\,\gamma\;\;\}\ ,
 $$
 where $\Flag_{v,\,w}\gamma\,:=\,\Flag_v\gamma\,\cap\,
 \theta(\,\Flag_w\gamma\,)$ is the set of flags adjacent to $v$ on
 edges to $w$.

 \begin{Definition}[Colored Trivalent Graphs]
 \hfill\label{cg}\break
  A colored trivalent graph is a trivalent graph $\gamma$ endowed
  with a coloring of its edges by colors red and black or equivalently
  by a $\theta$--invariant coloring $\col:\,\Flag\,\gamma\longrightarrow
  \{\,\red,\,\black\,\}$ of its flags in the sense $\col\,\circ\,\theta
  \,=\,\col$ such that every vertex of $\gamma$ is adjacent to exactly
  one red flag.  
 \end{Definition}

 \noindent
 The attaching map $\at$ of a colored trivalent graph $\gamma$ restricts by
 definition to a bijection $\at_\red:\,\col^{-1}(\,\red\,)\longrightarrow\Vert
 \,\gamma$ between the set of red flags of $\gamma$ and its set of vertices.
 In turn the flag involution $\theta$ induces a fix point free involution
 $\theta_\red\,:=\,\at_\red\,\circ\,\theta\,\circ\,\at_\red^{-1}$ on the
 set of vertices of the graph $\gamma$. Evidently the red edges in a
 colored trivalent graph are completely determined by the fix point free
 involution $\theta_\red$, hence we may think of a colored trivalent graph
 $\gamma$ as the bivalent graph $\gamma_\black$ we obtain by removing all red
 edges from $\gamma$ endowed with the fix point free involution $\theta_\red$
 on its set of vertices. In diagrams we will depict the red and black edges of
 a colored graph by dotted and solid lines respectively, the six graphs
 \begin{equation}\label{examp}
  \raise-17pt
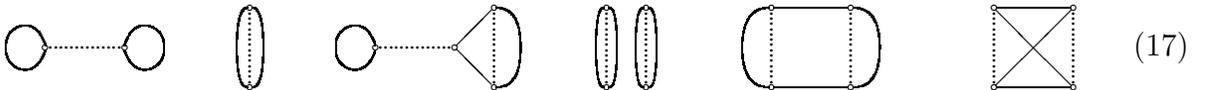
\hbox{\begin{picture}(60,40)
   \put(15,20){\circle2}
   \put(45,20){\circle2}
   \multiput(17,20)(+2, 0){14}{\circle*1}
   \qbezier( 0,20)( 0,26)( 5,28)
   \qbezier( 5,28)(12,30)(15,21)
   \qbezier(60,20)(60,26)(55,28)
   \qbezier(55,28)(48,30)(45,21)
   \qbezier( 0,20)( 0,14)( 5,12)
   \qbezier( 5,12)(12,10)(15,19)
   \qbezier(60,20)(60,14)(55,12)
   \qbezier(55,12)(48,10)(45,19)
  \end{picture}
  \qquad
  \begin{picture}(10,40)
   \put( 5, 5){\circle2}
   \put( 5,35){\circle2}
   \multiput( 5, 7)( 0,+2){14}{\circle*1}
   \qbezier( 4, 5)( 0, 5)( 0,20)
   \qbezier( 0,20)( 0,35)( 4,35)
   \qbezier( 6, 5)(10, 5)(10,20)
   \qbezier(10,20)(10,35)( 6,35)
  \end{picture}
  \qquad
  \begin{picture}(70,40)
   \put(15,20){\circle2}
   \put(45,20){\circle2}
   \multiput(17,20)(+2, 0){14}{\circle*1}
   \put(60,35){\circle2}
   \put(60, 5){\circle2}
   \multiput(60, 7)( 0,+2){14}{\circle*1}
   \put(45.7,20.7){\line(+1,+1){13.3}}
   \put(45.7,19.3){\line(+1,-1){13.3}}
   \qbezier( 0,20)( 0,26)( 5,28)
   \qbezier( 5,28)(12,30)(15,21)
   \qbezier( 0,20)( 0,14)( 5,12)
   \qbezier( 5,12)(12,10)(15,19)
   \qbezier(61, 5)(70, 5)(70,20)
   \qbezier(70,20)(70,35)(61,35)
  \end{picture}
  \qquad
  \begin{picture}(25,40)
   \put( 5, 5){\circle2}
   \put( 5,35){\circle2}
   \put(20, 5){\circle2}
   \put(20,35){\circle2}
   \multiput( 5, 7)( 0,+2){14}{\circle*1}
   \multiput(20, 7)( 0,+2){14}{\circle*1}
   \qbezier( 4, 5)( 1, 5)( 1,20)
   \qbezier( 1,20)( 1,35)( 4,35)
   \qbezier( 6, 5)( 9, 5)( 9,20)
   \qbezier( 9,20)( 9,35)( 6,35)
   \qbezier(19, 5)(16, 5)(16,20)
   \qbezier(16,20)(16,35)(19,35)
   \qbezier(21, 5)(24, 5)(24,20)
   \qbezier(24,20)(24,35)(21,35)
  \end{picture}
  \qquad
  \begin{picture}(60,40)
   \put(15, 5){\circle2}
   \put(15,35){\circle2}
   \put(45, 5){\circle2}
   \put(45,35){\circle2}
   \multiput(15, 7)( 0,+2){14}{\circle*1}
   \multiput(45, 7)( 0,+2){14}{\circle*1}
   \put(16, 5){\line(+1, 0){28}}
   \put(16,35){\line(+1, 0){28}}
   \qbezier(14, 5)( 4, 5)( 4,20)
   \qbezier( 4,20)( 4,35)(14,35)
   \qbezier(46, 5)(56, 5)(56,20)
   \qbezier(56,20)(56,35)(46,35)
  \end{picture}
  \qquad\quad
  \begin{picture}(30,40)
   \put( 0, 5){\circle2}
   \put( 0,35){\circle2}
   \put(30, 5){\circle2}
   \put(30,35){\circle2}
   \multiput( 0, 7)( 0,+2){14}{\circle*1}
   \multiput(30, 7)( 0,+2){14}{\circle*1}
   \put( 1, 5){\line(+1, 0){28}}
   \put( 1,35){\line(+1, 0){28}}
   \put( 0.7,5.7){\line(+1,+1){28.5}}
   \put(29.3,5.7){\line(-1,+1){28.5}}
  \end{picture}}
 \end{equation}
 say are among the simplest colored trivalent graphs. Homomorphisms and
 isomorphisms of colored trivalent graphs $\varphi:\,\gamma\longrightarrow
 \hat\gamma$ are required to preserve the flag coloring $\col\,=\,\hat\col\,
 \circ\,\varphi_\Flag$ of course. With this definition in place we can
 consider the sets of isomorphism classes
 $$
  \Gamma^n
  \;\;:=\;\;
  \{\;\;[\,\gamma\,]\;\;|\;\;\gamma\textrm{\ colored trivalent
   graph with $2n$ vertices}\;\;\}
 $$
 of colored trivalent graphs with $2n,\,n\,\in\,\N_0,$ vertices and their
 union $\Gamma^\bullet\,:=\,\bigcup_{n\,\in\,\N_0}\Gamma^n$, which is a graded
 commutative monoid under the multiplication induced by the disjoint union
 $\dot\cup$
 $$
  \Vert(\,\gamma\,\dot\cup\,\hat\gamma\,)
  \;\;:=\;\;
  \Vert\;\gamma\;\;\dot\cup\;\;\Vert\;\hat\gamma
  \qquad\qquad
  \Flag(\,\gamma\,\dot\cup\,\hat\gamma\,)
  \;\;:=\;\;
  \Flag\;\gamma\;\;\dot\cup\;\;\Flag\;\hat\gamma
 $$
 of colored graphs with the obvious definitions for the flag involution
 $\theta$, the attaching map $\at$ and the coloring $\col$. On the level of
 graphs the disjoint union may or may not be commutative and associative
 depending on the model of set theory in use, this minor nuisance however
 disappears on the level of isomorphism classes of colored trivalent graphs.
 Up to isomorphism for example the fourth graph in diagram (\ref{examp}) equals
 the disjoint union $\gamma\,\dot\cup\,\gamma$ of the second graph $\gamma$
 with itself. The unit element of the monoid $\Gamma^\bullet$ is represented
 by the unique empty colored trivalent graph $\gamma_{\mathrm{empty}}\,:=\,
 (\,\emptyset,\,\emptyset;\,\theta,\,\at,\,\col\,)$ without vertices or
 flags at all.

 \begin{Definition}[Algebra of Colored Trivalent Graphs]
 \hfill\label{ga}\break
  The algebra of colored trivalent graphs is the graded convolution algebra
  of the monoid $\Gamma^\bullet$ over $\R$. Put differently the algebra of
  colored trivalent graphs is the free graded vector space
  $$
   \A^\bullet
   \;\;:=\;\;
   \R\;\Gamma^\bullet
  $$
  over $\R$ generated by the graded set $\Gamma^\bullet\,=\,\bigcup\Gamma^n$
  with multiplication given by the $\R$--bilinear extension of the disjoint
  union multiplication of the commutative monoid $\Gamma^\bullet$. In passing
  we remark that the grading on $\A^\bullet$ is completely determined by the
  number operator derivation:
  $$
   N:\;\;\A^\bullet\;\longrightarrow\;\A^\bullet,\qquad[\,\gamma\,]
   \;\longmapsto\;\frac{\#\,\Vert\;\gamma}2\;[\,\gamma\,]\ .
  $$
 \end{Definition}

 \noindent
 From an algebraic point of view the algebra of colored trivalent graphs is
 not particularly interesting, because the underlying commutative monoid
 $\Gamma^\bullet$ is evidently the free commutative monoid generated by the
 subset $\Gamma^\bullet_\conn\,\subseteq\,\Gamma^\bullet$ of isomorphism
 classes of connected colored trivalent graphs. In consequence the inclusion
 $\R\Gamma^\bullet_\conn\,\subseteq\,\A^\bullet$ induces an algebra isomorphism
 $$
  \S^\bullet(\;\R\,\Gamma_\conn^\bullet\;)\;\stackrel\cong\longrightarrow\;
  \A^\bullet
 $$
 exhibiting $\A^\bullet$ as the free polynomial algebra generated by the graded
 vector space $\R\,\Gamma_\conn^\bullet$. From the differential geometric
 point of view it is much more interesting that every isomorphism class
 $[\,\gamma\,]$ of colored trivalent graphs with $2n$ vertices defines a
 homogeneous polynomial
 $$
  [\,\gamma\,]:\;\;\Curv^-T\;\longrightarrow\;\R,
  \qquad R\;\longmapsto\;[\,\gamma\,](\;R\;)\ ,
 $$
 of degree $n$ on the space of algebraic curvature tensors on an arbitrary
 euclidean vector space:

 \begin{Definition}[Stable Curvature Invariants]
 \hfill\label{slci}\break
  Every isomorphism class of colored trivalent graphs $\gamma$ with $2n$
  vertices gives rise to a homogeneous polynomial $[\,\gamma\,]:\,\Curv^-T
  \longrightarrow\R$ of degree $n\,\in\,\N_0$ for algebraic curvature tensors
  over an arbitrary euclidean vector space $T$. More precisely $\gamma$
  evaluates on an algebraic curvature tensor $R\,\in\,\Curv^-T$ to an
  iterated sum over an orthonormal basis $E_1,\,\ldots,\,E_m$ of $T$
  $$
   [\;\gamma\;](\;R\;)
   \;\;:=\;\;
   \sum_{\mu:\,\Edge\,\gamma_\black\longrightarrow\{1,\ldots,m\}}
   \prod_{{\scriptstyle\{\,v_+,\,v_-\,\}\,\subseteq\,\Vert\,\gamma}
    \atop{\scriptstyle\theta_\red(v_+)\,=\,v_-}}
   \Sec(\;E_{\overline\mu(f^1_+)},\;E_{\overline\mu(f^2_+)};
   \;E_{\overline\mu(f^1_-)},\;E_{\overline\mu(f^2_-)}\;)\ ,
  $$
  where $\Sec\,:=\,\Phi^+R$ denotes the corresponding sectional curvature
  tensor. Sum and product extend over all maps $\mu:\,\Edge\,\gamma_\black
  \longrightarrow\{1,\ldots,m\}$ and all orbits $\{v_+,v_-\}\,\subseteq\,
  \Vert\,\gamma$ of the fixed point free involution $\theta_\red$ of
  $\Vert\,\gamma$. The flags $f^1_+,\,f^2_+$ and $f^1_-,\,f^2_-$ in this
  formula denote the pairs of black flags adjacent to $v_+$ and $v_-$, while
  $\overline\mu:\,\Flag\,\gamma_\black\longrightarrow\{1,\ldots,m\}$ refers
  to the composition of $\mu$ with the canonical projection $\Flag\,
  \gamma_\black\longrightarrow\Edge\,\gamma_\black$.
 \end{Definition}

 \noindent
 In order to simplify the definition of the polynomial $[\,\gamma\,]$ we
 have omitted the usual sign factors $\pm$ appearing in iterated sums over
 orthonormal bases in cases the scalar product $g$ of the euclidean vector
 space $T$ is not positive definite. The polynomial $[\,\gamma\,]$ is of
 course well defined independent of the choices made in labelling the two
 vertices in an orbit of $\theta_\red$ and in labelling the two black flags
 adjacent to each vertex. Extending the preceding construction linearly we
 may associate a polynomial on $\Curv^-T$ to every element of the algebra
 $\A^\bullet$ of colored trivalent graphs, more precisely we obtain a
 homomorphism of graded algebras
 $$
  \Inv_{(\,T,\,g\,)}:\;\;\A^\bullet\;\longrightarrow\;
  \S^\bullet(\,\Curv^-T\,)^*,\qquad[\,\gamma\,]\;\longmapsto\;
  \Big(\;\;R\;\longmapsto\;[\,\gamma\,](\;R\;)\;\;\Big)\ ,
 $$
 because both the sum and the product split over a disjoint union of graphs to
 provide for:
 $$
  [\,\gamma\,\dot\cup\,\hat\gamma\,](\;R\;)
  \;\;=\;\;
  [\,\gamma\,](\;R\;)\;\cdot\;[\,\hat\gamma\,](\;R\;)\ .
 $$
 A discussion of the dependence of the algebra homomorphism
 $\Inv_{(\,T,\,g\,)}$ on the euclidean vector space $T$ necessarily involves
 the category $\KM$ of curvature models. Objects in this category are triples
 $(\,T,\,g,\,R\,)$ describing algebraic curvature tensors $R$ on euclidean
 vector spaces $(\,T,\,g\,)$, morphisms are the adjoints $F^*:\,T
 \longrightarrow\hat T$ of isometric maps $F:\,\hat T\longrightarrow T$
 between the euclidean vector spaces $T$ and $\hat T$ satisfying either
 of the two equivalent constraints
 \begin{equation}\label{sfc}
  \begin{array}{lcl}
   \;\widehat R(\;F^*X,\,F^*Y;\,F^*U,\,F^*V\;)
   &=&
   \;R(\;X,\,Y;\,U,\,V\;)
   \\[2pt]
   \widehat\Sec(\;F^*X,\,F^*Y;\,F^*U,\,F^*V\;)
   &=&
   \Sec(\;X,\,Y;\,U,\,V\;)
  \end{array}
 \end{equation}
 for all $X,\,Y,\,U,\,V\,\in\,T$. For every such morphism $(\,T,\,g,\,R\,)
 \longrightarrow(\,\hat T,\,\hat g,\,\hat R\,)$ we have equality
 \begin{equation}\label{inv}
  [\,\gamma\,](\;R\;)
  \;\;=\;\;
  [\,\gamma\,](\;\hat R\;)\ ,
 \end{equation}
 because the adjoint $F^*:\,T\longrightarrow\hat T$ of an isometric map
 $F:\,\hat T\longrightarrow T$ factors into the orthogonal projection to the
 regular subspace $\im\;F\,\subseteq\,T$ and the isometry $F^{-1}:\,\im\;F
 \longrightarrow\hat T$. In consequence we can choose an orthonormal basis
 $E_1,\,\ldots,\,E_m$ for the euclidean vector space $T$ in such a way that
 $F^*E_1,\,\ldots,\,F^*E_{\hat m}$ is an orthonormal basis for $\hat T$,
 while $F^*E_\mu\,=\,0$ for all $\mu\,>\,\hat m$. Due to the constraint
 (\ref{sfc}) imposed on morphisms in the category $\KM$ we find
 $$
  \widehat\Sec(\;F^*E_{\overline\mu(f^1_+)},\;F^*E_{\overline\mu(f^2_+)};
  \;F^*E_{\overline\mu(f^1_-)},\;F^*E_{\overline\mu(f^2_-)}\;)
  \;\;=\;\;  
  \Sec(\;E_{\overline\mu(f^1_+)},\;E_{\overline\mu(f^2_+)};
  \;E_{\overline\mu(f^1_-)},\;E_{\overline\mu(f^2_-)}\;)
 $$
 for this particular choice of orthonormal bases and all $\mu:\,\Edge\,
 \gamma_\black\longrightarrow\{1,\ldots,\hat m\}$, while
 $$
  \prod_{{\scriptstyle\{\,v_+,\,v_-\,\}\,\subseteq\,\Vert\,\gamma}
   \atop{\scriptstyle\theta_\red(v_+)\,=\,v_-}}
  \Sec(\;E_{\overline\mu(f^1_+)},\;E_{\overline\mu(f^2_+)};
  \;E_{\overline\mu(f^1_-)},\;E_{\overline\mu(f^2_-)}\;)
  \;\;=\;\;
  0
 $$
 for all $\mu:\,\Edge\,\gamma_\black\longrightarrow\{1,\ldots,m\}$ with maximum
 larger than $\hat m$. A particular case of the invariance (\ref{inv}) of the
 polynomial $[\,\gamma\,]$ under morphisms in the category $\KM$ occurs for
 the adjoints $F^*\,=\,F^{-1}$ of self isometries $F:\,T\longrightarrow T$
 of a given euclidean vector space $T$. In this case the invariance (\ref{inv})
 reads $[\,\gamma\,](\,F\star\hat R\,)\,=\,[\,\gamma\,](\,\hat R\,)$ for all
 algebraic curvature tensors $\hat R\,\in\,\Curv^-T$ in terms of the natural
 representation $\star$ of the orthogonal group $\OG(\,T,\,g\,)$ of isometries
 of $T$ on $\Curv^-T$; in other words the homomorphism of graded algebras
 $$
  \Inv_{(\,T,\,g\,)}:\;\;\A^\bullet\;\longrightarrow\;
  [\;\S^\bullet(\,\Curv^-T\,)^*\,]^{\OG(\,T,\,g\,)}
 $$
 take values in the subalgebra of $\OG(\,T,\,g\,)$--invariant polynomials on
 $\Curv^-T$. A pretty similar argument to the one used to prove the invariance
 (\ref{inv}) implies that the polynomial associated to a linear combination
 $[\,\gamma\,]\,\in\,\R\Gamma^\bullet_\conn$ of connected colored trivalent
 graphs is additive
 $$
  [\,\gamma\,](\;R\,\oplus\,\hat R\;)
  \;\;=\;\;
  [\,\gamma\,](\;R\;)\;+\;[\,\gamma\,](\;\hat R\;)
 $$
 under Cartesian products of curvature models. For general elements of
 the algebra $\A^\bullet$ of colored trivalent graphs we can thus use the
 comultiplication $\Delta:\,\A^\bullet\longrightarrow\A^\bullet\,\otimes\,
 \A^\bullet$ induced by the algebra isomorphism $\A^\bullet\,\cong
 \,\S^\bullet(\,\R\Gamma^\bullet_\conn\,)$ to find $[\,\gamma\,]
 (\,R\oplus\hat R\,)\,=\,[\,\Delta\gamma\,](\,R,\,\hat R\,)$.

 \pfill
 Perhaps the best way to think about Definition \ref{slci} is as a set of
 Feynman rules for some unspecified field theory, which allow us to evaluate
 the isomorphism class $[\,\gamma\,]$ of a colored trivalent graph $\gamma$
 on the sectional curvature tensor $\Sec\,:=\,\Phi^+R$ associated to an
 algebraic curvature tensor $R\,\in\,\Curv^-T$. These Feynman rules stipulate
 a summation over an orthonormal basis for every black edge and the
 multiplication of the red edge interactions:
 \begin{equation}\label{feyn1}
  \raise-19pt\hbox{\begin{picture}(54,45)
   \put(28,13){\circle2}
   \put(28,33){\circle2}
   \put(10, 4){\line(+2,+1){17}}
   \put(10,42){\line(+2,-1){17}}
   \put(46, 4){\line(-2,+1){17}}
   \put(46,42){\line(-2,-1){17}}
   \multiput(28,15)( 0,+2){9}{\circle*1}
   \put( 0,39){$X$}
   \put(48,39){$Y$}
   \put( 0, 0){$V$}
   \put(48, 0){$U$}
  \end{picture}}
  \;\;\widehat=\;\;
  \Sec(\;X,\;Y;\;U,\;V\;)\ .
 \end{equation}
 In light of this Feynman rule description of Definition \ref{slci} we
 readily observe that the algebra homomorphism $\Inv_{(\,T,\,g\,)}$ can not
 be injective, because the Feynman interpretation of
 \begin{equation}\label{ihx}
  \raise-17pt\hbox{\begin{picture}(40,40)
   \put(20,10){\circle2}
   \put(20,30){\circle2}
   \multiput(20,12)( 0,+2){9}{\circle*1}
   \put( 1, 1){\line(+2,+1){18}}
   \put(39, 1){\line(-2,+1){18}}
   \put( 1,39){\line(+2,-1){18}}
   \put(39,39){\line(-2,-1){18}}
  \end{picture}}
  \;\;\;+\;\;\;
  \raise-17pt\hbox{\begin{picture}(40,40)
   \put(10,20){\circle2}
   \put(30,20){\circle2}
   \multiput(12,20)(+2, 0){9}{\circle*1}
   \put( 1, 1){\line(+1,+2){9}}
   \put(39, 1){\line(-1,+2){9}}
   \put( 1,39){\line(+1,-2){9}}
   \put(39,39){\line(-1,-2){9}}
  \end{picture}}
  \;\;\;+\;\;\;
  \raise-17pt\hbox{\begin{picture}(40,40)
   \put(10,20){\circle2}
   \put(30,20){\circle2}
   \multiput(12,20)(+2, 0){9}{\circle*1}
   \put( 1, 1){\line(+1,+2){9}}
   \put(39, 1){\line(-1,+2){9}}
   \put( 1,39){\line(+3,-2){28}}
   \put(39,39){\line(-3,-2){28}}
  \end{picture}}
 \end{equation}
 is exactly the left hand side of the first Bianchi identity (\ref{1BI'}) for
 sectional curvature tensors. Whenever colored trivalent graphs $\gamma_I,\,
 \gamma_H$ and $\gamma_X$ are isomorphic except for a pair of vertices
 connected by a red edge and four black flags attaching to these two
 vertices according to
 \begin{equation}\label{ihxex}
  \gamma_I:\quad
  \raise-17pt\hbox{\begin{picture}(40,40)
   \put(20,10){\circle2}
   \put(20,30){\circle2}
   \multiput(20,12)( 0,+2){9}{\circle*1}
   \put( 1, 1){\line(+2,+1){18}}
   \put(39, 1){\line(-2,+1){18}}
   \put( 1,39){\line(+2,-1){18}}
   \put(39,39){\line(-2,-1){18}}
  \end{picture}}
  \qquad\qquad\qquad
  \gamma_H:\quad
  \raise-17pt\hbox{\begin{picture}(40,40)
   \put(10,20){\circle2}
   \put(30,20){\circle2}
   \multiput(12,20)(+2, 0){9}{\circle*1}
   \put( 1, 1){\line(+1,+2){9}}
   \put(39, 1){\line(-1,+2){9}}
   \put( 1,39){\line(+1,-2){9}}
   \put(39,39){\line(-1,-2){9}}
  \end{picture}}
  \qquad\qquad\qquad
  \gamma_X:\quad
  \raise-17pt\hbox{\begin{picture}(40,40)
   \put(10,20){\circle2}
   \put(30,20){\circle2}
   \multiput(12,20)(+2, 0){9}{\circle*1}
   \put( 1, 1){\line(+1,+2){9}}
   \put(39, 1){\line(-1,+2){9}}
   \put( 1,39){\line(+3,-2){28}}
   \put(39,39){\line(-3,-2){28}}
  \end{picture}}\ ,
 \end{equation}
 then $\Inv_{(\,T,\,g\,)}([\,\gamma_I\,]+[\,\gamma_H\,]+[\,\gamma_X\,])\,=\,0$.
 The subspace spanned by these IHX--relations
 \begin{equation}\label{span}
  \langle\;\IHX\;\rangle
  \;\;:=\;\;
  \span_\R\{\;\;[\,\gamma_I\,]\,+\,[\,\gamma_H\,]\,+\,[\,\gamma_X\,]
  \;\;|\;\;\gamma_I,\,\gamma_H,\,\gamma_X\textrm{\ isomorphic except for
  (\ref{ihxex})}\;\;\}
 \end{equation}
 is by construction a homogeneous ideal in the graph algebra $\A^\bullet$,
 because an IHX--relation for some red edge remains an IHX--relation at the
 same red edge after taking the disjoint union product with another graph. In
 consequence the algebra homomorphism $\Inv_{(\,T,\,g\,)}$ factorizes through
 the canonical projection to the reduced algebra of colored trivalent graphs
 \begin{equation}\label{rga}
  \overline\A^\bullet
  \;\;:=\;\;
  \A^\bullet/_{\displaystyle\langle\;\IHX\;\rangle}
 \end{equation}
 and an algebra homomorphism $\overline\Inv_{(\,T,\,g\,)}:\,\overline
 \A^\bullet\longrightarrow[\,\S^\bullet(\,\Curv^-T\,)^*\,]^{\OG(\,T,\,g\,)}$.
 One of the simplest examples of a congruence modulo the ideal $\langle\,
 \IHX\,\rangle$ is given by the IHX--relation:
 \begin{equation}\label{ihq}
  \raise-9pt\hbox{\begin{picture}(24,24)
   \put( 1, 1){\circle2}
   \put( 1,23){\circle2}
   \put(23, 1){\circle2}
   \put(23,23){\circle2}
   \multiput( 1, 3)( 0,+2){10}{\circle*1}
   \multiput(23, 3)( 0,+2){10}{\circle*1}
   \put(1.8, 0.5){\line(+1, 0){20.4}}
   \put(1.8, 1.5){\line(+1, 0){20.4}}
   \put(1.8,22.5){\line(+1, 0){20.4}}
   \put(1.8,23.5){\line(+1, 0){20.4}}
  \end{picture}}
  \;+\;2\;
  \raise-9pt\hbox{\begin{picture}(24,24)
   \put( 1, 1){\circle2}
   \put( 1,23){\circle2}
   \put(23, 1){\circle2}
   \put(23,23){\circle2}
   \multiput( 1, 3)( 0,+2){10}{\circle*1}
   \multiput(23, 3)( 0,+2){10}{\circle*1}
   \put( 2, 1){\line(+1, 0){20}}
   \put( 2,23){\line(+1, 0){20}}
   \put( 1.7,1.7){\line(+1,+1){20.4}}
   \put(22.3,1.7){\line(-1,+1){20.4}}
  \end{picture}}
  \;\;=\;\;
  \raise-12pt\hbox{\begin{picture}(40,28)
   \put(39, 7){\circle2}
   \put(39,21){\circle2}
   \multiput(39, 9)( 0,+2){6}{\circle*1}
   \put( 1, 1){\circle*2}
   \put(27, 1){\circle*2}
   \put( 1,27){\circle*2}
   \put(27,27){\circle*2}
   \qbezier(39.0, 5.9)(34.0,-12.0)( 1.0, 1.0)
   \qbezier(38.2, 6.4)(32.0,- 2.0)(27.0, 1.0)
   \qbezier(39.0,22.1)(34.0, 42.0)( 1.0,27.0)
   \qbezier(38.2,21.6)(32.0, 30.0)(27.0,27.0)
   \put(14, 7){\circle2}
   \put(14,21){\circle2}
   \multiput(14, 9)( 0,+2){6}{\circle*1}
   \put( 1, 1){\line(+2,+1){11.8}}
   \put(27, 1){\line(-2,+1){11.8}}
   \put( 1,27){\line(+2,-1){11.8}}
   \put(27,27){\line(-2,-1){11.8}}
  \end{picture}}
  \;+\;
  \raise-12pt\hbox{\begin{picture}(40,28)
   \put(39, 7){\circle2}
   \put(39,21){\circle2}
   \multiput(39, 9)( 0,+2){6}{\circle*1}
   \put( 1, 1){\circle*2}
   \put(27, 1){\circle*2}
   \put( 1,27){\circle*2}
   \put(27,27){\circle*2}
   \qbezier(39.0, 5.9)(34.0,-12.0)( 1.0, 1.0)
   \qbezier(38.2, 6.4)(32.0,- 2.0)(27.0, 1.0)
   \qbezier(39.0,22.1)(34.0, 42.0)( 1.0,27.0)
   \qbezier(38.2,21.6)(32.0, 30.0)(27.0,27.0)
   \put( 7,14){\circle2}
   \put(21,14){\circle2}
   \multiput( 9,14)(+2, 0){6}{\circle*1}
   \put( 1, 1){\line(+1,+2){5.9}}
   \put(27, 1){\line(-1,+2){5.9}}
   \put( 1,27){\line(+1,-2){5.9}}
   \put(27,27){\line(-1,-2){5.9}}
  \end{picture}}
  \;+\;
  \raise-12pt\hbox{\begin{picture}(40,28)
   \put(39, 7){\circle2}
   \put(39,21){\circle2}
   \multiput(39, 9)( 0,+2){6}{\circle*1}
   \put( 1, 1){\circle*2}
   \put(27, 1){\circle*2}
   \put( 1,27){\circle*2}
   \put(27,27){\circle*2}
   \qbezier(39.0, 5.9)(34.0,-12.0)( 1.0, 1.0)
   \qbezier(38.2, 6.4)(32.0,- 2.0)(27.0, 1.0)
   \qbezier(39.0,22.1)(34.0, 42.0)( 1.0,27.0)
   \qbezier(38.2,21.6)(32.0, 30.0)(27.0,27.0)
   \put( 7,14){\circle2}
   \put(21,14){\circle2}
   \multiput( 9,14)(+2, 0){6}{\circle*1}
   \put( 1, 1){\line(+1,+2){ 5.9}}
   \put(27, 1){\line(-1,+2){ 5.9}}
   \put( 1,27){\line(+3,-2){18.9}}
   \put(27,27){\line(-3,-2){18.9}}
  \end{picture}}
  \;\;\equiv\;\;
  0\ .
 \end{equation}
 Another useful family of congruences modulo $\langle\,\IHX\,\rangle$ is
 described by the Feynman rules:
 \begin{equation}\label{ihex}
  \raise-12pt\hbox{\begin{picture}(64,30)
   \put(20,15){\circle2}
   \put(40,15){\circle2}
   \put(41,16){\line(+1,+1){10}}\put(53,21){$X$}
   \put(41,14){\line(+1,-1){10}}\put(53, 1){$Y$}
   \multiput(22,15)(+2, 0){9}{\circle*1}
   \qbezier( 4,15)( 4,21)(10,22)
   \qbezier(10,22)(15,23)(19,16)
   \qbezier( 4,15)( 4, 9)(10, 8)
   \qbezier(10, 8)(15, 7)(19,14)
  \end{picture}}
  \;\;\equiv\;\;
  -\;2\,
  \raise-12pt\hbox{\begin{picture}(40,30)
   \put(10,25){\circle2}
   \put(10, 5){\circle2}
   \put(11,25){\line(+1, 0){14}}\put(27,21){$X$}
   \put(11, 5){\line(+1, 0){14}}\put(27, 1){$Y$}
   \multiput(10, 7)(0,+2){9}{\circle*1}
   \qbezier( 9,25)( 2,24)(2,15)
   \qbezier( 9, 5)( 2, 6)(2,15)
  \end{picture}}
  \;\;\widehat=\;\;
  4\,\Ric(\,X,\,Y\,)\ .
 \end{equation}
 Concerning the structure of the reduced graph algebra we remark
 $$
  \A^\bullet
  \;\;\cong\;\;
  \S^\bullet(\,\R\Gamma^\bullet_\conn\,)
  \qquad\Longrightarrow\qquad
  \overline\A^\bullet
  \;\;\cong\;\;
  \S^\bullet\Big(\,\R\Gamma^\bullet_\conn/_{\displaystyle\R\Gamma^\bullet_\conn
  \,\cap\,\langle\;\IHX\;\rangle}\;\Big)
 $$
 by decomposing a colored trivalent graph $\gamma$ into its connected
 components as before. In order to establish a much stronger result we will
 make use of the bivalent black subgraph $\gamma_\black$ obtained from a
 colored trivalent graph $\gamma$ by removing all its red edges. Consider
 for the moment an arbitrary red edge in a connected colored trivalent
 graph $\gamma$:
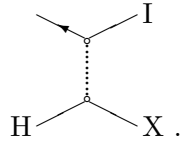
 $$
  \begin{picture}(40,44)
   \put(20, 9){\circle2}
   \put(20,31){\circle2}
   \multiput(20,11)( 0,+2){10}{\circle*1}
   \put(21, 8){\line(+2,-1){18}}
   \put(19, 8){\line(-2,-1){18}}
   \put(21,32){\line(+2,+1){18}}
   \put(19,32){\line(-2,+1){18}}
   \put( 9,37){\vector(-2,+1){0}}
   \put(-9,-5){\small{H}}
   \put(41,-5){\small{X}\ .}
   \put(41,38){\small{I}}
  \end{picture}
 $$
 After leaving this red edge along the top left flag we will necessarily
 return to this red edge at some point or other along a different black
 flag, because every finite bivalent graph like $\gamma_\black$ equals the
 disjoint union of cycles of length $\geq1$. Depending on the black flag of
 first return exactly one of the three configurations in the IHX--relation
 corresponding to the chosen red edge will have one black cycle more than
 the other two configurations. Say the H configuration will have one cycle
 more than both the I and X configurations, if we return first along the
 black flag marked with H, analogous considerations apply for a first
 return along the black flag marked with I or X.

 Given a connected graph $\gamma$ such that $\gamma_\black$ has at least two
 cycles there necessarily exists a red edge connecting two different black
 cycles. In turn the IHX--relation for such a red edge becomes a congruence
 $[\,\gamma\,]\,\equiv\,-([\,\hat\gamma\,]+[\,\tilde\gamma\,])$ modulo the
 ideal $\langle\,\IHX\,\rangle$ with connected graphs $\hat\gamma$ and
 $\tilde\gamma$ such that $\hat\gamma_\black$ and $\tilde\gamma_\black$
 both have one cycle less than $\gamma_\black$. Repeating this process
 with $\hat\gamma$ and $\tilde\gamma$ we eventually end up with a congruence
 $[\,\gamma\,]\,\equiv\,\pm([\,\hat\gamma_1\,]+\ldots+[\,\hat\gamma_r\,])$
 modulo the ideal $\langle\,\IHX\,\rangle$, in which all graphs
 $\hat\gamma_1,\,\ldots,\,\hat\gamma_r$ have connected black subgraphs:
 
 \begin{Corollary}[Generators for the Reduced Graph Algebra]
 \hfill\label{gga}\break
  The reduced graph algebra $\overline\A^\bullet$, the quotient of the
  algebra $\A^\bullet$ of colored trivalent graphs modulo the ideal
  $\langle\,\IHX\,\rangle$ of IHX--relations, is generated as
  an algebra, albeit not freely generated, by the classes $[\,\gamma\,]$
  of colored trivalent graphs $\gamma$ with connected bivalent black
  subgraph $\gamma_\black$:
  $$
   \overline\A^\bullet
   \;\;:=\;\;
   \A^\bullet/_{\displaystyle\langle\;\IHX\;\rangle}
   \;\;=\;\;
   \langle\,\{\;[\,\gamma\,]\,+\,\langle\,\IHX\,\rangle\;|\;
   \textrm{$\gamma$ colored trivalent graph, $\gamma_\black$ connected}
   \;\}\,\rangle\ .
  $$
 \end{Corollary}

 \noindent
 In consequence the list of generators of the reduced graph algebra
 $\overline\A^\bullet$ up to degree $3$ reads:
 \begin{equation}\label{gen6}
  \raise-9pt
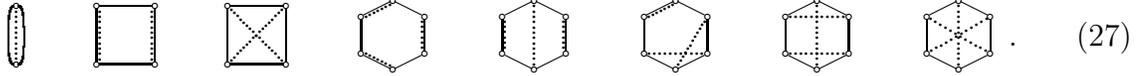
\hbox{\begin{picture}(10,26)
   \put( 5, 2){\circle2}
   \put( 5,24){\circle2}
   \qbezier(4.0, 2.0)( 2.0, 2.0)( 2.0,12.0)
   \qbezier(4.0,24.0)( 2.0,24.0)( 2.0,12.0)
   \qbezier(6.1, 2.0)( 8.0, 2.0)( 8.0,12.0)
   \qbezier(6.1,24.0)( 8.0,24.0)( 8.0,12.0)
   \multiput( 5.0, 4.0)( 0.0,+2.0){10}{\circle*1}
  \end{picture}}
  \qquad
  \raise-9pt\hbox{\begin{picture}(26,26)
   \put( 2, 2){\circle2}
   \put(24, 2){\circle2}
   \put( 2,24){\circle2}
   \put(24,24){\circle2}
   \put( 2.0, 3.0){\line( 0,+1){20}}
   \put(24.0, 3.0){\line( 0,+1){20}}
   \put( 3.0, 2.0){\line(+1, 0){20}}
   \put( 3.0,24.0){\line(+1, 0){20}}
   \multiput( 3.2, 3.9)( 0.0,+2.0){10}{\circle*1}
   \multiput(22.8, 3.9)( 0.0,+2.0){10}{\circle*1}
  \end{picture}}
  \qquad
  \raise-9pt\hbox{\begin{picture}(26,26)
   \put( 2, 2){\circle2}
   \put(24, 2){\circle2}
   \put( 2,24){\circle2}
   \put(24,24){\circle2}
   \put( 2.0, 3.0){\line( 0,+1){20}}
   \put(24.0, 3.0){\line( 0,+1){20}}
   \put( 3.0, 2.0){\line(+1, 0){20}}
   \put( 3.0,24.0){\line(+1, 0){20}}
   \multiput( 3.4, 3.4)(+1.2,+1.2){17}{\circle*1}
   \multiput(22.6, 3.4)(-1.2,+1.2){17}{\circle*1}
  \end{picture}}
  \qquad
  \raise-11pt\hbox{\begin{picture}(30,30)
   \put(15, 2){\circle2}
   \put(15,28){\circle2}
   \put( 3, 8){\circle2}
   \put(27, 8){\circle2}
   \put( 3,22){\circle2}
   \put(27,22){\circle2}
   \put( 3.0, 9.0){\line( 0,+1){12}}
   \put(27.0, 9.0){\line( 0,+1){12}}
   \put(15.9, 2.5){\line(+2,+1){10}}
   \put(14.1, 2.5){\line(-2,+1){10}}
   \put(15.9,27.5){\line(+2,-1){10}}
   \put(14.1,27.5){\line(-2,-1){10}}
   \multiput( 4.9,21.9)(+1.5,+0.75){7}{\circle*1}
   \multiput( 4.9, 8.1)(+1.5,-0.75){7}{\circle*1}
   \multiput(26.1,10.1)( 0.0,+2.0){6}{\circle*1}
  \end{picture}}
  \qquad
  \raise-11pt\hbox{\begin{picture}(30,30)
   \put(15, 2){\circle2}
   \put(15,28){\circle2}
   \put( 3, 8){\circle2}
   \put(27, 8){\circle2}
   \put( 3,22){\circle2}
   \put(27,22){\circle2}
   \put( 3.0, 9.0){\line( 0,+1){12}}
   \put(27.0, 9.0){\line( 0,+1){12}}
   \put(15.9, 2.5){\line(+2,+1){10}}
   \put(14.1, 2.5){\line(-2,+1){10}}
   \put(15.9,27.5){\line(+2,-1){10}}
   \put(14.1,27.5){\line(-2,-1){10}}
   \multiput( 3.9,10.1)( 0.0,+2.0){6}{\circle*1}
   \multiput(15.0, 4.0)( 0.0,+2.0){12}{\circle*1}
   \multiput(26.1,10.1)( 0.0,+2.0){6}{\circle*1}
  \end{picture}}
  \qquad
  \raise-11pt\hbox{\begin{picture}(30,30)
   \put(15, 2){\circle2}
   \put(15,28){\circle2}
   \put( 3, 8){\circle2}
   \put(27, 8){\circle2}
   \put( 3,22){\circle2}
   \put(27,22){\circle2}
   \put( 3.0, 9.0){\line( 0,+1){12}}
   \put(27.0, 9.0){\line( 0,+1){12}}
   \put(15.9, 2.5){\line(+2,+1){10}}
   \put(14.1, 2.5){\line(-2,+1){10}}
   \put(15.9,27.5){\line(+2,-1){10}}
   \put(14.1,27.5){\line(-2,-1){10}}
   \multiput( 4.9,21.9)(+1.5,+0.75){7}{\circle*1}
   \multiput(16.1, 3.8)(+0.9,+1.5){12}{\circle*1}
   \multiput( 5.0, 8.0)(+2.0, 0.0){11}{\circle*1}
  \end{picture}}
  \qquad
  \raise-11pt\hbox{\begin{picture}(30,30)
   \put(15, 2){\circle2}
   \put(15,28){\circle2}
   \put( 3, 8){\circle2}
   \put(27, 8){\circle2}
   \put( 3,22){\circle2}
   \put(27,22){\circle2}
   \put( 3.0, 9.0){\line( 0,+1){12}}
   \put(27.0, 9.0){\line( 0,+1){12}}
   \put(15.9, 2.5){\line(+2,+1){10}}
   \put(14.1, 2.5){\line(-2,+1){10}}
   \put(15.9,27.5){\line(+2,-1){10}}
   \put(14.1,27.5){\line(-2,-1){10}}
   \multiput(15.0, 4.0)( 0.0,+2.0){12}{\circle*1}
   \multiput( 4.9, 8.2)(+2.0, 0.0){11}{\circle*1}
   \multiput( 4.9,21.8)(+2.0, 0.0){11}{\circle*1}
  \end{picture}}
  \qquad
  \raise-11pt\hbox{\begin{picture}(30,30)
   \put(15, 2){\circle2}
   \put(15,28){\circle2}
   \put( 3, 8){\circle2}
   \put(27, 8){\circle2}
   \put( 3,22){\circle2}
   \put(27,22){\circle2}
   \put( 3.0, 9.0){\line( 0,+1){12}}
   \put(27.0, 9.0){\line( 0,+1){12}}
   \put(15.9, 2.5){\line(+2,+1){10}}
   \put(14.1, 2.5){\line(-2,+1){10}}
   \put(15.9,27.5){\line(+2,-1){10}}
   \put(14.1,27.5){\line(-2,-1){10}}
   \multiput(15.0, 4.0)( 0.0,+2.0){12}{\circle*1}
   \multiput( 4.7, 9.0)(+1.9,+1.1){12}{\circle*1}
   \multiput(25.3, 9.0)(-1.9,+1.1){12}{\circle*1}
  \end{picture}}\ .
 \end{equation}
 Interestingly there exist no further IHX--relations between these generators
 with connected black subgraphs up to degree $3$. In degree $4$ however there
 exist exactly two independent relations between the $17$ isomorphism classes
 of colored trivalent graphs with $8$ vertices and connected black subgraph in
 the reduced graph algebra $\overline\A^\bullet$ induced by the congruences
 \begin{eqnarray*}
  \lefteqn{\raise-17pt\hbox{\begin{picture}(40,40)(0,0)
   \put(11.7, 0.0){\circle2}
   \put(28.3, 0.0){\circle2}
   \put( 0.0,11.7){\circle2}
   \put( 0.0,28.3){\circle2}
   \put(11.7,40.0){\circle2}
   \put(28.3,40.0){\circle2}
   \put(40.0,11.7){\circle2}
   \put(40.0,28.3){\circle2}
   \put(12.7, 0.0){\line(+1, 0){14.6}}
   \put(12.7,40.0){\line(+1, 0){14.6}}
   \put( 0.0,12.7){\line( 0,+1){14.6}}
   \put(40.0,12.7){\line( 0,+1){14.6}}
   \put(10.9, 0.6){\line(-1,+1){10.2}}
   \put(29.1, 0.6){\line(+1,+1){10.2}}
   \put(10.9,39.4){\line(-1,-1){10.2}}
   \put(29.1,39.4){\line(+1,-1){10.2}}
   \multiput(13.4, 1.0)(+1.9, 0.0){8}{\circle*1}
   \multiput(11.1,37.7)(-0.6,-1.5){17}{\circle*1}
   \multiput(28.9,37.7)(+0.6,-1.5){17}{\circle*1}
   \multiput( 1.8,28.3)(+1.9, 0.0){20}{\circle*1}
  \end{picture}}
  \;-\;
  \raise-17pt\hbox{\begin{picture}(40,40)(0,0)
   \put(11.7, 0.0){\circle2}
   \put(28.3, 0.0){\circle2}
   \put( 0.0,11.7){\circle2}
   \put( 0.0,28.3){\circle2}
   \put(11.7,40.0){\circle2}
   \put(28.3,40.0){\circle2}
   \put(40.0,11.7){\circle2}
   \put(40.0,28.3){\circle2}
   \put(12.7, 0.0){\line(+1, 0){14.6}}
   \put(12.7,40.0){\line(+1, 0){14.6}}
   \put( 0.0,12.7){\line( 0,+1){14.6}}
   \put(40.0,12.7){\line( 0,+1){14.6}}
   \put(10.9, 0.6){\line(-1,+1){10.2}}
   \put(29.1, 0.6){\line(+1,+1){10.2}}
   \put(10.9,39.4){\line(-1,-1){10.2}}
   \put(29.1,39.4){\line(+1,-1){10.2}}
   \multiput(13.4, 1.0)(+1.9,+0.0){8}{\circle*1}
   \multiput(13.1,38.1)(+1.2,-1.2){22}{\circle*1}
   \multiput( 2.2,12.5)(+1.7,+0.7){22}{\circle*1}
   \multiput( 2.2,28.5)(+1.6,+0.7){16}{\circle*1}
  \end{picture}}}
  \qquad
  &&
  \\[4pt]
  &=&
  \raise-17pt\hbox{\begin{picture}(40,40)(0,0)
   \put(20.0,40.0){\circle2}
   \put(33.0, 0.0){\circle2}
   \put( 7.0, 0.0){\circle2}
   \put(41.8,25.5){\circle2} 
   \put(-1.8,25.5){\circle2}
   \put(20.0,29.0){\circle2}
   \put(12.2,16.5){\circle2}
   \put(27.8,16.5){\circle2}
   \put( 8.0, 0.0){\line(+1, 0){24.0}}
   \put( 6.5, 0.8){\line(-1,+3){ 7.9}}
   \put(33.5, 0.8){\line(+1,+3){ 7.9}}
   \put(19.2,39.5){\line(-3,-2){20.0}}
   \put(20.8,39.5){\line(+3,-2){20.0}}
   \put(19.4,28.2){\line(-2,-3){ 7.0}}
   \put(20.6,28.2){\line(+2,-3){ 7.0}}
   \put(13.2,16.5){\line(+1, 0){13.6}}
   \multiput( 8.7, 1.0)(+1.9,+0.0){13}{\circle*1}
   \multiput(20.0,30.9)( 0.0,+1.8){5}{\circle*1}
   \multiput(10.5,17.7)(-1.5,+1.0){8}{\circle*1}
   \multiput(29.5,17.7)(+1.5,+1.0){8}{\circle*1}
   \put( 9.4, 9.7){$\scriptscriptstyle A$}
   \put(24.6, 9.7){$\scriptscriptstyle B$}
   \put(21.0,32.4){$\scriptscriptstyle*$}
  \end{picture}}
  \;\;+\;
  \raise-17pt\hbox{\begin{picture}(40,40)(0,0)
   \put(11.7, 0.0){\circle2}
   \put(28.3, 0.0){\circle2}
   \put( 0.0,11.7){\circle2}
   \put( 0.0,28.3){\circle2}
   \put(11.7,40.0){\circle2}
   \put(28.3,40.0){\circle2}
   \put(40.0,11.7){\circle2}
   \put(40.0,28.3){\circle2}
   \put(12.7, 0.0){\line(+1, 0){14.6}}
   \put(12.7,40.0){\line(+1, 0){14.6}}
   \put( 0.0,12.7){\line( 0,+1){14.6}}
   \put(40.0,12.7){\line( 0,+1){14.6}}
   \put(10.9, 0.6){\line(-1,+1){10.2}}
   \put(29.1, 0.6){\line(+1,+1){10.2}}
   \put(10.9,39.4){\line(-1,-1){10.2}}
   \put(29.1,39.4){\line(+1,-1){10.2}}
   \multiput(13.4, 1.0)(+1.9, 0.0){8}{\circle*1}
   \multiput(11.1,37.7)(-0.6,-1.5){17}{\circle*1}
   \multiput(28.9,37.7)(+0.6,-1.5){17}{\circle*1}
   \multiput( 1.9,28.3)(+1.9, 0.0){20}{\circle*1}
   \put( 8.7,42.4){$\scriptscriptstyle A$}
   \put(26.0,42.4){$\scriptscriptstyle B$}
  \end{picture}}
  \;\;+\;
  \raise-17pt\hbox{\begin{picture}(40,40)(0,0)
   \put(11.7, 0.0){\circle2}
   \put(28.3, 0.0){\circle2}
   \put( 0.0,11.7){\circle2}
   \put( 0.0,28.3){\circle2}
   \put(11.7,40.0){\circle2}
   \put(28.3,40.0){\circle2}
   \put(40.0,11.7){\circle2}
   \put(40.0,28.3){\circle2}
   \put(12.7, 0.0){\line(+1, 0){14.6}}
   \put(12.7,40.0){\line(+1, 0){14.6}}
   \put( 0.0,12.7){\line( 0,+1){14.6}}
   \put(40.0,12.7){\line( 0,+1){14.6}}
   \put(10.9, 0.6){\line(-1,+1){10.2}}
   \put(29.1, 0.6){\line(+1,+1){10.2}}
   \put(10.9,39.4){\line(-1,-1){10.2}}
   \put(29.1,39.4){\line(+1,-1){10.2}}
   \multiput(13.4, 1.0)(+1.9, 0.0){8}{\circle*1}
   \multiput(12.9,38.1)(+1.2,-1.2){22}{\circle*1}
   \multiput(27.1,38.1)(-1.2,-1.2){22}{\circle*1}
   \multiput( 1.9,28.3)(+1.9, 0.0){20}{\circle*1}
   \put( 8.7,42.4){$\scriptscriptstyle B$}
   \put(26.0,42.4){$\scriptscriptstyle A$}
  \end{picture}}
  \;-\;\; 
  \raise-17pt\hbox{\begin{picture}(40,40)(0,0)
   \put(20.0,40.0){\circle2}
   \put(33.0, 0.0){\circle2}
   \put( 7.0, 0.0){\circle2}
   \put(41.8,25.5){\circle2} 
   \put(-1.8,25.5){\circle2}
   \put(20.0,29.0){\circle2}
   \put(12.2,16.5){\circle2}
   \put(27.8,16.5){\circle2}
   \put( 8.0, 0.0){\line(+1, 0){24.0}}
   \put( 6.5, 0.8){\line(-1,+3){ 7.8}}
   \put(33.5, 0.8){\line(+1,+3){ 7.9}}
   \put(19.2,39.5){\line(-3,-2){20.0}}
   \put(20.8,39.5){\line(+3,-2){20.0}}
   \put(19.4,28.2){\line(-2,-3){ 7.0}}
   \put(20.6,28.2){\line(+2,-3){ 7.0}}
   \put(13.2,16.5){\line(+1, 0){13.6}}
   \multiput( 8.7, 1.0)(+1.9,+0.0){13}{\circle*1}
   \multiput(20.0,30.9)( 0.0,+1.8){5}{\circle*1}
   \multiput(10.5,17.7)(-1.5,+1.0){8}{\circle*1}
   \multiput(29.5,17.7)(+1.5,+1.0){8}{\circle*1}
   \put( 5.4,12.3){$\scriptscriptstyle A$}
   \put(12.5,28.1){$\scriptscriptstyle B$}
   \put(31.0,21.7){$\scriptscriptstyle*$}
  \end{picture}}
  \;\;-\;
  \raise-17pt\hbox{\begin{picture}(40,40)(0,0)
   \put(11.7, 0.0){\circle2}
   \put(28.3, 0.0){\circle2}
   \put( 0.0,11.7){\circle2}
   \put( 0.0,28.3){\circle2}
   \put(11.7,40.0){\circle2}
   \put(28.3,40.0){\circle2}
   \put(40.0,11.7){\circle2}
   \put(40.0,28.3){\circle2}
   \put(12.7, 0.0){\line(+1, 0){14.6}}
   \put(12.7,40.0){\line(+1, 0){14.6}}
   \put( 0.0,12.7){\line( 0,+1){14.6}}
   \put(40.0,12.7){\line( 0,+1){14.6}}
   \put(10.9, 0.6){\line(-1,+1){10.2}}
   \put(29.1, 0.6){\line(+1,+1){10.2}}
   \put(10.9,39.4){\line(-1,-1){10.2}}
   \put(29.1,39.4){\line(+1,-1){10.2}}
   \multiput(13.4, 1.0)(+1.9,+0.0){8}{\circle*1}
   \multiput(13.1,38.1)(+1.2,-1.2){22}{\circle*1}
   \multiput( 2.2,12.5)(+1.7,+0.7){22}{\circle*1}
   \multiput( 2.2,28.5)(+1.6,+0.7){16}{\circle*1}
   \put(29.0,41.5){$\scriptscriptstyle B$}
   \put(41.0,28.7){$\scriptscriptstyle A$}
  \end{picture}}
  \;-\;
  \raise-17pt\hbox{\begin{picture}(40,40)(0,0)
   \put(11.7, 0.0){\circle2}
   \put(28.3, 0.0){\circle2}
   \put( 0.0,11.7){\circle2}
   \put( 0.0,28.3){\circle2}
   \put(11.7,40.0){\circle2}
   \put(28.3,40.0){\circle2}
   \put(40.0,11.7){\circle2}
   \put(40.0,28.3){\circle2}
   \put(12.7, 0.0){\line(+1, 0){14.6}}
   \put(12.7,40.0){\line(+1, 0){14.6}}
   \put( 0.0,12.7){\line( 0,+1){14.6}}
   \put(40.0,12.7){\line( 0,+1){14.6}}
   \put(10.9, 0.6){\line(-1,+1){10.2}}
   \put(29.1, 0.6){\line(+1,+1){10.2}}
   \put(10.9,39.4){\line(-1,-1){10.2}}
   \put(29.1,39.4){\line(+1,-1){10.2}}
   \multiput(13.4, 1.0)(+1.9, 0.0){8}{\circle*1}
   \multiput(12.9,38.1)(+1.2,-1.2){22}{\circle*1}
   \multiput(27.1,38.1)(-1.2,-1.2){22}{\circle*1}
   \multiput( 1.9,28.3)(+1.9, 0.0){20}{\circle*1}
   \put(29.0,41.5){$\scriptscriptstyle A$}
   \put(41.0,28.7){$\scriptscriptstyle B$}
  \end{picture}}
  \;\;\equiv\;\;
  0
 \end{eqnarray*}
 and
 \begin{eqnarray*}
  \lefteqn{\raise-17pt\hbox{\begin{picture}(40,40)(0,0)
   \put(11.7, 0.0){\circle2}
   \put(28.3, 0.0){\circle2}
   \put( 0.0,11.7){\circle2}
   \put( 0.0,28.3){\circle2}
   \put(11.7,40.0){\circle2}
   \put(28.3,40.0){\circle2}
   \put(40.0,11.7){\circle2}
   \put(40.0,28.3){\circle2}
   \put(12.7, 0.0){\line(+1, 0){14.6}}
   \put(12.7,40.0){\line(+1, 0){14.6}}
   \put( 0.0,12.7){\line( 0,+1){14.6}}
   \put(40.0,12.7){\line( 0,+1){14.6}}
   \put(10.9, 0.6){\line(-1,+1){10.2}}
   \put(29.1, 0.6){\line(+1,+1){10.2}}
   \put(10.9,39.4){\line(-1,-1){10.2}}
   \put(29.1,39.4){\line(+1,-1){10.2}}
   \multiput(11.7, 2.0)( 0.0,+1.9){20}{\circle*1}
   \multiput(28.3, 2.0)( 0.0,+1.9){20}{\circle*1}
   \multiput( 2.0,11.7)(+1.9, 0.0){20}{\circle*1}
   \multiput( 2.0,28.3)(+1.9, 0.0){20}{\circle*1}
  \end{picture}}
  \;+\;
  \raise-17pt\hbox{\begin{picture}(40,40)(0,0)
   \put(11.7, 0.0){\circle2}
   \put(28.3, 0.0){\circle2}
   \put( 0.0,11.7){\circle2}
   \put( 0.0,28.3){\circle2}
   \put(11.7,40.0){\circle2}
   \put(28.3,40.0){\circle2}
   \put(40.0,11.7){\circle2}
   \put(40.0,28.3){\circle2}
   \put(12.7, 0.0){\line(+1, 0){14.6}}
   \put(12.7,40.0){\line(+1, 0){14.6}}
   \put( 0.0,12.7){\line( 0,+1){14.6}}
   \put(40.0,12.7){\line( 0,+1){14.6}}
   \put(10.9, 0.6){\line(-1,+1){10.2}}
   \put(29.1, 0.6){\line(+1,+1){10.2}}
   \put(10.9,39.4){\line(-1,-1){10.2}}
   \put(29.1,39.4){\line(+1,-1){10.2}}
   \multiput(12.4, 2.0)(+0.8,+1.9){20}{\circle*1}
   \multiput(27.6, 2.0)(-0.8,+1.9){20}{\circle*1}
   \multiput( 2.0,11.7)(+1.9, 0.0){20}{\circle*1}
   \multiput( 2.0,28.3)(+1.9, 0.0){20}{\circle*1}
  \end{picture}}
  \;-\;
  \raise-17pt\hbox{\begin{picture}(40,40)(0,0)
   \put(11.7, 0.0){\circle2}
   \put(28.3, 0.0){\circle2}
   \put( 0.0,11.7){\circle2}
   \put( 0.0,28.3){\circle2}
   \put(11.7,40.0){\circle2}
   \put(28.3,40.0){\circle2}
   \put(40.0,11.7){\circle2}
   \put(40.0,28.3){\circle2}
   \put(12.7, 0.0){\line(+1, 0){14.6}}
   \put(12.7,40.0){\line(+1, 0){14.6}}
   \put( 0.0,12.7){\line( 0,+1){14.6}}
   \put(40.0,12.7){\line( 0,+1){14.6}}
   \put(10.9, 0.6){\line(-1,+1){10.2}}
   \put(29.1, 0.6){\line(+1,+1){10.2}}
   \put(10.9,39.4){\line(-1,-1){10.2}}
   \put(29.1,39.4){\line(+1,-1){10.2}}
   \multiput( 2.3,11.1)(+1.5,-0.6){17}{\circle*1}
   \multiput( 2.3,28.9)(+1.5,+0.6){17}{\circle*1}
   \multiput(38.1,12.9)(-1.2,+1.2){22}{\circle*1}
   \multiput(38.1,27.1)(-1.2,-1.2){22}{\circle*1}
  \end{picture}}
  \;-\;
  \raise-17pt\hbox{\begin{picture}(40,40)(0,0)
   \put(11.7, 0.0){\circle2}
   \put(28.3, 0.0){\circle2}
   \put( 0.0,11.7){\circle2}
   \put( 0.0,28.3){\circle2}
   \put(11.7,40.0){\circle2}
   \put(28.3,40.0){\circle2}
   \put(40.0,11.7){\circle2}
   \put(40.0,28.3){\circle2}
   \put(12.7, 0.0){\line(+1, 0){14.6}}
   \put(12.7,40.0){\line(+1, 0){14.6}}
   \put( 0.0,12.7){\line( 0,+1){14.6}}
   \put(40.0,12.7){\line( 0,+1){14.6}}
   \put(10.9, 0.6){\line(-1,+1){10.2}}
   \put(29.1, 0.6){\line(+1,+1){10.2}}
   \put(10.9,39.4){\line(-1,-1){10.2}}
   \put(29.1,39.4){\line(+1,-1){10.2}}
   \multiput( 2.0,11.7)(+1.9, 0.0){20}{\circle*1}
   \multiput(12.9, 1.9)(+1.2,+1.2){22}{\circle*1}
   \multiput(27.6, 2.0)(-0.8,+1.9){20}{\circle*1}
   \multiput( 2.3,28.9)(+1.5,+0.6){17}{\circle*1}
  \end{picture}}}
  \qquad
  &&
  \\[10pt]
  &=&
  \raise-17pt\hbox{\begin{picture}(40,40)(0,0)
   \put(20.0,40.0){\circle2}
   \put(33.0, 0.0){\circle2}
   \put( 7.0, 0.0){\circle2}
   \put(41.8,25.5){\circle2} 
   \put(-1.8,25.5){\circle2}
   \put(20.0,21.0){\circle2}
   \put(12.2, 8.5){\circle2}
   \put(27.8, 8.5){\circle2}
   \put( 8.0, 0.0){\line(+1, 0){24.0}}
   \put( 6.5, 0.8){\line(-1,+3){ 7.8}}
   \put(33.5, 0.8){\line(+1,+3){ 7.9}}
   \put(19.2,39.5){\line(-3,-2){20.0}}
   \put(20.8,39.5){\line(+3,-2){20.0}}
   \put(19.4,20.2){\line(-2,-3){ 7.0}}
   \put(20.6,20.2){\line(+2,-3){ 7.0}}
   \put(13.2, 8.5){\line(+1, 0){13.6}}
   \multiput(20.0,22.9)( 0.0,+1.9){9}{\circle*1}
   \multiput( 0.0,25.5)(+1.9, 0.0){22}{\circle*1}
   \multiput(11.2, 6.6)(-0.7,-1.2){5}{\circle*1}
   \multiput(28.8, 6.6)(+0.8,-1.2){5}{\circle*1}
   \put(21.0,30.0){$\scriptscriptstyle*$}
   \put( 6.8,10.5){$\scriptscriptstyle A$}
   \put(27.5,10.5){$\scriptscriptstyle B$}
  \end{picture}}
  \;\;+\;
  \raise-17pt\hbox{\begin{picture}(40,40)(0,0)
   \put(11.7, 0.0){\circle2}
   \put(28.3, 0.0){\circle2}
   \put( 0.0,11.7){\circle2}
   \put( 0.0,28.3){\circle2}
   \put(11.7,40.0){\circle2}
   \put(28.3,40.0){\circle2}
   \put(40.0,11.7){\circle2}
   \put(40.0,28.3){\circle2}
   \put(12.7, 0.0){\line(+1, 0){14.6}}
   \put(12.7,40.0){\line(+1, 0){14.6}}
   \put( 0.0,12.7){\line( 0,+1){14.6}}
   \put(40.0,12.7){\line( 0,+1){14.6}}
   \put(10.9, 0.6){\line(-1,+1){10.2}}
   \put(29.1, 0.6){\line(+1,+1){10.2}}
   \put(10.9,39.4){\line(-1,-1){10.2}}
   \put(29.1,39.4){\line(+1,-1){10.2}}
   \multiput(11.7, 2.0)( 0.0,+1.9){20}{\circle*1}
   \multiput(28.3, 2.0)( 0.0,+1.9){20}{\circle*1}
   \multiput( 2.0,11.7)(+1.9, 0.0){20}{\circle*1}
   \multiput( 2.0,28.3)(+1.9, 0.0){20}{\circle*1}
   \put( 8.7,42.4){$\scriptscriptstyle A$}
   \put(26.0,42.4){$\scriptscriptstyle B$}
  \end{picture}}
  \;\;+\;
  \raise-17pt\hbox{\begin{picture}(40,40)(0,0)
   \put(11.7, 0.0){\circle2}
   \put(28.3, 0.0){\circle2}
   \put( 0.0,11.7){\circle2}
   \put( 0.0,28.3){\circle2}
   \put(11.7,40.0){\circle2}
   \put(28.3,40.0){\circle2}
   \put(40.0,11.7){\circle2}
   \put(40.0,28.3){\circle2}
   \put(12.7, 0.0){\line(+1, 0){14.6}}
   \put(12.7,40.0){\line(+1, 0){14.6}}
   \put( 0.0,12.7){\line( 0,+1){14.6}}
   \put(40.0,12.7){\line( 0,+1){14.6}}
   \put(10.9, 0.6){\line(-1,+1){10.2}}
   \put(29.1, 0.6){\line(+1,+1){10.2}}
   \put(10.9,39.4){\line(-1,-1){10.2}}
   \put(29.1,39.4){\line(+1,-1){10.2}}
   \multiput(12.4, 2.0)(+0.8,+1.9){20}{\circle*1}
   \multiput(27.6, 2.0)(-0.8,+1.9){20}{\circle*1}
   \multiput( 2.0,11.7)(+1.9, 0.0){20}{\circle*1}
   \multiput( 2.0,28.3)(+1.9, 0.0){20}{\circle*1}
   \put( 8.7,42.4){$\scriptscriptstyle B$}
   \put(26.0,42.4){$\scriptscriptstyle A$}
  \end{picture}}
  \;-\;\; 
  \raise-17pt\hbox{\begin{picture}(40,40)(0,0)
   \put(20.0,40.0){\circle2}
   \put(33.0, 0.0){\circle2}
   \put( 7.0, 0.0){\circle2}
   \put(41.8,25.5){\circle2} 
   \put(-1.8,25.5){\circle2}
   \put(20.0,21.0){\circle2}
   \put(12.2, 8.5){\circle2}
   \put(27.8, 8.5){\circle2}
   \put( 8.0, 0.0){\line(+1, 0){24.0}}
   \put( 6.5, 0.8){\line(-1,+3){ 7.8}}
   \put(33.5, 0.8){\line(+1,+3){ 7.9}}
   \put(19.2,39.5){\line(-3,-2){20.0}}
   \put(20.8,39.5){\line(+3,-2){20.0}}
   \put(19.4,20.2){\line(-2,-3){ 7.0}}
   \put(20.6,20.2){\line(+2,-3){ 7.0}}
   \put(13.2, 8.5){\line(+1, 0){13.6}}
   \multiput(20.0,22.9)( 0.0,+1.9){9}{\circle*1}
   \multiput( 0.0,25.5)(+1.9, 0.0){22}{\circle*1}
   \multiput(11.2, 6.6)(-0.7,-1.2){5}{\circle*1}
   \multiput(28.8, 6.6)(+0.8,-1.2){5}{\circle*1}
   \put(30.5, 4.3){$\scriptscriptstyle*$}
   \put( 6.4,10.2){$\scriptscriptstyle A$}
   \put(12.4,19.9){$\scriptscriptstyle B$}
  \end{picture}}
  \;\;-\;
  \raise-17pt\hbox{\begin{picture}(40,40)(0,0)
   \put(11.7, 0.0){\circle2}
   \put(28.3, 0.0){\circle2}
   \put( 0.0,11.7){\circle2}
   \put( 0.0,28.3){\circle2}
   \put(11.7,40.0){\circle2}
   \put(28.3,40.0){\circle2}
   \put(40.0,11.7){\circle2}
   \put(40.0,28.3){\circle2}
   \put(12.7, 0.0){\line(+1, 0){14.6}}
   \put(12.7,40.0){\line(+1, 0){14.6}}
   \put( 0.0,12.7){\line( 0,+1){14.6}}
   \put(40.0,12.7){\line( 0,+1){14.6}}
   \put(10.9, 0.6){\line(-1,+1){10.2}}
   \put(29.1, 0.6){\line(+1,+1){10.2}}
   \put(10.9,39.4){\line(-1,-1){10.2}}
   \put(29.1,39.4){\line(+1,-1){10.2}}
   \multiput( 2.3,11.1)(+1.5,-0.6){17}{\circle*1}
   \multiput( 2.3,28.9)(+1.5,+0.6){17}{\circle*1}
   \multiput(38.1,12.9)(-1.2,+1.2){22}{\circle*1}
   \multiput(38.1,27.1)(-1.2,-1.2){22}{\circle*1}
   \put(27.1,-5.0){$\scriptscriptstyle A$}
   \put(41.0, 8.5){$\scriptscriptstyle B$}
  \end{picture}}
  \;-\;
  \raise-17pt\hbox{\begin{picture}(40,40)(0,0)
   \put(11.7, 0.0){\circle2}
   \put(28.3, 0.0){\circle2}
   \put( 0.0,11.7){\circle2}
   \put( 0.0,28.3){\circle2}
   \put(11.7,40.0){\circle2}
   \put(28.3,40.0){\circle2}
   \put(40.0,11.7){\circle2}
   \put(40.0,28.3){\circle2}
   \put(12.7, 0.0){\line(+1, 0){14.6}}
   \put(12.7,40.0){\line(+1, 0){14.6}}
   \put( 0.0,12.7){\line( 0,+1){14.6}}
   \put(40.0,12.7){\line( 0,+1){14.6}}
   \put(10.9, 0.6){\line(-1,+1){10.2}}
   \put(29.1, 0.6){\line(+1,+1){10.2}}
   \put(10.9,39.4){\line(-1,-1){10.2}}
   \put(29.1,39.4){\line(+1,-1){10.2}}
   \multiput( 2.0,11.7)(+1.9, 0.0){20}{\circle*1}
   \multiput(12.9, 1.9)(+1.2,+1.2){22}{\circle*1}
   \multiput(27.6, 2.0)(-0.8,+1.9){20}{\circle*1}
   \multiput( 2.3,28.9)(+1.5,+0.6){17}{\circle*1}
   \put(28.2,-5.0){$\scriptscriptstyle B$}
   \put(40.6, 8.5){$\scriptscriptstyle A$}
  \end{picture}}
  \;\;\equiv\;\;
  0
  \\[-5pt]
  &&
 \end{eqnarray*}
 in the graph algebra $\A^\bullet$ modulo the ideal of IHX--relations. For
 the convenience of the reader we have marked the relevant red edge for these
 IHX--relations with $*$, in addition we have traced the images of a pair of
 interesting vertices marking them with the letters $A$ and $B$.

 \pfill
 Before closing this section let us briefly discuss a modification of the
 construction of the algebra $\A^\bullet$ of colored trivalent graphs and
 its quotient $\overline\A^\bullet$ by the IHX--relations, which allows us
 to use the curvature tensor $R$ directly in the emerging graphical calculus
 for stable curvature invariants. For this purpose we need to allow tetravalent
 besides trivalent vertices, moreover all these tetravalent vertices need to
 come along with an orientation and an unordered partition of the set of
 adjacent flags into two pairs of flags:

 \begin{Definition}[Extended Colored Trivalent Graphs]
 \hfill\label{ectg}\break
 An extended colored trivalent graph is a graph $\gamma$ with tri-- and
 tetravalent vertices endowed with a coloring of its edges in red \& black
 as before such that every trivalent vertex is adjacent to exactly one red
 flag and vice versa. Moreover the set $\Flag_v\gamma\,=\,\{\,f_1,f_2,f_3,
 f_4\,\}$ of black flags adjacent to each tetravalent vertex $v\,\in\,\Vert
 \,\gamma$ is endowed with an orientation and an unordered partition
 $\Flag_v\gamma\,=\,\{\,f_1,f_2\,\}\,\dot\cup\,\{\,f_3,f_4\,\}$ into
 two pairs of flags.
 \end{Definition}

 \noindent
 Concerning homomorphisms and isomorphisms we stipulate that
 every homomorphism between extended colored trivalent graphs
 $\varphi:\,\gamma\longrightarrow\hat\gamma$ needs to respect the
 $\theta$--invariant coloring $\col:\,\Flag\,\gamma\longrightarrow
 \{\red,\black\}$ of the flags and the orientation as well as
 the partition of the sets of flags adjacent to tetravalent vertices. In
 diagrams we will depict tetravalent vertices by circles with horizontal
 or vertical chords indicating the partition of flags into two pairs:
 $$
  \raise-22pt
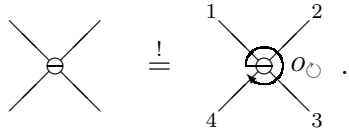
\hbox{\begin{picture}(50,50)
   \put( 25  ,25  ){\circle6}
   \put( 22  ,25  ){\line(+1, 0){6}}
   \put( 27.2,27.2){\line(+1,+1){15}}
   \put( 22.8,27.2){\line(-1,+1){15}}
   \put( 27.2,22.8){\line(+1,-1){15}}
   \put( 22.8,22.8){\line(-1,-1){15}}
  \end{picture}}
  \;\;\stackrel!=\;\;
  \raise-22pt\hbox{\begin{picture}(50,50)
   \put( 25  ,25  ){\circle6}
   \put( 18  ,25  ){\line(+1, 0){10}}
   \put( 35  ,23  ){$o_\circlearrowright$}
   \put( 27.2,27.2){\line(+1,+1){15}}
   \put( 22.8,27.2){\line(-1,+1){15}}
   \put( 27.2,22.8){\line(+1,-1){15}}
   \put( 22.8,22.8){\line(-1,-1){15}}
   \put( 18.3,21.9){\vector(-1,+1){0}}
   \put( 43,43){$\scriptstyle2$}
   \put(  3,43){$\scriptstyle1$}
   \put( 43, 2){$\scriptstyle3$}
   \put(  3, 2){$\scriptstyle4$}
   \qbezier(18.0,25.0)(18.0,28.2)(20.1,29.9)
   \qbezier(20.1,29.9)(21.8,32.0)(25.0,32.0) 
   \qbezier(25.0,32.0)(28.2,32.0)(29.9,29.9)
   \qbezier(29.9,29.9)(32.0,28.2)(32.0,25.0)
   \qbezier(32.0,25.0)(32.0,21.8)(29.9,20.1)
   \qbezier(29.9,20.1)(28.2,18.0)(25.0,18.0)
   \qbezier(25.0,18.0)(21.8,18.0)(20.1,20.1)
  \end{picture}}\ .
 $$
 Unless stated otherwise we will tacitly assume that the orientation
 on the set of flags adjacent to such a tetravalent vertex is the variant
 $o_\circlearrowright$ of the so called blackboard orientation: Starting
 on one of the two end points of the chord and continuing counterclockwise
 around the vertex results in a labelling $L$ of the adjacent flags
 representing $o_\circlearrowright\,:=\,[\,L,\,+1\,]$. Note that this
 orientation does not depend on which end of the chord we begin with.
 
 In passing we remark that the joint stabilizer of the unordered partition
 $\{\{1,2\},\{3,4\}\}$ into pairs and the tautological orientation of the
 set $\{1,2,3,4\}$ in the symmetric group $S_4$ equals the normal Kleinian
 Four subgroup $K\,\subseteq\,S_4$. Instead of requiring both an orientation
 and a partition into pairs we may hence require alternatively that the set
 of flags adjacent to each tetravalent vertex $v$ of the graph $\gamma$ is
 decorated by an equivalence class of a labelling $L:\,\Flag_v\gamma\stackrel
 \cong\longrightarrow\{\,1,2,3,4\,\}$ modulo postcomposition with elements of
 $K$. In this alternative formulation of Definition \ref{ectg} the six
 different decorations on a tetravalent vertex read
 $$
  
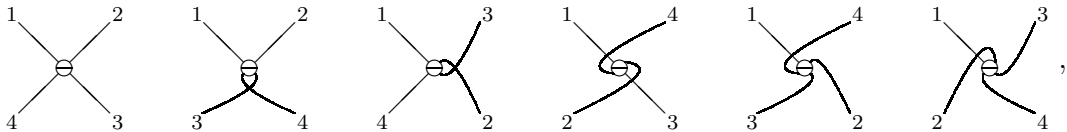
\begin{picture}(400,50)
   \put( 25  ,25  ){\circle6}
   \put( 22  ,25  ){\line(+1, 0){6}}
   \put( 27.2,27.2){\line(+1,+1){15}}
   \put( 22.8,27.2){\line(-1,+1){15}}
   \put( 27.2,22.8){\line(+1,-1){15}}
   \put( 22.8,22.8){\line(-1,-1){15}}
   \put( 43,43){$\scriptstyle2$}
   \put(  3,43){$\scriptstyle1$}
   \put( 43, 2){$\scriptstyle3$}
   \put(  3, 2){$\scriptstyle4$}
   \put( 95  ,25  ){\circle6}
   \put( 92  ,25  ){\line(+1, 0){6}}
   \put( 97.2,27.2){\line(+1,+1){15}}
   \put( 92.8,27.2){\line(-1,+1){15}}
   \qbezier( 97.2,22.8)(101.0,18.0)( 77.2, 7.8)
   \qbezier( 92.8,22.8)( 89.0,18.0)(112.2, 7.8)
   \put(113,43){$\scriptstyle2$}
   \put( 73,43){$\scriptstyle1$}
   \put(113, 2){$\scriptstyle4$}
   \put( 73, 2){$\scriptstyle3$}
   \put(165  ,25  ){\circle6}
   \put(162  ,25  ){\line(+1, 0){6}}
   \qbezier(167.2,27.2)(172.0,31.0)(182.2, 7.8)
   \put(162.8,27.2){\line(-1,+1){15}}
   \qbezier(167.2,22.8)(172.0,19.0)(182.2,42.2)
   \put(162.8,22.8){\line(-1,-1){15}}
   \put(183,43){$\scriptstyle3$}
   \put(143,43){$\scriptstyle1$}
   \put(183, 2){$\scriptstyle2$}
   \put(143, 2){$\scriptstyle4$}
   \put(235  ,25  ){\circle6}
   \put(232  ,25  ){\line(+1, 0){6}}
   \qbezier(232.8,22.8)(216.0,25.0)(252.2,42.2)
   \put(232.8,27.2){\line(-1,+1){15}}
   \put(237.2,22.8){\line(+1,-1){15}}
   \qbezier(237.2,27.2)(254.0,25.0)(217.8, 7.8)
   \put(253,43){$\scriptstyle4$}
   \put(213,43){$\scriptstyle1$}
   \put(253, 2){$\scriptstyle3$}
   \put(213, 2){$\scriptstyle2$}
   \put(305  ,25  ){\circle6}
   \put(302  ,25  ){\line(+1, 0){6}}
   \qbezier(302.8,22.8)(286.0,25.0)(322.2,42.2)
   \put(302.8,27.2){\line(-1,+1){15}}
   \qbezier(307.2,27.2)(311.0,32.0)(322.2, 7.8)
   \qbezier(307.2,22.8)(312.0,19.0)(287.2, 7.8)
   \put(323,43){$\scriptstyle4$}
   \put(283,43){$\scriptstyle1$}
   \put(323, 2){$\scriptstyle2$}
   \put(283, 2){$\scriptstyle3$}
   \put(375  ,25  ){\circle6}
   \put(372  ,25  ){\line(+1, 0){6}}
   \qbezier(377.2,22.8)(382.0,19.0)(392.2,42.2)
   \put(372.8,27.2){\line(-1,+1){15}}
   \qbezier(372.8,22.8)(369.0,18.0)(392.2, 7.8)
   \qbezier(377.2,27.2)(375.0,44.0)(357.8, 7.8)
   \put(393,43){$\scriptstyle3$}
   \put(353,43){$\scriptstyle1$}
   \put(393, 2){$\scriptstyle4$}
   \put(353, 2){$\scriptstyle2$}
   \put(401,25){,}
  \end{picture}
 $$
 where the numbers indicate a representative labelling $L$ under
 postcomposition with elements of $K\,\subseteq\,S_4$. Mimicking the
 construction of the algebra $\A^\bullet$ of colored trivalent graphs
 and its quotient $\overline\A^\bullet$ by the IHX--relations we
 consider the set of isomorphism classes
 $$
  \Gamma^n_\ext
  \;\;:=\;\;
  \{\;\;[\,\gamma\,]\;\;|\;\;\gamma\textrm{\ extended colored trivalent graph
   with\ }n\,=\,{\textstyle\frac12}\,\#\Flag\,\gamma-\#\Vert\,\gamma\;\;\}
 $$
 of extended colored trivalent graphs with $2n$ vertices counting all
 tetravalent vertices twice. Unlike $\A^\bullet$ the algebra $\A^\bullet_\ext$
 of extended colored trivalent graphs is not the convolution algebra
 $\R\,\Gamma^\bullet_\ext$ associated to the commutative monoid
 $\Gamma^\bullet_\ext$, but its quotient by the ideal $\langle\,
 \mathrm{O}\,\rangle\,\subseteq\,\R\,\Gamma^\bullet_\ext$ spanned
 by a change of orientation for the set of flags adjacent to some
 tetravalent vertex:
 \begin{equation}\label{ox}
  \raise-22pt\hbox{\begin{picture}(50,50)
   \put( 25  ,25  ){\circle6}
   \put( 22  ,25  ){\line(+1, 0){6}}
   \put( 27.2,27.2){\line(+1,+1){15}}
   \put( 22.8,27.2){\line(-1,+1){15}}
   \put( 27.2,22.8){\line(+1,-1){15}}
   \put( 22.8,22.8){\line(-1,-1){15}}
   \put( 31  ,21  ){$+o$}
  \end{picture}}
  \;\;\;+\;
  \raise-22pt\hbox{\begin{picture}(50,50)
   \put( 25  ,25  ){\circle6}
   \put( 22  ,25  ){\line(+1, 0){6}}
   \put( 27.2,27.2){\line(+1,+1){15}}
   \put( 22.8,27.2){\line(-1,+1){15}}
   \put( 27.2,22.8){\line(+1,-1){15}}
   \put( 22.8,22.8){\line(-1,-1){15}}
   \put( 31  ,21  ){$-o$\ .}
  \end{picture}}
 \end{equation}
 Endowed with the $\R$--bilinear extension of the disjoint union the quotient
 $\A^\bullet_\ext\,:=\,\R\,\Gamma^\bullet_\ext/\langle\,\mathrm{O}\,\rangle$
 becomes a graded commutative algebra with its own version $\langle\,\IHX\,
 \rangle\,\subseteq\,\A^\bullet_\ext$ of the ideal of IHX--relations defined
 for all red edges by equations (\ref{ihxex}) and (\ref{span}) as before. The
 quotient
 $$
  \overline\A^\bullet_\ext
  \;\;=\;\;
  \A^\bullet_\ext/_{\displaystyle\langle\;\IHX\;\rangle}
 $$
 is the reduced algebra of extended colored trivalent graphs. In order to
 extend the definition of stable curvature invariants from $\A^\bullet$ to
 $\A^\bullet_\ext$ we specify the additional vertex interaction
 \begin{equation}\label{feyn2}
  \raise-17pt\hbox{\begin{picture}(53,42)
   \put(26,21){\circle6}
   \put(28.2,23.2){\line(+1,+1){14}}\put(43,33){$Y$}
   \put(23.8,23.2){\line(-1,+1){14}}\put( 0,33){$X$}
   \put(28.2,18.8){\line(+1,-1){14}}\put(43, 1){$U$}
   \put(23.8,18.8){\line(-1,-1){14}}\put( 1, 1){$V$}
   \put(23,21){\line(+1, 0){6}}
  \end{picture}}
  \;\;\widehat=\;\;
  R(\;X,\;Y;\;U,\;V\;)
 \end{equation}
 besides the red edge interaction (\ref{feyn1}). The entire point of this
 construction is that the extended graph algebra $\A^\bullet_\ext$ comes
 along with a surjective algebra homomorphism
 $$
  \Phi^-:\;\;\A^\bullet_\ext\;\longrightarrow\;\A^\bullet\ ,
 $$
 which descends to a surjective algebra homomorphism $\overline\Phi^-:\,
 \overline\A^\bullet_\ext\longrightarrow\overline\A^\bullet$ between the
 respective quotients by the ideal of IHX--relations. More precisely the
 algebra homomorphism $\Phi^-$ expands every tetravalent vertex into a pair
 of trivalent vertices connected by a red edge
 \begin{equation}\label{rex}
  \raise-17pt\hbox{\begin{picture}(40,40)
   \put(20,20){\circle6}
   \put(22.2,22.2){\line(+1,+1){19}}
   \put(17.8,22.2){\line(-1,+1){19}}
   \put(22.2,17.8){\line(+1,-1){19}}
   \put(17.8,17.8){\line(-1,-1){19}}
   \put(17,20){\line(+1, 0){6}}
  \end{picture}}
  \qquad\stackrel{\Phi^-}\longmapsto\qquad
  \frac16\;\raise-17pt\hbox{\begin{picture}(40,40)
   \put( 9,20){\circle2}
   \put(31,20){\circle2}
   \multiput(11,20)(+2, 0){10}{\circle*1}
   \put( 8,21){\line(-1,+2){10}}
   \put( 8,19){\line(-1,-2){10}}
   \put(32,21){\line(+1,+2){10}}
   \put(32,19){\line(+1,-2){10}}
  \end{picture}}
  \;-\;\frac16\;\raise-17pt\hbox{\begin{picture}(40,40)
   \put( 9,20){\circle2}
   \put(31,20){\circle2}
   \multiput(11,20)(+2, 0){10}{\circle*1}
   \put(30,21){\line(-3,+2){30}}
   \put(10,21){\line(+3,+2){30}}
   \put( 8,19){\line(-1,-2){10}}
   \put(32,19){\line(+1,-2){10}}
  \end{picture}}
 \end{equation}
 in accordance with the Feynman rules interpretation of the identity
 $R\,=\,\Phi^-\Sec$, namely:
 $$
  R(\;X,Y;\,U,V\;)
  \;\;=\;\;
  \frac16\;\Sec(\;X,V;\,Y,U\;)\;-\;\frac16\;\Sec(\;X,U;\,Y,V\;)\ .
 $$
 An example of the usefulness of considering $\A^\bullet_\ext$ in addition to
 $\A^\bullet$ is the congruence
 \begin{eqnarray*}
  \Phi^-\left(\;\frac14\,\raise-15pt\hbox{\begin{picture}(60,36)
   \put(10,18){\circle6}
   \put(50,18){\circle6}
   \put( 7,18){\line(+1,0){6}}
   \put(47,18){\line(+1,0){6}}
   \qbezier(12.2,20.2)(69,44)(52.2,20.2)
   \qbezier( 7.8,20.2)(-9,44)(47.8,20.2)
   \qbezier(12.2,15.8)(69,-8)(52.2,15.8)
   \qbezier( 7.8,15.8)(-9,-8)(47.8,15.8)
  \end{picture}}\;\right)
  &=&
  \frac1{144}\left(
  \raise-15pt\hbox{\begin{picture}(60,36)
   \put( 5,18){\circle2}
   \put(21,18){\circle2}
   \put(39,18){\circle2}
   \put(55,18){\circle2}
   \multiput( 7,18)(+2, 0){7}{\circle*1}
   \multiput(41,18)(+2, 0){7}{\circle*1}
   \qbezier( 5.0,19.0)(10.0,33.0)(22.0,33.0)
   \qbezier(22.0,33.0)(34.0,33.0)(39.0,19.0)
   \qbezier(21.0,19.0)(26.0,33.0)(38.0,33.0)
   \qbezier(38.0,33.0)(50.0,33.0)(55.0,19.0)
   \qbezier( 5.0,17.0)(10.0, 3.0)(22.0, 3.0)
   \qbezier(22.0, 3.0)(34.0, 3.0)(39.0,17.0)
   \qbezier(21.0,17.0)(26.0, 3.0)(38.0, 3.0)
   \qbezier(38.0, 3.0)(50.0, 3.0)(55.0,17.0)
  \end{picture}}
  -
  \raise-15pt\hbox{\begin{picture}(60,36)
   \put( 5,18){\circle2}
   \put(21,18){\circle2}
   \put(39,18){\circle2}
   \put(55,18){\circle2}
   \multiput( 7,18)(+2, 0){7}{\circle*1}
   \multiput(41,18)(+2, 0){7}{\circle*1}
   \qbezier( 5.0,19.0)(10.0,33.0)(22.0,33.0)
   \qbezier(22.0,33.0)(34.0,33.0)(39.0,19.0)
   \qbezier(21.0,19.0)(26.0,33.0)(38.0,33.0)
   \qbezier(38.0,33.0)(50.0,33.0)(55.0,19.0)
   \qbezier(21.8,17.2)(30.0, 2.0)(38.2,17.2)
   \qbezier(30.0, 1.0)(50.0, 2.0)(54.4,17.1)
   \qbezier( 5.5,17.0)(10.0, 2.0)(30.0, 1.0)
  \end{picture}}
  -
  \raise-15pt\hbox{\begin{picture}(60,36)
   \put( 5,18){\circle2}
   \put(21,18){\circle2}
   \put(39,18){\circle2}
   \put(55,18){\circle2}
   \multiput( 7,18)(+2, 0){7}{\circle*1}
   \multiput(41,18)(+2, 0){7}{\circle*1}
   \qbezier(21.8,18.8)(30.0,34.0)(38.2,18.8)
   \qbezier( 5.5,19.0)(10.0,34.0)(30.0,35.0)
   \qbezier(30.0,35.0)(50.0,34.0)(54.4,19.0)
   \qbezier( 5.0,17.0)(10.0, 3.0)(22.0, 3.0)
   \qbezier(22.0, 3.0)(34.0, 3.0)(39.0,17.0)
   \qbezier(21.0,17.0)(26.0, 3.0)(38.0, 3.0)
   \qbezier(38.0, 3.0)(50.0, 3.0)(55.0,17.0)
  \end{picture}}
  +
  \raise-15pt\hbox{\begin{picture}(60,36)
   \put( 5,18){\circle2}
   \put(21,18){\circle2}
   \put(39,18){\circle2}
   \put(55,18){\circle2}
   \multiput( 7,18)(+2, 0){7}{\circle*1}
   \multiput(41,18)(+2, 0){7}{\circle*1}
   \qbezier(21.8,18.8)(30.0,34.0)(38.2,18.8)
   \qbezier( 5.5,19.0)(10.0,34.0)(30.0,35.0)
   \qbezier(30.0,35.0)(50.0,34.0)(54.4,19.0)
   \qbezier(21.8,17.2)(30.0, 2.0)(38.2,17.2)
   \qbezier(30.0, 1.0)(50.0, 2.0)(54.4,17.1)
   \qbezier( 5.5,17.0)(10.0, 2.0)(30.0, 1.0)
  \end{picture}}\right)
  \\[4pt]
  &=&
  \frac1{72}\;
  \raise-8pt\hbox{\begin{picture}(24,24)
   \put( 0, 0){\circle2}
   \put( 0,24){\circle2}
   \put(24, 0){\circle2}
   \put(24,24){\circle2}
   \multiput( 0, 2)( 0,+2){11}{\circle*1}
   \multiput(24, 2)( 0,+2){11}{\circle*1}
   \put(0.8,-0.7){\line(+1, 0){22.5}}
   \put(0.8, 0.7){\line(+1, 0){22.5}}
   \put(0.8,23.3){\line(+1, 0){22.5}}
   \put(0.8,24.7){\line(+1, 0){22.5}}
  \end{picture}}
  \;-\;
  \frac1{72}\;
  \raise-8pt\hbox{\begin{picture}(24,24)
   \put( 0, 0){\circle2}
   \put( 0,24){\circle2}
   \put(24, 0){\circle2}
   \put(24,24){\circle2}
   \multiput( 0, 2)( 0,+2){11}{\circle*1}
   \multiput(24, 2)( 0,+2){11}{\circle*1}
   \put( 1.0, 0.0){\line(+1, 0){22}}
   \put( 1.0,24.0){\line(+1, 0){22}}
   \put( 0.7, 0.7){\line(+1,+1){22.5}}
   \put(23.3, 0.7){\line(-1,+1){22.5}}
  \end{picture}}
  \;\;\equiv\;\;
  \frac1{48}\;
  \raise-8pt\hbox{\begin{picture}(24,24)
   \put( 0, 0){\circle2}
   \put( 0,24){\circle2}
   \put(24, 0){\circle2}
   \put(24,24){\circle2}
   \multiput( 0, 2)( 0,+2){11}{\circle*1}
   \multiput(24, 2)( 0,+2){11}{\circle*1}
   \put(0.8,-0.7){\line(+1, 0){22.5}}
   \put(0.8, 0.7){\line(+1, 0){22.5}}
   \put(0.8,23.3){\line(+1, 0){22.5}}
   \put(0.8,24.7){\line(+1, 0){22.5}}
  \end{picture}}
 \end{eqnarray*}
 modulo the ideal $\langle\,\IHX\,\rangle$ based on the congruence (\ref{ihq}).
 Using the Feynman rules (\ref{feyn1}) and (\ref{feyn2}) to work out the
 curvature invariants associated to the left and right hand sides we find
 $$
  \begin{array}{lclcl}
   \left[\;{\displaystyle\frac14}\;
    \raise-13pt\hbox{\begin{picture}(30,30)
     \put( 4,15){\circle4}
     \put(26,15){\circle4}
     \put( 2,15){\line(+1,0){4}}
     \put(24,15){\line(+1,0){4}}
     \qbezier( 2.7,16.4)(-6.0,34.0)(24.7,16.4)
     \qbezier( 2.7,13.6)(-6.0,-4.0)(24.7,13.6)
     \qbezier( 5.3,16.4)(36.0,34.0)(27.3,16.4)
     \qbezier( 5.3,13.6)(36.0,-4.0)(27.3,13.6)
    \end{picture}}\;
   \right](\;R\;)
   &=&
   {\displaystyle\frac14\;\sum_{\mu,\,\nu,\,\alpha,\,\beta\,=\,1}^m
   R(\,E_\mu,E_\nu;\,E_\alpha,E_\beta\,)^2}
   &=&
   \;\;g_{\LP^2T^*\otimes\LP^2T^*}(\;\,R,\;\,R\;\,)
   \\[16pt]
   \left[\;\,{\displaystyle\frac14}\;\,
    \raise-8pt\hbox{\begin{picture}(24,24)
     \put( 0, 0){\circle2}
     \put( 0,24){\circle2}
     \put(24, 0){\circle2}
     \put(24,24){\circle2}
     \multiput( 0, 2)( 0,+2){11}{\circle*1}
     \multiput(24, 2)( 0,+2){11}{\circle*1}
     \put(0.8,-0.7){\line(+1, 0){22.5}}
     \put(0.8, 0.7){\line(+1, 0){22.5}}
     \put(0.8,23.3){\line(+1, 0){22.5}}
     \put(0.8,24.7){\line(+1, 0){22.5}}
    \end{picture}}\;\;\right](\;R\;)
   &=&
   {\displaystyle\frac14\;\sum_{\mu,\,\nu,\,\alpha,\,\beta\,=\,1}^m
   \Sec(\,E_\mu,E_\nu;E_\alpha,E_\beta\,)^2}
   &=&
   g_{\S^2T^*\otimes\S^2T^*}(\,\Sec,\,\Sec\,)
  \end{array}
 $$
 and conclude that the isomorphisms $\Phi^+$ and $\Phi^-$ of Section
 \ref{tensors} are essentially isometries:  
 \begin{equation}\label{eqnorm}
  g_{\S^2T^*\otimes\S^2T^*}(\;\Sec,\;\Sec\;)
  \;\;=\;\;
  12\;g_{\LP^2T^*\otimes\LP^2T^*}(\;R,\;R\;)\ .
 \end{equation}
\section{Explicit Values of Curvature Invariants}
\label{values}
 In order to illustrate the rather abstract construction of stable curvature
 invariants from colored trivalent graphs discussed in Section \ref{graphs}
 we want to present a combinatorial method to calculate explicitly the values
 of these stable curvature invariants on all space form algebraic curvature
 tensors, this is on all algebraic curvature tensors of constant sectional
 curvature. Central to this presentation is a family $\delta_m$ of homogeneous
 derivations of degree $-1$ of the graded graph algebras $\A^\bullet$ and
 $\overline\A^\bullet$ indexed by a formal dimension parameter $m$.

 \pfill
 The space of algebraic curvature tensors on a euclidean vector space $T$
 of dimension $m\,\geq\,2$ contains a unique curvature tensor up to scale
 invariant under the isometry group $\OG(\,T,\,g\,)$, namely the algebraic
 curvature tensor $R^{\OG(\,T,\,g\,)}\,\in\,\Curv^-T$ of constant sectional
 curvature $1$:
 $$
  R^{\OG(\,T,\,g\,)}(\;X,\,Y;\,U,\,V\;)
  \;\;:=\;\;
  -\;g_{\LP^2T}(\;X\wedge Y,\;U\wedge V\;)\ .
 $$
 The corresponding sectional curvature tensor $\Sec^{\OG(\,T,\,g\,)}\,\in\,
 \Curv^+T$ is given by:
 \begin{equation}\label{sec1}
  \Sec^{\OG(\,T,\,g\,)}(X,Y;U,V)
  \;\;=\;\;
  4\,g(X,Y)g(U,V)\,-\,2\,g(X,U)g(Y,V)\,-\,2\,g(X,V)g(Y,U)\ .
 \end{equation}
 The Ricci tensor of the algebraic curvature tensor $R^{\OG(\,T,\,g\,)}$ of
 constant sectional curvature $1$ is the simple multiple $\Ric\,=\,(m-1)\,g$
 of the scalar product $g$ and so its scalar curvature equals $\kappa\,=\,
 m(m-1)$. In passing we remark that $R^{\OG(\,T,\,g\,)}$ can also be defined
 in terms of the Nomizu--Kulkarni product of equation (\ref{nk}), more
 precisely $R^{\OG(\,T,\,g\,)}\,:=\,-\frac12\,g\times g$.

 \begin{Definition}[Curvature Derivation]
 \hfill\label{cd}\break
  The curvature derivation in formal dimension $m\,\in\,\N_0$ is the
  homogeneous derivation
  $$
   \delta_m:\;\;\A^\bullet\;\longrightarrow\;\A^{\bullet-1},
   \qquad[\,\gamma\,]\;\longmapsto\;\delta_m[\,\gamma\,]\ ,
  $$
  of degree $-1$ of the graded algebra $\A^\bullet$ of colored trivalent
  graphs characterized by
  $$
   \Big(\;\delta_m[\,\gamma\,]\;\Big)[\;R\;]
   \;\;=\;\;
   \left.\frac d{dt}\right|_0[\,\gamma\,]
   \Big(\,R\,+\,tR^{\OG(\,T,\,g\,)}\,\Big)
  $$
  for every algebraic curvature tensor $R$ on a euclidean vector space $T$
  of dimension $m$.
 \end{Definition}

 \noindent
 Due to the construction of stable curvature invariants in Definition
 \ref{slci} the value $[\,\gamma\,](\,R\,)$ of such an invariant on an
 algebraic curvature tensor $R$ on a euclidean vector space $T$ is an
 iterated sum over an orthonormal basis of $T$, whose parts are products
 over the red edges of the colored trivalent graph $\gamma$. Replacing
 $R$ by $R+tR^{\OG(\,T,\,g\,)}$ and taking the derivative $\left.\frac d{dt}
 \right|_0$ we thus obtain a sum over the red edges of $\gamma$, in which
 the sectional curvature tensor factor $\Sec\,\in\,\Curv^+T$ associated
 to $R$ has been replaced for this particular red edge by the factor
 $\Sec^{\OG(\,T,\,g\,)}$. In consequence the curvature derivation $\delta_m$
 is the homogeneous derivation of $\A^\bullet$ of degree $-1$, which expands
 each red edge in turn according to equation (\ref{sec1})
 \begin{equation}\label{cdex}
  \raise-12pt\hbox{\begin{picture}(30,30)
   \put(15, 8){\circle2}
   \put(15,22){\circle2}
   \multiput(15,10)( 0,+2){6}{\circle*1}
   \put( 1, 1){\line(+2,+1){13}}
   \put(29, 1){\line(-2,+1){13}}
   \put( 1,29){\line(+2,-1){13}}
   \put(29,29){\line(-2,-1){13}}
  \end{picture}}
  \qquad\stackrel{\delta_m}\longmapsto\qquad
  4\;\raise-12pt\hbox{\begin{picture}(30,30)
   \qbezier( 1, 1)(15,15)(29, 1)
   \qbezier( 1,29)(15,15)(29,29)
  \end{picture}}
  \;-\;
  2\;\raise-12pt\hbox{\begin{picture}(30,30)
   \put( 1, 1){\line(+1,+1){14.2}}
   \put(29,29){\line(-1,-1){12.4}}
   \put(29, 1){\line(-1,+1){12}}
   \put( 1,29){\line(+1,-1){12}}
   \qbezier(13,17)(15,20)(17,13)
  \end{picture}}
  \;-\;
  2\;\raise-12pt\hbox{\begin{picture}(30,30)
   \qbezier( 1, 1)(15,15)( 1,29)
   \qbezier(29, 1)(15,15)(29,29)
  \end{picture}}
 \end{equation}
 and sums all results together with the proviso that every black circle
 without any vertices at all occurred in the process needs to be interpreted as
 the scalar factor $m\,=\,\sum g(dE^\up_\mu,E_\mu)$. Needless to say this is
 the only way in which the dimension $m\,\in\,\N_0$ of the euclidean vector
 space $T$ enters in the definition of the curvature derivation $\delta_m$,
 in other words the derivation $\delta_m$ is at most a quadratic polynomial
 in the formal dimension $m$.

 Evidently every IHX--relation of the form (\ref{ihx}) vanishes under
 the expansion (\ref{cdex}) of the curvature derivation, for this reason
 $\delta_m$ maps the ideal $\langle\,\IHX\,\rangle\,\subseteq\,\A^\bullet$
 of all IHX--relations to itself. In turn the curvature derivation descends
 to a homogeneous derivation of degree $-1$ of the quotient algebra
 $\overline\A^\bullet$ of colored trivalent graphs modulo IHX--relations.
 Calculating $\delta_m$ for all the generators (\ref{gen6}) of the reduced
 graph algebra $\overline\A^\bullet$ up to degree $3$ we find easily
 \begin{equation}\label{cdv1}
  \delta_m\,\raise-9pt\hbox{\begin{picture}(10,26)
   \put( 5, 2){\circle2}
   \put( 5,24){\circle2}
   \qbezier(4.0, 2.0)( 2.0, 2.0)( 2.0,12.0)
   \qbezier(4.0,24.0)( 2.0,24.0)( 2.0,12.0)
   \qbezier(6.1, 2.0)( 8.0, 2.0)( 8.0,12.0)
   \qbezier(6.1,24.0)( 8.0,24.0)( 8.0,12.0)
   \multiput( 5.0, 4.0)( 0.0,+2.0){10}{\circle*1}
  \end{picture}}
  \;\;=\;\;
  -\,2\,m(m-1)
  \qquad
  \delta_m\,\raise-9pt\hbox{\begin{picture}(26,26)
   \put( 2, 2){\circle2}
   \put(24, 2){\circle2}
   \put( 2,24){\circle2}
   \put(24,24){\circle2}
   \put( 2.0, 3.0){\line( 0,+1){20}}
   \put(24.0, 3.0){\line( 0,+1){20}}
   \put( 3.0, 2.0){\line(+1, 0){20}}
   \put( 3.0,24.0){\line(+1, 0){20}}
   \multiput( 3.2, 3.9)( 0.0,+2.0){10}{\circle*1}
   \multiput(22.8, 3.9)( 0.0,+2.0){10}{\circle*1}
  \end{picture}}
  \;\;=\;\;
  -\,4\,(m-1)\,\raise-9pt\hbox{\begin{picture}(10,26)
   \put( 5, 2){\circle2}
   \put( 5,24){\circle2}
   \qbezier(4.0, 2.0)( 2.0, 2.0)( 2.0,12.0)
   \qbezier(4.0,24.0)( 2.0,24.0)( 2.0,12.0)
   \qbezier(6.1, 2.0)( 8.0, 2.0)( 8.0,12.0)
   \qbezier(6.1,24.0)( 8.0,24.0)( 8.0,12.0)
   \multiput( 5.0, 4.0)( 0.0,+2.0){10}{\circle*1}
  \end{picture}}
  \qquad
  \delta_m\,\raise-9pt\hbox{\begin{picture}(26,26)
   \put( 2, 2){\circle2}
   \put(24, 2){\circle2}
   \put( 2,24){\circle2}
   \put(24,24){\circle2}
   \put( 2.0, 3.0){\line( 0,+1){20}}
   \put(24.0, 3.0){\line( 0,+1){20}}
   \put( 3.0, 2.0){\line(+1, 0){20}}
   \put( 3.0,24.0){\line(+1, 0){20}}
   \multiput( 3.4, 3.4)(+1.2,+1.2){17}{\circle*1}
   \multiput(22.6, 3.4)(-1.2,+1.2){17}{\circle*1}
  \end{picture}}
  \;\;=\;\;
  +\,12\,\raise-9pt\hbox{\begin{picture}(10,26)
   \put( 5, 2){\circle2}
   \put( 5,24){\circle2}
   \qbezier(4.0, 2.0)( 2.0, 2.0)( 2.0,12.0)
   \qbezier(4.0,24.0)( 2.0,24.0)( 2.0,12.0)
   \qbezier(6.1, 2.0)( 8.0, 2.0)( 8.0,12.0)
   \qbezier(6.1,24.0)( 8.0,24.0)( 8.0,12.0)
   \multiput( 5.0, 4.0)( 0.0,+2.0){10}{\circle*1}
  \end{picture}}
 \end{equation}
 for the generators of degree $1$ and $2$ as well as:
 \begin{equation}\label{cdv}
  \begin{array}{lclcllcl}
   \delta_m\;\raise-11pt\hbox{\begin{picture}(30,30)
    \put(15, 2){\circle2}
    \put(15,28){\circle2}
    \put( 3, 8){\circle2}
    \put(27, 8){\circle2}
    \put( 3,22){\circle2}
    \put(27,22){\circle2}
    \put( 3.0, 9.0){\line( 0,+1){12}}
    \put(27.0, 9.0){\line( 0,+1){12}}
    \put(15.9, 2.5){\line(+2,+1){10}}
    \put(14.1, 2.5){\line(-2,+1){10}}
    \put(15.9,27.5){\line(+2,-1){10}}
    \put(14.1,27.5){\line(-2,-1){10}}
    \multiput( 4.9,21.9)(+1.5,+0.75){7}{\circle*1}
    \multiput( 4.9, 8.1)(+1.5,-0.75){7}{\circle*1}
    \multiput(26.1,10.1)( 0.0,+2.0){6}{\circle*1}
   \end{picture}}
   &=&
   -&\!\!\!6\,(m-1)&
   \!\!\!\raise-9pt\hbox{\begin{picture}(26,26)
    \put( 2, 2){\circle2}
    \put(24, 2){\circle2}
    \put( 2,24){\circle2}
    \put(24,24){\circle2}
    \put( 2.0, 3.0){\line( 0,+1){20}}
    \put(24.0, 3.0){\line( 0,+1){20}}
    \put( 3.0, 2.0){\line(+1, 0){20}}
    \put( 3.0,24.0){\line(+1, 0){20}}
    \multiput( 3.2, 3.9)( 0.0,+2.0){10}{\circle*1}
    \multiput(22.8, 3.9)( 0.0,+2.0){10}{\circle*1}
   \end{picture}}
   &&&
   \\[7pt]
   \delta_m\;\raise-11pt\hbox{\begin{picture}(30,30)
    \put(15, 2){\circle2}
    \put(15,28){\circle2}
    \put( 3, 8){\circle2}
    \put(27, 8){\circle2}
    \put( 3,22){\circle2}
    \put(27,22){\circle2}
    \put( 3.0, 9.0){\line( 0,+1){12}}
    \put(27.0, 9.0){\line( 0,+1){12}}
    \put(15.9, 2.5){\line(+2,+1){10}}
    \put(14.1, 2.5){\line(-2,+1){10}}
    \put(15.9,27.5){\line(+2,-1){10}}
    \put(14.1,27.5){\line(-2,-1){10}}
    \multiput( 3.9,10.1)( 0.0,+2.0){6}{\circle*1}
    \multiput(15.0, 4.0)( 0.0,+2.0){12}{\circle*1}
    \multiput(26.1,10.1)( 0.0,+2.0){6}{\circle*1}
   \end{picture}}
   &=&
   -&\!\!\!(\,4m-6\,)&
   \!\!\!\raise-9pt\hbox{\begin{picture}(26,26)
    \put( 2, 2){\circle2}
    \put(24, 2){\circle2}
    \put( 2,24){\circle2}
    \put(24,24){\circle2}
    \put( 2.0, 3.0){\line( 0,+1){20}}
    \put(24.0, 3.0){\line( 0,+1){20}}
    \put( 3.0, 2.0){\line(+1, 0){20}}
    \put( 3.0,24.0){\line(+1, 0){20}}
    \multiput( 3.2, 3.9)( 0.0,+2.0){10}{\circle*1}
    \multiput(22.8, 3.9)( 0.0,+2.0){10}{\circle*1}
   \end{picture}}
   &-&2&
   \!\!\!\raise-9pt\hbox{\begin{picture}(26,26)
   \put( 6, 2){\circle2}
   \put( 6,24){\circle2}
   \put(20, 2){\circle2}
   \put(20,24){\circle2}
   \qbezier( 5.0, 2.0)( 3.0, 2.0)( 3.0,12.0)
   \qbezier( 5.0,24.0)( 3.0,24.0)( 3.0,12.0)
   \qbezier( 7.1, 2.0)( 9.0, 2.0)( 9.0,12.0)
   \qbezier( 7.1,24.0)( 9.0,24.0)( 9.0,12.0)
   \qbezier(19.0, 2.0)(17.0, 2.0)(17.0,12.0)
   \qbezier(19.0,24.0)(17.0,24.0)(17.0,12.0)
   \qbezier(21.1, 2.0)(23.0, 2.0)(23.0,12.0)
   \qbezier(21.1,24.0)(23.0,24.0)(23.0,12.0)
   \multiput( 6.0, 4.0)( 0.0,+2.0){10}{\circle*1}
   \multiput(20.0, 4.0)( 0.0,+2.0){10}{\circle*1}
  \end{picture}}
   \\[7pt]
   \delta_m\;\raise-11pt\hbox{\begin{picture}(30,30)
    \put(15, 2){\circle2}
    \put(15,28){\circle2}
    \put( 3, 8){\circle2}
    \put(27, 8){\circle2}
    \put( 3,22){\circle2}
    \put(27,22){\circle2}
    \put( 3.0, 9.0){\line( 0,+1){12}}
    \put(27.0, 9.0){\line( 0,+1){12}}
    \put(15.9, 2.5){\line(+2,+1){10}}
    \put(14.1, 2.5){\line(-2,+1){10}}
    \put(15.9,27.5){\line(+2,-1){10}}
    \put(14.1,27.5){\line(-2,-1){10}}
    \multiput( 4.9,21.9)(+1.5,+0.75){7}{\circle*1}
    \multiput(16.1, 3.8)(+0.9,+1.5){12}{\circle*1}
    \multiput( 5.0, 8.0)(+2.0, 0.0){11}{\circle*1}
   \end{picture}}
   &=&
   +&12&
   \!\!\!\raise-9pt\hbox{\begin{picture}(26,26)
    \put( 2, 2){\circle2}
    \put(24, 2){\circle2}
    \put( 2,24){\circle2}
    \put(24,24){\circle2}
    \put( 2.0, 3.0){\line( 0,+1){20}}
    \put(24.0, 3.0){\line( 0,+1){20}}
    \put( 3.0, 2.0){\line(+1, 0){20}}
    \put( 3.0,24.0){\line(+1, 0){20}}
    \multiput( 3.2, 3.9)( 0.0,+2.0){10}{\circle*1}
    \multiput(22.8, 3.9)( 0.0,+2.0){10}{\circle*1}
   \end{picture}}
   &-&\!\!\!2\,(m-1)&
   \!\!\!\raise-9pt\hbox{\begin{picture}(26,26)
    \put( 2, 2){\circle2}
    \put(24, 2){\circle2}
    \put( 2,24){\circle2}
    \put(24,24){\circle2}
    \put( 2.0, 3.0){\line( 0,+1){20}}
    \put(24.0, 3.0){\line( 0,+1){20}}
    \put( 3.0, 2.0){\line(+1, 0){20}}
    \put( 3.0,24.0){\line(+1, 0){20}}
    \multiput( 3.4, 3.4)(+1.2,+1.2){17}{\circle*1}
    \multiput(22.6, 3.4)(-1.2,+1.2){17}{\circle*1}
   \end{picture}}
   \\[7pt]
   \delta_m\;\raise-11pt\hbox{\begin{picture}(30,30)
    \put(15, 2){\circle2}
    \put(15,28){\circle2}
    \put( 3, 8){\circle2}
    \put(27, 8){\circle2}
    \put( 3,22){\circle2}
    \put(27,22){\circle2}
    \put( 3.0, 9.0){\line( 0,+1){12}}
    \put(27.0, 9.0){\line( 0,+1){12}}
    \put(15.9, 2.5){\line(+2,+1){10}}
    \put(14.1, 2.5){\line(-2,+1){10}}
    \put(15.9,27.5){\line(+2,-1){10}}
    \put(14.1,27.5){\line(-2,-1){10}}
    \multiput(15.0, 4.0)( 0.0,+2.0){12}{\circle*1}
    \multiput( 4.9, 8.2)(+2.0, 0.0){11}{\circle*1}
    \multiput( 4.9,21.8)(+2.0, 0.0){11}{\circle*1}
   \end{picture}}
   &=&
   +&12&
   \!\!\!\raise-9pt\hbox{\begin{picture}(26,26)
    \put( 2, 2){\circle2}
    \put(24, 2){\circle2}
    \put( 2,24){\circle2}
    \put(24,24){\circle2}
    \put( 2.0, 3.0){\line( 0,+1){20}}
    \put(24.0, 3.0){\line( 0,+1){20}}
    \put( 3.0, 2.0){\line(+1, 0){20}}
    \put( 3.0,24.0){\line(+1, 0){20}}
    \multiput( 3.2, 3.9)( 0.0,+2.0){10}{\circle*1}
    \multiput(22.8, 3.9)( 0.0,+2.0){10}{\circle*1}
   \end{picture}}
   &+&6&
   \!\!\!\raise-9pt\hbox{\begin{picture}(26,26)
    \put( 2, 2){\circle2}
    \put(24, 2){\circle2}
    \put( 2,24){\circle2}
    \put(24,24){\circle2}
    \put( 2.0, 3.0){\line( 0,+1){20}}
    \put(24.0, 3.0){\line( 0,+1){20}}
    \put( 3.0, 2.0){\line(+1, 0){20}}
    \put( 3.0,24.0){\line(+1, 0){20}}
    \multiput( 3.4, 3.4)(+1.2,+1.2){17}{\circle*1}
    \multiput(22.6, 3.4)(-1.2,+1.2){17}{\circle*1}
   \end{picture}}
   \\[7pt]
   \delta_m\;\raise-11pt\hbox{\begin{picture}(30,30)
    \put(15, 2){\circle2}
    \put(15,28){\circle2}
    \put( 3, 8){\circle2}
    \put(27, 8){\circle2}
    \put( 3,22){\circle2}
    \put(27,22){\circle2}
    \put( 3.0, 9.0){\line( 0,+1){12}}
    \put(27.0, 9.0){\line( 0,+1){12}}
    \put(15.9, 2.5){\line(+2,+1){10}}
    \put(14.1, 2.5){\line(-2,+1){10}}
    \put(15.9,27.5){\line(+2,-1){10}}
    \put(14.1,27.5){\line(-2,-1){10}}
    \multiput(15.0, 4.0)( 0.0,+2.0){12}{\circle*1}
    \multiput( 4.7, 9.0)(+1.9,+1.1){12}{\circle*1}
    \multiput(25.3, 9.0)(-1.9,+1.1){12}{\circle*1}
   \end{picture}}
   &=&
   -&6&
   \!\!\!\raise-9pt\hbox{\begin{picture}(26,26)
    \put( 2, 2){\circle2}
    \put(24, 2){\circle2}
    \put( 2,24){\circle2}
    \put(24,24){\circle2}
    \put( 2.0, 3.0){\line( 0,+1){20}}
    \put(24.0, 3.0){\line( 0,+1){20}}
    \put( 3.0, 2.0){\line(+1, 0){20}}
    \put( 3.0,24.0){\line(+1, 0){20}}
    \multiput( 3.2, 3.9)( 0.0,+2.0){10}{\circle*1}
    \multiput(22.8, 3.9)( 0.0,+2.0){10}{\circle*1}
   \end{picture}}
   &+&24&
   \!\!\!\raise-9pt\hbox{\begin{picture}(26,26)
    \put( 2, 2){\circle2}
    \put(24, 2){\circle2}
    \put( 2,24){\circle2}
    \put(24,24){\circle2}
    \put( 2.0, 3.0){\line( 0,+1){20}}
    \put(24.0, 3.0){\line( 0,+1){20}}
    \put( 3.0, 2.0){\line(+1, 0){20}}
    \put( 3.0,24.0){\line(+1, 0){20}}
    \multiput( 3.4, 3.4)(+1.2,+1.2){17}{\circle*1}
    \multiput(22.6, 3.4)(-1.2,+1.2){17}{\circle*1}
   \end{picture}}\ .
  \end{array}
 \end{equation}
 It should be pointed out that the curvature derivation allows us to calculate
 explicitly the values of stable curvature invariants on algebraic curvature
 tensors $R\,\in\,\Curv^-T$ of constant sectional curvature on every euclidean
 vector space $T$ of dimension $m\,\geq\,2$ using the formula
 \begin{equation}\label{val}
  [\,\gamma\,](\;R\;)
  \;\;=\;\;
  \frac1{n!}\,\Big(\frac\kappa{m(m-1)}\Big)^n\,\delta^n_m[\,\gamma\,]
  \;\;\in\;\;
  \overline\A^0
  \;\;=\;\;
  \R
 \end{equation}
 valid for all isomorphism classes $[\,\gamma\,]\,\in\,\Gamma^n$ of colored
 trivalent graphs $\gamma$ of degree $n\,\in\,\N_0$, where $\kappa$
 denotes the scalar curvature of $R$. In fact this formula is a direct
 consequence of the homogeneity of degree $n$ of the polynomial $[\,\gamma\,]$
 and the ensuing recursion relation:
 $$
  \Big(\;\delta_m[\,\gamma\,]\;\Big)(\;R^{\OG(\,T,\,g\,)}\;)
  \;\;=\;\;
  \left.\frac d{dt}\right|_0(\,1+t\,)^n\,[\,\gamma\,](\;R^{\OG(\,T,\,g\,)}\;)
  \;\;=\;\;
  n\,[\,\gamma\,](\;R^{\OG(\,T,\,g\,)}\;)\ .
 $$
 Based on equations (\ref{cdv1}) and (\ref{cdv}) we find say for every
 space form curvature tensor $R$:
 $$
  \left[\raise-11pt\hbox{\begin{picture}(30,30)
   \put(15, 2){\circle2}
   \put(15,28){\circle2}
   \put( 3, 8){\circle2}
   \put(27, 8){\circle2}
   \put( 3,22){\circle2}
   \put(27,22){\circle2}
   \put( 3.0, 9.0){\line( 0,+1){12}}
   \put(27.0, 9.0){\line( 0,+1){12}}
   \put(15.9, 2.5){\line(+2,+1){10}}
   \put(14.1, 2.5){\line(-2,+1){10}}
   \put(15.9,27.5){\line(+2,-1){10}}
   \put(14.1,27.5){\line(-2,-1){10}}
   \multiput(15.0, 4.0)( 0.0,+2.0){12}{\circle*1}
   \multiput( 4.9, 8.2)(+2.0, 0.0){11}{\circle*1}
   \multiput( 4.9,21.8)(+2.0, 0.0){11}{\circle*1}
  \end{picture}}\right](\;R\;)
  \;\;=\;\;
  +\;8\,\frac{(2m-5)\,\kappa^3}{m^2\,(m-1)^2}
  \qquad\qquad
  \left[\raise-11pt\hbox{\begin{picture}(30,30)
   \put(15, 2){\circle2}
   \put(15,28){\circle2}
   \put( 3, 8){\circle2}
   \put(27, 8){\circle2}
   \put( 3,22){\circle2}
   \put(27,22){\circle2}
   \put( 3.0, 9.0){\line( 0,+1){12}}
   \put(27.0, 9.0){\line( 0,+1){12}}
   \put(15.9, 2.5){\line(+2,+1){10}}
   \put(14.1, 2.5){\line(-2,+1){10}}
   \put(15.9,27.5){\line(+2,-1){10}}
   \put(14.1,27.5){\line(-2,-1){10}}
   \multiput(15.0, 4.0)( 0.0,+2.0){12}{\circle*1}
   \multiput( 4.7, 9.0)(+1.9,+1.1){12}{\circle*1}
   \multiput(25.3, 9.0)(-1.9,+1.1){12}{\circle*1}
  \end{picture}}\right](\;R\;)
  \;\;=\;\;
  -\;8\,\frac{(m+11)\,\kappa^3}{m^2\,(m-1)^2}\ .
 $$
\section{Pfaffian and Moment Polynomials}
\label{special}
 All characteristic numbers of a compact Riemannian manifold are integrated
 polynomial invariants of the curvature tensor, the corresponding polynomial
 however is not a stable curvature invariant in the sense of this article.
 The only stable characteristic number among the classical Euler and Pontryagin
 numbers turns out to be the Euler characteristic, in this section we will
 identity the corresponding elements in the algebra $\A^\bullet$ of colored
 trivalent graphs, the Pfaffian polynomials $(\,\pf_n\,)_{n\,\in\,\N_0}$.
 Besides the Euler characteristic we will study the normalized moment
 polynomials $(\,\Psi^\circ_n\,)_{n\,\in\,\N_0}$ in this section, which
 calculate the moments $(\,\Psi_n\,)_{n\,\in\,\N_0}$ of the sectional
 curvature considered as a random variable on the Graßmannian of planes
 up to a normalization constant depending on the dimension.

 \pfill
 The sequence of Pfaffian polynomials is a sequence $(\,\pf_n\,)_{n\,\in\,
 \N_0}$ of elements of the algebra $\A^\bullet$ of colored trivalent graphs
 defined for $n\,=\,0$ by $\pf_0\,=\,\1$ and in positive degrees $n\,>\,0$ by:
 $$
  \pf_n(\;R\;)
  \;\;:=\;\;
  \frac1{12^n\,n!}\;\sum_{\sigma\,\in\,S_{2n}}(\,\sgn\,\sigma\,)
  \sum_{\mu_1,\,\ldots,\,\mu_{2n}\,=\,1}^m\;\prod_{r\,=\,1}^n\;
  \Sec(\;E_{\mu_{2r-1}},\;dE^\up_{\mu_{\sigma(2r-1)}};\;E_{\mu_{2r}},
  \;dE^\up_{\mu_{\sigma(2r)}}\,)\ .
 $$
 Right from this definition we can read off the corresponding element
 of the graph algebra $\A^\bullet$ 
 \begin{equation}\label{pfg}
  \pf_n
  \;\;=\;\;
  \frac1{12^n\,n!}
  \sum_{\sigma\,\in\,S_{2n}}
  (\,\sgn\;\sigma\,)\;
  \underbrace{\raise-12pt\hbox{\begin{picture}(234,30)
   \put(  4,15){\circle2}
   \put( 30,15){\circle2}
   \put( 54,15){\circle2}
   \put( 80,15){\circle2}
   \put(104,15){\circle2}
   \put(130,15){\circle2}
   \put(204,15){\circle2}
   \put(230,15){\circle2}
   \multiput(  6,15)(+2, 0){12}{\circle*1}
   \multiput( 56,15)(+2, 0){12}{\circle*1}
   \multiput(106,15)(+2, 0){12}{\circle*1}
   \multiput(206,15)(+2, 0){12}{\circle*1}
   \put(  4,16){\line( 0,+1){8}}
   \put(  4,14){\line( 0,-1){8}}
   \put( 30,16){\line( 0,+1){8}}
   \put( 30,14){\line( 0,-1){8}}
   \put( 54,16){\line( 0,+1){8}}
   \put( 54,14){\line( 0,-1){8}}
   \put( 80,16){\line( 0,+1){8}}
   \put( 80,14){\line( 0,-1){8}}
   \put(104,16){\line( 0,+1){8}}
   \put(104,14){\line( 0,-1){8}}
   \put(130,16){\line( 0,+1){8}}
   \put(130,14){\line( 0,-1){8}}
   \put(204,16){\line( 0,+1){8}}
   \put(204,14){\line( 0,-1){8}}
   \put(230,16){\line( 0,+1){8}}
   \put(230,14){\line( 0,-1){8}}
   \put(160,15){$\ldots$}
   \put(  2, 0){$\scriptscriptstyle1$}
   \put( -2,27){$\scriptscriptstyle\sigma(1)$}
   \put( 28, 0){$\scriptscriptstyle2$}
   \put( 24,27){$\scriptscriptstyle\sigma(2)$}
   \put( 52, 0){$\scriptscriptstyle3$}
   \put( 48,27){$\scriptscriptstyle\sigma(3)$}
   \put( 78, 0){$\scriptscriptstyle4$}
   \put( 74,27){$\scriptscriptstyle\sigma(4)$}
   \put(102, 0){$\scriptscriptstyle5$}
   \put( 98,27){$\scriptscriptstyle\sigma(5)$}
   \put(128, 0){$\scriptscriptstyle6$}
   \put(124,27){$\scriptscriptstyle\sigma(6)$}
   \put(195, 0){$\scriptscriptstyle2n-1$}
   \put(189,27){$\scriptscriptstyle\sigma(2n-1)$}
   \put(226, 0){$\scriptscriptstyle2n$}   
   \put(223,27){$\scriptscriptstyle\sigma(2n)$}
  \end{picture}}}_n\ ,
 \end{equation}
 where it is understood that a black edge runs between every pair of flags
 indexed by the same integer, one flag above and one below. According to
 the Theorem of Chern--Gauß--Bonnet the Euler characteristic of every
 compact Riemannian manifold $M$ of even dimension $m\,\in\,2\N_0$ with
 Riemannian metric $g$ can be written as the integrated curvature invariant
 \begin{equation}\label{cgb}
  \chi(\;M\;)
  \;\;=\;\;
  \frac1{(\,2\pi\,)^{\frac m2}}\;
  \int_M\pf_{\frac m2}(\;R\;)\,|\,\vol_g\,|\ .
 \end{equation}
 associated to the sequence $(\,\pf_n\,)_{n\,\in\,\N_0}$. For Riemannian
 surfaces $M$ for example we find
 $$
  \pf_1
  \;\;=\;\;
  \frac1{12}\;\Big(\;+\;
   \raise-9pt\hbox{\begin{picture}(40,24)(0,0)
    \put( 7,12){\circle2}
    \put(33,12){\circle2}
    \multiput( 9,12)(+2, 0){12}{\circle*1}
    \put( 7,11){\line( 0,-1){8}}
    \put(33,11){\line( 0,-1){8}}
    \put( 7,13){\line( 0,+1){8}}
    \put(33,13){\line( 0,+1){8}}
    \qbezier( 7,21)( 3,24)( 2,17)
    \qbezier( 2,17)( 1,12)( 2, 7)
    \qbezier( 7, 3)( 3, 0)( 2, 7)
    \qbezier(33,21)(37,24)(38,17)
    \qbezier(38,17)(39,12)(38, 7)
    \qbezier(33, 3)(37, 0)(38, 7)
   \end{picture}}
   \;-\;
   \raise-9pt\hbox{\begin{picture}(30,24)(0,0)
    \put( 2,12){\circle2}
    \put(28,12){\circle2}
    \multiput( 4,12)(+2, 0){12}{\circle*1}
    \put( 2,11){\line( 0,-1){8}}
    \put(28,11){\line( 0,-1){8}}
    \put( 2,13){\line( 0,+1){8}}
    \put(28,13){\line( 0,+1){8}}
    \qbezier( 2,21)(10,26)(15,12)
    \qbezier( 2, 3)(10,-2)(15,12)
    \qbezier(28,21)(20,26)(15,12)
    \qbezier(28, 3)(20,-2)(15,12)
   \end{picture}}
  \;\;\Big)
  \;\;\equiv\;\;
  -\;\frac14\;
  \raise-9pt\hbox{\begin{picture}(30,24)(0,0)
   \put( 2,12){\circle2}
   \put(28,12){\circle2}
   \multiput( 4,12)(+2, 0){12}{\circle*1}
   \qbezier( 2.7,11.3)( 8.0, 6.0)(15.0, 6.0)
   \qbezier(27.3,11.3)(22.0, 6.0)(15.0, 6.0)
   \qbezier( 2.7,12.7)( 8.0,18.0)(15.0,18.0)
   \qbezier(27.3,12.7)(22.0,18.0)(15.0,18.0)
  \end{picture}}
  \;\;\widehat=\;\;
  \frac\kappa2
 $$
 using the basic IHX--relation (\ref{ihex}) and conclude that $\pf_1(\,R\,)$
 equals the Gaußian curvature of the surface $M$; keep in mind that $\kappa$
 denotes the scalar curvature throughout this article. In order to relate
 equation (\ref{cgb}) to the Theorem of Chern--Gauß--Bonnet in higher
 dimensions we recall that it involves in its standard formulation the
 pointwise Berezin integral
 $$
  [\;\;]_{g,\,o}:\;\;\LP^\bullet T^*M\,\otimes\,\LP^\circ TM
  \;\longrightarrow\;\LP^\bullet T^*M,\qquad\xi\,\otimes\,\X
  \;\longmapsto\;\<\,\vol_{g,\,o},\,\X\,>\,\xi\ ,
 $$
 which depends on both the Riemannian metric $g$ and an additional orientation
 $o$ via the Riemannian volume form $\vol_{g,\,o}\,\in\,\Gamma(\,\LP^mT^*M\,)$.
 Considering the curvature as a bivector valued $2$--form $R\,\in\,\Gamma(\,
 \LP^2T^*M\otimes\LP^2TM\,)$ like in equation (\ref{2v}) we may exponentiate it
 pointwise in the algebra bundle $\LP\,T^*M\otimes\LP\,TM$ and take the full
 Berezin integral of the result to obtain
 \begin{equation}\label{chern}
  \chi(\;M\;)
  \;\;=\;\;
  \int_{(\,M,\,o\,)}
  \Big[\;\exp\Big(\;-\,\frac R{2\pi}\;\Big)\;\Big]_{g,\,o}
 \end{equation}
 according to Chern (\cite{ch},\cite{bgv},\cite{wa}). One way to define
 the oriented integration of the Pfaffian differential form $\mathrm{Pf}
 (\,\frac R{2\pi}\,)\,:=\,[\,\exp(\,-\frac R{2\pi}\,)\,]_{g,\,o}$ over the
 oriented manifold $M$ is to think of the orientation as the real line bundle
 isomorphism $o:\,\LP^mT^*M\longrightarrow\vartheta^1M$ determined by
 $o(\,\vol_{g,\,o}\,)\,=\,|\,\vol_g\,|$ to turn the Pfaffian differential
 form into a multiple of the volume density $|\,\vol_g\,|\,\in\,
 \Gamma(\,\vartheta^1M\,)$. Combining the pointwise Berezin integral
 $[\;]_{g,\,o}$ used to define $\mathrm{Pf}(\,\frac R{2\pi}\,)$ with
 this orientation vector bundle isomorphism we obtain the vector bundle
 homomorphism
 $$
  [\;\;]_g:\;\;\LP^\bullet T^*M\,\otimes\,\LP^\circ TM
  \;\stackrel{[\;]_{g,\,o}}\longrightarrow\;\LP^\bullet T^*M
  \;\stackrel o\longrightarrow\;\vartheta^1M,\qquad
  \xi\,\otimes\,\X\;\longmapsto\;
  \<\,\xi,\,\pr_{\LP^mTM}\X\,>\,|\,\vol_g\,|\ ,
 $$
 where $\pr_{\LP^mTM}:\,\LP^\circ TM\longrightarrow\LP^mTM$ denotes the
 projection to the top dimensional term. This vector bundle homomorphism
 however is defined even for non--orientable manifolds $M$ and so we can
 rewrite the Theorem of Chern--Gauß--Bonnet as an unoriented integral identity
 $$
  \chi(\;M\;)
  \;\;=\;\;
  \frac1{(\,2\pi\,)^{\frac m2}}\;\int_M\Big[\;\frac 1{(\frac m2)!}
  \,\underbrace{(-R)\,\wedge\,\ldots\,\wedge\,(-R)}_{\frac m2}\;\Big]_g\ ,
 $$
 in fact the projection of $\exp(\,-\frac R{2\pi}\,)$ to the top dimensional
 term leaves us with the $\frac m2$--th power of $R$ only. Replacing the
 algebraic curvature tensor $R\,=\,\Phi^-\Sec$ in its expansion (\ref{2v})
 as a bivector valued $2$--form by the corresponding sectional curvature
 tensor $\Sec$ we find 
 \begin{eqnarray*}
  \lefteqn{-\;R}
  &&
  \\
  &=&
  +\;\frac1{24}\,\sum_{\mu,\,\nu,\,\alpha,\,\beta\,=\,1}^m
  \Big(\,\Sec(E_\mu,E_\alpha;E_\nu,E_\beta)\,-\,
  \Sec(E_\mu,E_\beta;E_\nu,E_\alpha)\,\Big)
  dE_\mu\wedge dE_\nu\,\otimes\,dE^\up_\alpha\wedge dE^\up_\beta
  \\
  &=&
  +\;\frac1{12}\,\sum_{{\scriptstyle\mu_1,\,\mu_2\,=\,1}
   \atop{\scriptstyle\alpha_1,\,\alpha_2\,=\,1}}^m
  \Sec(\,E_{\mu_1},\,E_{\alpha_1};\,E_{\mu_2},\,E_{\alpha_2}\,)\;
  dE_{\mu_1}\wedge dE_{\mu_2}\,\otimes\,dE^\up_{\alpha_1}\wedge
  dE^\up_{\alpha_2}
 \end{eqnarray*}
 and conclude that $\frac1{n!}(-R)^n\,\in\,\Gamma(\,\LP^{2n}T^*M\otimes
 \LP^{2n}TM\,)$ can be expanded for all $n\,\in\,\N_0$ into:
 \begin{eqnarray*}
  \frac1{n!}\underbrace{(-R)\wedge\ldots\wedge(-R)}_n
  &=&
  \frac1{12^n\,n!}
  \sum_{{\scriptstyle\mu_1,\,\ldots,\,\mu_{2n}\,=\,1}\atop
   {\scriptstyle\alpha_1,\,\ldots,\,\alpha_{2n}\,=\,1}}^m
  \prod_{r\,=\,1}^n
  \Sec(E_{\mu_{2r-1}},E_{\alpha_{2r-1}};E_{\mu_{2r}},E_{\alpha_{2r}})
  \\[-12pt]
  &&
  \qquad\qquad\qquad\qquad\qquad
  dE_{\mu_1}\wedge\ldots\wedge dE_{\mu_{2n}}\otimes
  dE^\up_{\alpha_1}\wedge\ldots\wedge dE^\up_{\alpha_{2n}}\ .
 \end{eqnarray*}
 In order to write the Theorem of Chern--Gauß--Bonnet in the form (\ref{cgb})
 it remains to observe
 $$
  \<dE_{\mu_1}\wedge\ldots\wedge dE_{\mu_{2n}},
  dE^\up_{\alpha_1}\wedge\ldots\wedge dE^\up_{\alpha_{2n}}>
  \;\;=\;\;
  \sum_{\sigma\,\in\,S_{2n}}(\,\sgn\;\sigma\,)
  \,dE_{\mu_{\sigma(1)}}(\,dE^\up_{\alpha_1}\,)\,\ldots\,
  \,dE_{\mu_{\sigma(2n)}}(\,dE^\up_{\alpha_{2n}}\,)
 $$
 and sum over $\alpha_1,\,\ldots,\,\alpha_{2n}\,=\,1,\,\ldots,\,m$ using
 the equality $\xi^\up\,=\,\sum_\alpha\xi(\,dE^\up_\alpha\,)\,E_\alpha$.
 Although the argument presented above establishes the equivalence of the
 versions (\ref{cgb}) and (\ref{chern}) of the Theorem of Chern--Gauß--Bonnet
 for all orientable compact manifolds, it does not per se prove the validity
 of the unoriented version (\ref{cgb}) for all compact manifolds. In case of
 doubts the reader may simply apply equation (\ref{cgb}) to the orientable
 $2$--fold cover of a compact manifold $M$ and divide both sides by $2$.
 Better still would be to rework Flanders' proof \cite{fl} of the Theorem of
 Chern--Gauß--Bonnet with densities instead of differential forms.
 
 \begin{Lemma}[Generating Series of Pfaffian Polynomials]
 \hfill\label{efp}\break
  Every isomorphism class of colored trivalent graphs occurs with non--zero
  coefficient in the Pfaffian $\pf_n\,\in\,\A^n$ of the appropriate degree
  $n\,\in\,\N_0$. In the power series completion of the graded algebra
  $\A^\bullet$ of colored trivalent graphs the total Pfaffian equals
  the exponential
  $$
   \sum_{n\,\geq\,0}\pf_n
   \;\;=\;\;
   \exp\left(\;\sum_{[\,\gamma\,]\,\in\,\Gamma^\bullet_\conn}
   \frac{(-1)^{e(\gamma)}}{6^{n(\gamma)}}\;
   \frac{2^{g(\gamma)}}{\#\overline\Aut\;\gamma}\;
   [\,\gamma\,]\;\right)\ ,
  $$
  where $n(\gamma)\,:=\,\frac12\,\#\Vert\,\gamma$ equals the degree of
  $\gamma$, while $e(\gamma)$ and $g(\gamma)$ are the numbers of cycles
  of the black subgraph $\gamma_\black$ of $\gamma$ of even and of length
  greater than $2$ respectively.
 \end{Lemma}

 \proof
 Consider to begin some $n\,\in\,\N$ and a non--void colored trivalent
 graph $\gamma$ with $2n$ vertices. Choosing a direction of cyclically
 traversing each cycle in the associated bivalent black subgraph
 $\gamma_\black$ allows us to define a successor permutation $\sigma\,\in\,
 S_{\Vert\,\gamma}$ of the set of vertices of $\gamma$, which sends every
 vertex to the next vertex on the same cycle in the chosen direction. The
 signature $\sgn\,\sigma\,=\,(-1)^{e(\gamma)}$ of the successor permutation
 reflects the parity of the number $e(\gamma)\,\in\,\N_0$ of even length
 cycles as always with permutations. In addition we chose a labelling
 $L:\,\Vert\,\gamma\longrightarrow\{\,1,\ldots,2n\,\}$ of the vertices
 of $\gamma$ such that red edges run exactly between the pairs $\{L^{-1}(1),
 L^{-1}(2)\},\,\{L^{-1}(3),L^{-1}(4)\},\,\ldots,\,\{L^{-1}(2n-1),L^{-1}(2n)\}$
 of vertices. With such a choice of labelling the original colored trivalent
 graph $\gamma$ turns out to be isomorphic to the summand graph in the sum
 (\ref{pfg}) associated to the permutation $L\circ\sigma\circ L^{-1}$ of the
 index set $\{\,1,\ldots,2n\,\}$ of the same signature.

 Each summand graph $\hat\gamma$ in the sum (\ref{pfg}) on the
 other hand comes along with a tautological labelling $L:\,\Vert\,\hat\gamma
 \longrightarrow\{\,1,\ldots,2n\,\}$ of its vertices, moreover all its
 edges come along with a distinguished direction from the flag below to
 the flag above. These edge directions assemble together to a direction
 for traversing each cycle in $\hat\gamma_\black$ such that the associated
 successor permutation agrees with the permutation $\sigma\,\in\,S_{2n}$
 indexing the summand graph $\hat\gamma$ in the first place. In consequence
 every colored trivalent graph $\gamma$ with $2n$ vertices occurs at least
 once up to isomorphism in the sum (\ref{pfg}) defining $\pf_n$ and all
 its occurrences in this sum share the same coefficient $\sgn\,\sigma\,=\,
 (-1)^{e(\gamma)}$ leaving no place for cancellations.

 \pfill
 In order to prove the second statement we want to count the number of times
 a given colored trivalent graph $\gamma$ is isomorphic to the summand graph
 $\hat\gamma$ indexed by a permutation $\sigma\,\in\,S_{2n}$. Fixing the
 labelling $L:\,\Vert\,\gamma\longrightarrow\{\,1,\ldots,2n\,\}$ for the
 moment as above we observe that reversing the direction of a cycle of
 $\gamma_\black$ replaces the successor permutation $\sigma$ by its inverse
 $\sigma^{-1}$ on the vertices of this cycle. A cycle is different from its
 inverse however unless the cycle is of length less than or equal to $2$,
 hence the different choices for the directions of the cycles with labelling
 $L$ fixed account for exactly $2^{g(\gamma)}$ summands in the sum (\ref{pfg}),
 where $g(\gamma)\,\in\,\N_0$ is the number of cycles of $\gamma_\black$ of
 length greater than $2$.

 Recall now that the labelling $L:\,\Vert\,\gamma\longrightarrow\{\,1,\ldots,
 2n\,\}$ was required to be such that the pairs of vertices $\{L^{-1}(1),
 L^{-1}(2)\},\,\{L^{-1}(3),L^{-1}(4)\},\,\ldots,\,\{L^{-1}(2n-1),L^{-1}(2n)\}$
 are connected by red edges. In general there are $2^n\,n!$ such labellings,
 however a different labelling $\hat L$ leads to different $2^{g(\gamma)}$
 summand graphs $\hat\gamma$ in the sum (\ref{pfg}) isomorphic to $\gamma$,
 if and only if $\hat L^{-1}\circ L\,\notin\,\overline\Aut\;\gamma$ fails to
 be a pure automorphism of $\gamma$. Hence the colored trivalent graph $\gamma$
 we consider occurs always with the same coefficient $(-1)^{e(\gamma)}$ in
 exactly
 $$
  \frac{2^n\,n!}{\#\,\overline\Aut\;\gamma}\;2^{g(\gamma)}
  \;\;=\;\;
  12^n\;n!\;
  \Big(\;\frac1{6^n}\;\frac{2^{g(\gamma)}}{\#\,\overline\Aut\;\gamma}\;\Big)
 $$
 different summands of the sum (\ref{pfg}) defining $\pf_n$; put differently
 we obtain for all $n\,\in\,\N_0$:
 $$
  \pf_n
  \;\;=\;\;
  \sum_{[\,\gamma\,]\,\in\,\Gamma^n}\frac{(-1)^{e(\gamma)}}{6^n}
  \;\frac{2^{g(\gamma)}}{\#\overline\Aut\;\gamma}\;[\,\gamma\,]\ .
 $$
 Contemplating this formula a bit the reader may easily verify that the total
 Pfaffian is the stipulated exponential in the power series completion of the
 graph algebra $\A^\bullet$.
 \qed

 \pfill
 Needless to say the formula for the total Pfaffian given in Lemma \ref{efp}
 can still be simplified using the IHX--relations. After a little bit of
 computation we obtain for the power series expansion of the total Pfaffian
 up to degree $3$ in the power series completion of $\overline\A^\bullet$:
 \begin{eqnarray*}
  \sum_{n\,\geq\,0}\pf_n
  &=&
  \exp\Big(\;
  -\,\frac14\,
  \raise-9pt\hbox{\begin{picture}(10,26)
   \put( 5, 2){\circle2}
   \put( 5,24){\circle2}
   \qbezier(4.0, 2.0)( 2.0, 2.0)( 2.0,12.0)
   \qbezier(4.0,24.0)( 2.0,24.0)( 2.0,12.0)
   \qbezier(6.1, 2.0)( 8.0, 2.0)( 8.0,12.0)
   \qbezier(6.1,24.0)( 8.0,24.0)( 8.0,12.0)
   \multiput( 5.0, 4.0)( 0.0,+2.0){10}{\circle*1}
  \end{picture}}
  \,-\,\frac18\,
  \raise-9pt\hbox{\begin{picture}(26,26)
   \put( 2, 2){\circle2}
   \put(24, 2){\circle2}
   \put( 2,24){\circle2}
   \put(24,24){\circle2}
   \put( 2.0, 3.0){\line( 0,+1){20}}
   \put(24.0, 3.0){\line( 0,+1){20}}
   \put( 3.0, 2.0){\line(+1, 0){20}}
   \put( 3.0,24.0){\line(+1, 0){20}}
   \multiput( 3.2, 3.9)( 0.0,+2.0){10}{\circle*1}
   \multiput(22.8, 3.9)( 0.0,+2.0){10}{\circle*1}
  \end{picture}}
  \,-\,\frac1{48}\,
  \raise-9pt\hbox{\begin{picture}(26,26)
   \put( 2, 2){\circle2}
   \put(24, 2){\circle2}
   \put( 2,24){\circle2}
   \put(24,24){\circle2}
   \put( 2.0, 3.0){\line( 0,+1){20}}
   \put(24.0, 3.0){\line( 0,+1){20}}
   \put( 3.0, 2.0){\line(+1, 0){20}}
   \put( 3.0,24.0){\line(+1, 0){20}}
   \multiput( 3.4, 3.4)(+1.2,+1.2){17}{\circle*1}
   \multiput(22.6, 3.4)(-1.2,+1.2){17}{\circle*1}
  \end{picture}}
  \\[-4pt]
  &&
  \quad\;-\,\frac1{24}\,
  \raise-11pt\hbox{\begin{picture}(30,30)
   \put(15, 2){\circle2}
   \put(15,28){\circle2}
   \put( 3, 8){\circle2}
   \put(27, 8){\circle2}
   \put( 3,22){\circle2}
   \put(27,22){\circle2}
   \put( 3.0, 9.0){\line( 0,+1){12}}
   \put(27.0, 9.0){\line( 0,+1){12}}
   \put(15.9, 2.5){\line(+2,+1){10}}
   \put(14.1, 2.5){\line(-2,+1){10}}
   \put(15.9,27.5){\line(+2,-1){10}}
   \put(14.1,27.5){\line(-2,-1){10}}
   \multiput( 4.9,21.9)(+1.5,+0.75){7}{\circle*1}
   \multiput( 4.9, 8.1)(+1.5,-0.75){7}{\circle*1}
   \multiput(26.1,10.1)( 0.0,+2.0){6}{\circle*1}
  \end{picture}}
  \,-\,\frac1{16}\,
  \raise-11pt\hbox{\begin{picture}(30,30)
   \put(15, 2){\circle2}
   \put(15,28){\circle2}
   \put( 3, 8){\circle2}
   \put(27, 8){\circle2}
   \put( 3,22){\circle2}
   \put(27,22){\circle2}
   \put( 3.0, 9.0){\line( 0,+1){12}}
   \put(27.0, 9.0){\line( 0,+1){12}}
   \put(15.9, 2.5){\line(+2,+1){10}}
   \put(14.1, 2.5){\line(-2,+1){10}}
   \put(15.9,27.5){\line(+2,-1){10}}
   \put(14.1,27.5){\line(-2,-1){10}}
   \multiput( 3.9,10.1)( 0.0,+2.0){6}{\circle*1}
   \multiput(15.0, 4.0)( 0.0,+2.0){12}{\circle*1}
   \multiput(26.1,10.1)( 0.0,+2.0){6}{\circle*1}
  \end{picture}}
  \,-\,\frac1{24}\,
  \raise-11pt\hbox{\begin{picture}(30,30)
   \put(15, 2){\circle2}
   \put(15,28){\circle2}
   \put( 3, 8){\circle2}
   \put(27, 8){\circle2}
   \put( 3,22){\circle2}
   \put(27,22){\circle2}
   \put( 3.0, 9.0){\line( 0,+1){12}}
   \put(27.0, 9.0){\line( 0,+1){12}}
   \put(15.9, 2.5){\line(+2,+1){10}}
   \put(14.1, 2.5){\line(-2,+1){10}}
   \put(15.9,27.5){\line(+2,-1){10}}
   \put(14.1,27.5){\line(-2,-1){10}}
   \multiput( 4.9,21.9)(+1.5,+0.75){7}{\circle*1}
   \multiput(16.1, 3.8)(+0.9,+1.5){12}{\circle*1}
   \multiput( 5.0, 8.0)(+2.0, 0.0){11}{\circle*1}
  \end{picture}}
  \,-\,\frac5{432}\,
  \raise-11pt\hbox{\begin{picture}(30,30)
   \put(15, 2){\circle2}
   \put(15,28){\circle2}
   \put( 3, 8){\circle2}
   \put(27, 8){\circle2}
   \put( 3,22){\circle2}
   \put(27,22){\circle2}
   \put( 3.0, 9.0){\line( 0,+1){12}}
   \put(27.0, 9.0){\line( 0,+1){12}}
   \put(15.9, 2.5){\line(+2,+1){10}}
   \put(14.1, 2.5){\line(-2,+1){10}}
   \put(15.9,27.5){\line(+2,-1){10}}
   \put(14.1,27.5){\line(-2,-1){10}}
   \multiput(15.0, 4.0)( 0.0,+2.0){12}{\circle*1}
   \multiput( 4.9, 8.2)(+2.0, 0.0){11}{\circle*1}
   \multiput( 4.9,21.8)(+2.0, 0.0){11}{\circle*1}
  \end{picture}}
  \,-\,\frac1{432}\,
  \raise-11pt\hbox{\begin{picture}(30,30)
   \put(15, 2){\circle2}
   \put(15,28){\circle2}
   \put( 3, 8){\circle2}
   \put(27, 8){\circle2}
   \put( 3,22){\circle2}
   \put(27,22){\circle2}
   \put( 3.0, 9.0){\line( 0,+1){12}}
   \put(27.0, 9.0){\line( 0,+1){12}}
   \put(15.9, 2.5){\line(+2,+1){10}}
   \put(14.1, 2.5){\line(-2,+1){10}}
   \put(15.9,27.5){\line(+2,-1){10}}
   \put(14.1,27.5){\line(-2,-1){10}}
   \multiput(15.0, 4.0)( 0.0,+2.0){12}{\circle*1}
   \multiput( 4.7, 9.0)(+1.9,+1.1){12}{\circle*1}
   \multiput(25.3, 9.0)(-1.9,+1.1){12}{\circle*1}
  \end{picture}}
  \,+\,\ldots\;\Big)\ .
 \end{eqnarray*}
   
 \begin{Remark}[Vanishing of Stable Curvature Invariants]
 \hfill\label{vsci}\break
  The Euler characteristic is of course multiplicative under taking
  Cartesian products of compact Riemannian manifolds, with the Theorem
  of Chern--Gauß--Bonnet in mind we anticipate that the Pfaffian of algebraic
  curvature tensors is multiplicative under taking direct sums
  $$
   \pf_{\frac{m+\hat m}2}(\;R\,\oplus\,\hat R\;)
   \;\;=\;\;
   \pf_{\frac m2}(\;R\;)\,\pf_{\frac{\hat m}2}(\;\hat R\;)
  $$
  of algebraic curvature tensors $R$ and $\hat R$ on euclidean vector spaces
  of even dimensions $m$ and $\hat m$ respectively. According to Lemma
  \ref{efp} the total Pfaffian is an exponential in the power series completion
  of $\A^\bullet$, hence the Pfaffian is in fact multiplicative as stipulated
  provided
  $$
   \pf_n(\;R\;)
   \;\;=\;\;
   0
  $$
  for every algebraic curvature tensor $R$ on a euclidean vector space $T$
  of dimension less than $2n$. Of course under this assumption already the
  $n$--th power $\frac1{n!}(-R)^n\,\in\,\LP^{2n}T^*\otimes\LP^{2n}T$ of the
  algebraic curvature tensor $R\,\in\,\Curv^-T$ considered as a bivector
  valued $2$--form vanishes.
 \end{Remark}
 
 \pfill
 Let us now come to a different set of interesting polynomials in algebraic
 curvature tensors, namely the normalized moment polynomials $\Psi^\circ_n$
 of degree $n\,\in\,\N_0$. Recall to begin with that every euclidean vector
 space $T$ with scalar product $g$ can be considered as a Riemannian
 manifold with translation invariant Riemannian metric $g$ and associated
 Laplace--Beltrami operator. For the purposes of this article we will take
 this to be the positive Laplacian
 $$
  \Delta_g
  \;\;:=\;\;
  -\;\Big(\;\frac{\partial^2}{\partial x^2_1}\;+\;\ldots\;+\;
  \frac{\partial^2}{\partial x^2_m}\;\Big)
 $$
 in a system of linear coordinates on $T$ provided by an orthonormal basis
 $x_1,\,\ldots,\,x_m$ of $T^*$:
 
 \begin{Definition}[Normalized Moment Polynomials]
 \hfill\label{mp}\break
  For every $n\,\in\,\N_0$ we define the normalized moment polynomial
  $\Psi^\circ_n$ as a homogeneous, stably invariant polynomial of degree
  $n$ in algebraic curvature tensors via the generating series:
  $$
   \sum_{n\,\geq\,0}\Psi^\circ_n(\,R\,)
   \;\;:=\;\;
   \left.\exp\Big(\,-\,\Delta_g\,\Big)\right|_0
   \left(\;\;X\;\longmapsto\;\exp\Big(\;\sum_{r\,>\,0}\frac1{2r}
   \,\tr(\;R_{\,\cdot\,,\,X}X\;)^r\;\Big)\;\;\right)\ .
  $$
 \end{Definition}

 \noindent
 Depending on the signature of the scalar product $g$ the positive Laplacian
 $\Delta_g$ may actually be the d'Alembert or wave operator $\Box$ of course.
 In order to justify calling $(\,\Psi^\circ_n\,)_{n\,\in\,\N_0}$ the normalized
 moment polynomials we consider a euclidean vector space $T$ with positive
 definite scalar product $g$. Under this assumption the sectional curvature
 function $\sec_R:\,\Gr_2T\longrightarrow\R$ associated to $R\,\in\,\Curv^-T$
 is well--defined on the Graßmannian of planes in $T$ by
 $$
  \sec_R(\;E\;)
  \;\;:=\;\;
  \frac{R(\;X,\,Y;\,Y,\,X\;)}{g_{\LP^2T}(\,X\wedge Y,\,X\wedge Y\,)}
  \;\;=\;\;
  \frac14\;\frac{\Sec(\;X,\,X;\,Y,\,Y\;)}{g_{\LP^2T}(\,X\wedge Y,\,X\wedge Y\,)}
 $$
 for linearly independent vectors $X,\,Y$ spanning $E\,=\,\span_\R\{\,X,\,
 Y\,\}\,\in\,\Gr_2T$. Endowing the Graßmannian with the Fubini--Study metric
 $g_\FS$ or measure $|\vol_\FS|$ we can think of the sectional curvature
 function as a random variable on $\Gr_2T$, its moments are then given by:
   
 \begin{Lemma}[Moments of Sectional Curvature \cite{w1}]
 \hfill\label{msc}\break
  In the positive definite case the normalized moment polynomials
  $\Psi^\circ_n$ calculate the moments $\Psi_n,\,n\,\in\,\N_0,$ of the
  sectional curvature function $\sec_R:\,\Gr_2T\longrightarrow\R$ considered
  as a random variable on the Graßmannian $\Gr_2T$ of planes in $T$ endowed
  with the Fubini--Study measure:
  $$
   \Psi_n(\,R\,)
   \;\;:=\;\;
   \frac1{\Vol(\,\Gr_2T,\,|\,\vol_\FS\,|\,)}
   \;\int_{\Gr_2T}\sec^n_R\;|\,\vol_\FS\,|
   \;\;\stackrel!=\;\;
   \frac{n!}{[\,m+2n-2\,]_{2n}}\;\Psi^\circ_n(\,R\,)\ .
  $$
  The normalization factor depends on the falling factorial
  $[\,x\,]_{2n}\,:=\,x(x-1)\ldots(x-2n+1)$.
 \end{Lemma}

 \noindent
 In passing we remark that it is possible to calculate the maximum and
 minimum of the sectional curvature function $\sec_R:\,\Gr_2T\longrightarrow
 \R$ associated to an algebraic curvature tensor $R\,\in\,\Curv^-T$ in the
 positive definite case from its moments $\Psi_n(\,R\,)$, more precisely
 \begin{eqnarray*}
  \sec_{\max}(\,R\,)
  \;\;:=\;\;
  \max_{E\,\in\,\Gr_2T}\sec_R(\,E\,)
  &=&
  -\;\LP\;+\;\lim_{N\,\to\,\infty}\left(\;
   \sum_{n\,=\,0}^N\;\;{N\choose n}\;\Psi_n(\,R\,)\;\LP^{N-n}
   \;\right)^{\frac1N}
  \\
  \sec_{\min}(\,R\,)
  \;\;:=\;\;
  \min_{E\,\in\,\Gr_2T}\sec_R(\,E\,)
  &=&
  -\;\LP\;+\;\lim_{N\,\to\,\infty}\left(\;
   \sum_{n\,=\,0}^N(-1)^n\,{N\choose n}\,\Psi_n(\,R\,)\,\LP^{N-n}
   \;\right)^{\frac1N}
 \end{eqnarray*}
 for all shift parameters $\LP\,\in\,\R$ satisfying $\LP\,\geq\,-\sec_{\min}
 (\,R\,)$ or $\LP\,\geq\,+\sec_{\max}(\,R\,)$ respectively. Needless to say
 both formulas are eventually formulas from probability theory for the
 essential supremum and infimum of an almost certainly bounded random
 variable.
 
 In difference to the moment polynomials $\Psi_n$ with their direct
 interpretation as moments the normalized moment polynomials $\Psi^\circ_n$
 turn out to be stable curvature invariants. In order to identify the
 corresponding elements of the algebra $\A^\bullet$ of colored trivalent
 graphs we use the consequence $R(U,X;X,V)\,=\,\frac14\,\Sec(X,X;U,V)\,=\,
 -\frac12\,\Sec(X,U;X,V)$ of the first Bianchi identity (\ref{1BI'}) for
 $\Sec$ in order to expand the Jacobi operator $U\longmapsto R_{U,\,X}X$
 into the sum
 $$
  R_{\;\cdot\,,\,X}X
  \;\;=\;\;
  \sum_{\mu\,=\,1}^mR(\;\cdot\;,\,X;\,X,\,dE^\up_\mu\,)\;E_\mu
  \;\;=\;\;
  -\;\frac12\;\sum_{\mu\,=\,1}^m\Sec(\;X,\;\cdot\;;\,X,\,dE^\up_\mu\;)\;E_\mu
 $$
 over a dual pair of bases $\{\,E_\mu\,\}$ and $\{\,dE_\mu\,\}$ for the
 euclidean vector space $T$ and its dual $T^*$. In turn the traces
 of the powers of the Jacobi operator become the iterated sums:
 \begin{eqnarray*}
  \lefteqn{\tr(\;R_{\;\cdot\,,\,X}X\;)^r}
  &&
  \\
  &=&
  \Big(\!-\frac12\Big)^r\!\!
  \sum_{\mu_1,\,\ldots,\,\mu_r\,=\,1}^m\!\!
  \Sec(\,X,E_{\mu_1};X,dE^\up_{\mu_2}\,)\,
  \Sec(\,X,E_{\mu_2};X,dE^\up_{\mu_3}\,)\,\ldots\,
  \Sec(\,X,E_{\mu_r};X,dE^\up_{\mu_1}\,)\ .
 \end{eqnarray*}
 Although we have used the exponential $\left.\exp(\,-\,\Delta_g\,)\right|_0$
 of the positive Laplacian on $T$ in Definition \ref{msc}, the rescaled
 exponential $\left.\exp(\,-\,\frac12\Delta_g\,)\right|_0$ has a rather
 compelling combinatorial description as a closure operation: Applying it
 to the diagonal polynomial $X\longmapsto a(X,\ldots,X)$ arising from a not
 necessarily symmetric form $a\,\in\,\bigotimes^{2n}T^*$ we obtain an
 iterated sum
 \begin{eqnarray*}
  \cls[\;a\;]
  &:=& 
  \left.\exp(\;-\,{\textstyle\frac12}\,\Delta_g\;)\right|_0
  \Big(\;X\,\longmapsto\,a(\,\underbrace{X,\,\ldots,\,X}_{2n}\,)\;\Big)
  \\
  &=&
  \sum_{{\scriptstyle\theta\,\in\,S_{2n}\mathrm{\;fix\;point\;free}}
  \atop{\scriptstyle\mathrm{involution\;of}\,\{\,1,\,\ldots,\,2n\,\}}}
  \sum_{{\scriptstyle\mu:\,\{\,1,\,\ldots,\,2n\,\}\longrightarrow\{\,1,
   \,\ldots,\,m\,\}}\atop{\scriptstyle\mu\,=\,\mu\,\circ\,\theta
   \mathrm{\;invariant\;under\;}\theta}}
  a(\;E_{\mu(1)},\,E_{\mu(2)},\,\ldots,\,E_{\mu(2n)}\;)
 \end{eqnarray*}
 over an orthonormal basis $E_1,\,\ldots,\,E_m$ of the euclidean vector space
 $T$, where the two sums extend over all fix point free involutions $\theta$
 of the index set $\{\,1,\,\ldots,\,2n\,\}$ and all $\theta$--invariant maps
 $\mu:\,\{\,1,\,\ldots,\,2n\,\}\longrightarrow\{\,1,\,\ldots,\,m\,\}$ in the
 sense $\mu\,=\,\mu\,\circ\,\theta$. For simplicity of the exposition we have
 pretended again that the scalar product $g$ is positive definite to omit the
 sign factors necessary otherwise. In terms of colored trivalent graphs we
 can hence write the Definition \ref{mp} of the generating series of the
 normalized moment polynomials in the form
 \begin{equation}\label{psi}
  \sum_{n\,\geq\,0}\Psi^\circ_n
  \;\;:=\;\;
  \cls\;\Bigg[\;\exp\Bigg(\;\sum_{r\,>\,0}
  \frac{(-1)^r}{2r}
  \underbrace{\raise-10pt\hbox{\begin{picture}(180,24)
   \put(  4,14){\circle2}
   \put( 24,14){\circle2}
   \put( 44,14){\circle2}
   \put( 64,14){\circle2}
   \put(116,14){\circle2}
   \put(136,14){\circle2}
   \put(156,14){\circle2}
   \put(176,14){\circle2}
   \put(  4,15){\line( 0,+1){8}}
   \put(  4,13){\line( 0,-1){5}}
   \put( 24,15){\line( 0,+1){8}}
   \put( 44,15){\line( 0,+1){8}}
   \put( 64,15){\line( 0,+1){8}}
   \put(116,15){\line( 0,+1){8}}
   \put(136,15){\line( 0,+1){8}}
   \put(156,15){\line( 0,+1){8}}
   \put(176,15){\line( 0,+1){8}}
   \put(176,13){\line( 0,-1){5}}
   \put( 26,11){\line(+1, 0){16}}
   \put( 66,11){\line(+1, 0){10}}
   \put(114,11){\line(-1, 0){10}}
   \put(154,11){\line(-1, 0){16}}
   \put( 10, 2){\line(+1, 0){160}}
   \put( 83,10.5){$\ldots$}
   \qbezier( 24.0,13.0)( 24.0,11.0)( 26.0,11.0)
   \qbezier( 44.0,13.0)( 44.0,11.0)( 42.0,11.0)
   \qbezier( 64.0,13.0)( 64.0,11.0)( 66.0,11.0)
   \qbezier(116.0,13.0)(116.0,11.0)(114.0,11.0)
   \qbezier(136.0,13.0)(136.0,11.0)(138.0,11.0)
   \qbezier(156.0,13.0)(156.0,11.0)(154.0,11.0)
   \qbezier(  4.0, 8.0)(  4.0, 5.4)(  5.7, 3.7)
   \qbezier(  5.7, 3.7)(  7.4, 2.0)( 10.0, 2.0)
   \qbezier(176.0, 8.0)(176.0, 5.4)(174.3, 3.7)
   \qbezier(174.3, 3.7)(172.6, 2.0)(170.0, 2.0)
   \multiput(  6,14)(2,0){9}{\circle*1}
   \multiput( 46,14)(2,0){9}{\circle*1}
   \multiput(118,14)(2,0){9}{\circle*1}
   \multiput(158,14)(2,0){9}{\circle*1}
  \end{picture}}}_{r\mathrm{\;pairs\;of\;vertices}}
  \;\Bigg)\;\Bigg]\ ,
 \end{equation}
 where the closure operation now sums over all ways to join up the open flags
 of the argument in pairs to form a colored trivalent graph. In this formula
 the difference between the rescaled exponential $\left.\exp(\,-\,\frac12\,
 \Delta_g\,)\right|_0$ and the actual exponential $\left.\exp(\,-\,\Delta_g\,)
 \right|_0$ of the Laplacian used in Definition \ref{msc} has been compensated
 by changing the original factor $(-\frac12)^r$ to $(-1)^r$.

 \begin{Lemma}[Generating Series of Normalized Moment Polynomials]
 \hfill\label{efpsi}\break
  In the power series completion of the algebra $\A^\bullet$ of colored
  trivalent graphs the generating series of the normalized moment polynomials
  $(\,\Psi^\circ_n\,)_{n\,\in\,\N_0}$ can be written as the exponential
  of a sum over all isomorphism classes of connected colored trivalent graphs
  $\gamma$ with the additional property that all cycles of the associated
  black subgraph $\gamma_\black$ are cycles of even length:
  $$
   \sum_{n\,\geq\,0}\Psi^\circ_n
   \;\;=\;\;
   \exp\left(\;
    \sum_{{\scriptstyle[\,\gamma\,]\,\in\,\Gamma^\bullet_\conn}\atop
     {\scriptstyle\gamma_\black\mathrm{\;even\;cycles}}}
    (-1)^{n(\gamma)}\frac{2^{g(\gamma)}}{\#\overline\Aut\;\gamma}\;
    [\,\gamma\,]\;\right)\ .
  $$
  As before $n(\gamma)\,:=\,\frac12\,\#\Vert\,\gamma$ and $g(\gamma)$ counts
  the cycles of $\gamma_\black$ of length greater than $2$.
 \end{Lemma}

 \proof
 The most important observation by far to understand the proof is that the
 black edges in the expansion (\ref{psi}) of the generating series of the
 normalized moment polynomials $(\,\Psi^\circ_n\,)_{n\,\in\,\N_0}$ either
 are ``old'' edges present already in the argument or are ``new'' edges joined
 up during the closure operation. In order to formalize this idea we define
 a tricoloring for a given colored trivalent graph $\gamma$ as an extension
 $\col_\ext:\,\Edge\,\gamma\longrightarrow\{\,\red,\,\black_\old,\,
 \black_\new\,\}$ of its coloring $\col:\,\Edge\,\gamma\longrightarrow
 \{\,\red,\,\black\,\}$ in the sense $\col\,=\,\pr\circ\col_\ext$ for the
 projection $\pr$ implied by notation such that the edges adjacent to every
 vertex are all colored differently. By construction every summand in the
 expansion (\ref{psi}) of the generating series $\sum\Psi^\circ_n$ is
 naturally a colored trivalent graph with a distinguished tricoloring
 $\col_\ext$.

 A necessary condition for the existence of a tricoloring for a colored
 trivalent graph $\gamma$ is that all cycles of the black subgraph
 $\gamma_\black$ associated to $\gamma$ are cycles of even length, because
 every tricoloring $\col_\ext$ necessarily colors the edges of $\gamma_\black$
 alternately in $\black_\old$ and $\black_\new$ along each cycle. Provided
 this necessary condition is met we can on the other hand color the edges
 alternately along each of the $e(\gamma)\,\in\,\N$ cycles of $\gamma_\black$
 starting with either $\black_\old$ or $\black_\new$ to obtain exactly
 $2^{e(\gamma)}$ different tricolorings for $\gamma$. In difference to
 the black subgraph $\gamma_\black$ associated to $\gamma$ the colored
 bivalent subgraph $\gamma_\old$ obtained by removing all edges colored
 in $\black_\new$ certainly depends on the chosen tricoloring $\col_\ext$.

 The set of all tricolorings of a colored trivalent graph $\gamma$ comes
 along with a natural action of the automorphism group $\Aut\,\gamma$ of
 $\gamma$ in such a way that the stabilizer subgroup of a tricoloring
 $\col_\ext$ agrees with the automorphism group $\Aut\,(\,\gamma,\col_\ext\,)$
 of $\gamma$ considered as a tricolored trivalent graph. Both the tricolored
 graph $(\,\gamma,\col_\ext\,)$ and the colored old subgraph $\gamma_\old$
 associated to $\gamma$ and a given tricoloring $\col_\ext$ have only pure
 automorphisms in the sense that all their automorphisms are completely
 determined by a permutation of their common set $\Vert\,\gamma$ of vertices,
 because the flags adjacent to each vertex are all colored differently:
 $$
  \Aut\,(\,\gamma,\,\col_\ext\,)
  \;\;\cong\;\;
  \overline\Aut\,(\,\gamma,\,\col_\ext\,)
  \qquad\qquad
  \Aut\,\gamma_\old
  \;\;\cong\;\;
  \overline\Aut\,\gamma_\old\ .
 $$
 Interestingly the same conclusion does not hold for the original colored
 trivalent graph $\gamma$, it may have trivial automorphisms fixing all
 vertices, but not all flags. More precisely an automorphism of a colored
 trivalent graph $\gamma$ fixing all its vertices can only swap the two
 flags or the two edges respectively in a short cycle of the black subgraph
 $\gamma_\black$ of length $1$ or $2$. Under the additional assumption that
 all cycles of $\gamma_\black$ have even length we thus obtain
 \begin{equation}\label{ta}
  \#\,\Aut_\circ\gamma
  \;\;=\;\;
  2^{e(\gamma)\,-\,g(\gamma)}\ ,
 \end{equation}
 where $g(\gamma)$ equals the number of cycles of the black subgraph
 $\gamma_\black$ of length greater than $2$.

 \pfill
 Coming back to the proof we consider in a first step the class of colored
 bivalent graphs with edges colored alternately with colors $\red$ and
 $\black_\old$. Due to the additional coloring the automorphism group of
 a connected colored bivalent graph acts simply transitively on its set
 of vertices, in turn the automorphism group of a colored bivalent graph
 with cycles up to length $2r$ and $d_1$ cycles of length $2$, $d_2$ cycles
 of length $4$ etc.~is a finite group of order
 $$
  \#\Aut\Big(\;\;\raise-5pt\hbox{\begin{picture}(28,16)
   \put( 4,14){\circle2}
   \put(24,14){\circle2}
   \put( 4,13){\line( 0,-1){5}}
   \put(24,13){\line( 0,-1){5}}
   \put(10, 2){\line(+1, 0){8}}
   \multiput( 6,14)(2,0){9}{\circle*1}
   \qbezier( 4.0, 8.0)( 4.0, 5.4)( 5.7, 3.7)
   \qbezier( 5.7, 3.7)( 7.4, 2.0)(10.0, 2.0)
   \qbezier(24.0, 8.0)(24.0, 5.4)(22.3, 3.7)
   \qbezier(22.3, 3.7)(20.6, 2.0)(18.0, 2.0)
  \end{picture}}^{d_1}\;\raise-5pt\hbox{\begin{picture}(68,16)
   \put( 4,14){\circle2}
   \put(24,14){\circle2}
   \put(44,14){\circle2}
   \put(64,14){\circle2}
   \put( 4,13){\line( 0,-1){5}}
   \put(64,13){\line( 0,-1){5}}
   \put(26,11){\line(+1, 0){16}}
   \put(10, 2){\line(+1, 0){48}}
   \multiput( 6,14)(2,0){9}{\circle*1}
   \multiput(46,14)(2,0){9}{\circle*1}
   \qbezier(24.0,13.0)(24.0,11.0)(26.0,11.0)
   \qbezier(44.0,13.0)(44.0,11.0)(42.0,11.0)
   \qbezier( 4.0, 8.0)( 4.0, 5.4)( 5.7, 3.7)
   \qbezier( 5.7, 3.7)( 7.4, 2.0)(10.0, 2.0)
   \qbezier(64.0, 8.0)(64.0, 5.4)(62.3, 3.7)
   \qbezier(62.3, 3.7)(60.6, 2.0)(58.0, 2.0)
  \end{picture}}^{d_2}\;\cdots\;
  \underbrace{\raise-5pt\hbox{\begin{picture}(100,16)
   \put( 4,14){\circle2}
   \put(24,14){\circle2}
   \put(76,14){\circle2}
   \put(96,14){\circle2}
   \put( 4,13){\line( 0,-1){5}}
   \put(96,13){\line( 0,-1){5}}
   \put(26,11){\line(+1, 0){10}}
   \put(74,11){\line(-1, 0){10}}
   \put(10, 2){\line(+1, 0){80}}
   \put(43,10.5){$\ldots$}
   \qbezier(24.0,13.0)(24.0,11.0)(26.0,11.0)
   \qbezier(76.0,13.0)(76.0,11.0)(74.0,11.0)
   \qbezier( 4.0, 8.0)( 4.0, 5.4)( 5.7, 3.7)
   \qbezier( 5.7, 3.7)( 7.4, 2.0)(10.0, 2.0)
   \qbezier(96.0, 8.0)(96.0, 5.4)(94.3, 3.7)
   \qbezier(94.3, 3.7)(92.6, 2.0)(90.0, 2.0)
   \multiput( 6,14)(+2, 0){9}{\circle*1}
   \multiput(78,14)(+2, 0){9}{\circle*1}
  \end{picture}}^{\hbox to0pt{$\scriptstyle d_r$\hss}}}_{r
  \mathrm{\;pairs\;of\;vertices}}\;\;\Big)
  \;\;=\;\;
  d_1!\,2^{d_1}\,d_2!\,4^{d_2}\,\ldots\,d_r!\,(2r)^{d_r}\ ,
 $$
 because we can permute all the cycles of length $2$, all the cycles of length
 $4$ etc.~In consequence the generating power series of the set of isomorphism
 classes of colored bivalent graphs reads: 
 $$
  \exp\Bigg(\;\sum_{r\,>\,0}
  \frac{(-1)^r}{2r}
  \underbrace{\raise-5pt\hbox{\begin{picture}(180,16)
   \put(  4,14){\circle2}
   \put( 24,14){\circle2}
   \put( 44,14){\circle2}
   \put( 64,14){\circle2}
   \put(116,14){\circle2}
   \put(136,14){\circle2}
   \put(156,14){\circle2}
   \put(176,14){\circle2}
   \put(  4,13){\line( 0,-1){5}}
   \put(176,13){\line( 0,-1){5}}
   \put( 26,11){\line(+1, 0){16}}
   \put( 66,11){\line(+1, 0){10}}
   \put(114,11){\line(-1, 0){10}}
   \put(154,11){\line(-1, 0){16}}
   \put( 10, 2){\line(+1, 0){160}}
   \put( 83,10.5){$\ldots$}
   \qbezier( 24.0,13.0)( 24.0,11.0)( 26.0,11.0)
   \qbezier( 44.0,13.0)( 44.0,11.0)( 42.0,11.0)
   \qbezier( 64.0,13.0)( 64.0,11.0)( 66.0,11.0)
   \qbezier(116.0,13.0)(116.0,11.0)(114.0,11.0)
   \qbezier(136.0,13.0)(136.0,11.0)(138.0,11.0)
   \qbezier(156.0,13.0)(156.0,11.0)(154.0,11.0)
   \qbezier(  4.0, 8.0)(  4.0, 5.4)(  5.7, 3.7)
   \qbezier(  5.7, 3.7)(  7.4, 2.0)( 10.0, 2.0)
   \qbezier(176.0, 8.0)(176.0, 5.4)(174.3, 3.7)
   \qbezier(174.3, 3.7)(172.6, 2.0)(170.0, 2.0)
   \multiput(  6,14)(2,0){9}{\circle*1}
   \multiput( 46,14)(2,0){9}{\circle*1}
   \multiput(118,14)(2,0){9}{\circle*1}
   \multiput(158,14)(2,0){9}{\circle*1}
  \end{picture}}}_{r\mathrm{\;pairs\;of\;vertices}}
  \;\Bigg)
  \;\;=\;\;
  \sum_{[\,\gamma_\old\,]}
  \frac{(-1)^{n(\gamma_\old)}}{\#\Aut\;\gamma_\old}[\,\gamma_\old\,]\ .
 $$
 By definition the closure operation in the expansion (\ref{psi}) of the
 generating power series of the normalized moment polynomials $\Psi^\circ_n$
 is a sum over all fixed point free involutions $\sigma$ of the vertex set
 $\Vert\,\gamma_\old$ of the colored bivalent graph argument $\gamma_\old$
 connecting all orbit pairs $\{\,v,\,\sigma(v)\,\}$ by an edge colored
 $\black_\new$ to obtain a tricolored trivalent graph $(\,\gamma,\,
 \col_\ext\,)$.

 The automorphism group $\Aut\,\gamma_\old\,\cong\,\overline\Aut\,\gamma_\old$
 of the colored bivalent graph $\gamma$ on the other hand acts on $\Vert\,
 \gamma_\old$ and in turn on the set of fixed point free involutions $\sigma$
 of $\Vert\,\gamma_\old$ by conjugation. The orbits of this action correspond
 bijectively to the isomorphism classes $[\,\gamma,\,\col_\ext\,]$ of
 tricolored trivalent graphs with old subgraph isomorphic to $\gamma_\old$,
 while the stabilizer of a given fixed point free involution $\sigma$ is
 essentially the automorphism group of the corresponding tricolored trivalent
 graph $(\,\gamma,\,\col_\ext\,)$. Counting the number of all fixed point free
 involutions $\tilde\sigma$ of $\Vert\,\gamma_\old$ resulting in a tricolored
 trivalent graph in the class $[\,\gamma,\,\col_\ext\,]$ thus amounts to
 calculating the length of the orbit of the fixed point free involution
 $\sigma$, in this way we obtain for the generating power series (\ref{psi})
 of the normalized moment polynomials:
 $$
  \sum_{n\,\geq\,0}\Psi^\circ_n
  \;\;=\;\;
  \cls\left[\;\sum_{[\,\gamma_\old\,]}
  \frac{(-1)^{n(\gamma_\old)}}{\#\Aut\;\gamma_\old}[\,\gamma_\old\,]\;\right]
  \;\;=\;\;
  \sum_{[\,\gamma,\,\col_\ext\,]}\frac{(-1)^{n(\gamma)}}{\#\Aut\;\gamma_\old}\;
  \frac{\#\Aut\;\gamma_\old}{\#\Aut\,(\,\gamma,\col_\ext\,)}\;[\,\gamma\,]\ .
 $$
 The problem with this expansion of the generating power series of the
 polynomials $\Psi^\circ_n$ is that the resulting sum is not effective,
 it sums multiples of the isomorphism class $[\,\gamma\,]$ of the colored
 trivalent graph underlying the tricolored isomorphism class $[\,\gamma,\,
 \col_\ext\,]$. In order to obtain an effective formula for the generating
 power series (\ref{psi}) we thus need to sum for a given isomorphism class
 $[\,\gamma\,]$ of a colored trivalent graph $\gamma$ over all possible
 isomorphism classes $[\,\gamma,\,\col_\ext\,]$ of tricolored trivalent
 graphs extending $\gamma$.

 A necessary and sufficient condition for the existence of some isomorphism
 class extension $[\,\gamma,\,\col_\ext\,]$ of a given colored trivalent graph
 $\gamma$ is that all cycles of the black subgraph $\gamma_\black$ are cycles
 of even length; in this case the different isomorphism classes correspond to
 the orbits of the automorphism group $\Aut\,\gamma$ of $\gamma$ in the set of
 tricolorings $\col_\ext$. Instead of summing over orbits of tricolorings it is
 more convenient to sum over tricolorings for $\gamma$ weighted by the inverse
 length of their respective orbits, in this way we arrive eventually at the
 formula:
 \begin{eqnarray*}
  \sum_{n\,\geq\,0}\Psi^\circ_n
  &=&
  \;\;\;\;\sum_{[\,\gamma,\,\col_\ext\,]}
  \frac{(-1)^{n(\gamma)}}{\#\Aut\,(\,\gamma,\col_\ext\,)}\;[\,\gamma\,]
  \\
  &=&
  \!\!\!\sum_{{\scriptstyle[\,\gamma\,]\,\in\,\Gamma^\bullet}
   \atop{\scriptstyle\gamma_\black\mathrm{\;even\;cycles}}}
  \sum_{\col_\ext\mathrm{\;tricoloring}}
  \left(\,\frac{\#\Aut\,\gamma}{\#\Aut(\,\gamma,\,\col_\ext\,)}\,\right)^{-1}
  \frac{(-1)^{n(\gamma)}}{\#\Aut\,(\,\gamma,\,\col_\ext\,)}\;[\,\gamma\,]
  \\
  &=&
  \!\!\!\sum_{{\scriptstyle[\,\gamma\,]\,\in\,\Gamma^\bullet}
   \atop{\scriptstyle\gamma_\black\mathrm{\;even\;cycles}}}
  (-1)^{n(\gamma)}\;\frac{2^{e(\gamma)}}{\#\Aut\;\gamma}\;[\,\gamma\,]
  \;\;=
  \sum_{{\scriptstyle[\,\gamma\,]\,\in\,\Gamma^\bullet}
   \atop{\scriptstyle\gamma_\black\mathrm{\;even\;cycles}}}
  (-1)^{n(\gamma)}
  \;\frac{2^{g(\gamma)}}{\#\overline\Aut\;\gamma}\;[\,\gamma\,]\ .
 \end{eqnarray*}
 For the last equality we have used equation (\ref{ta}) in the form
 $\#\Aut\,\gamma\,=\,2^{e(\gamma)-g(\gamma)}\,\#\overline\Aut\,\gamma$.
 Similarly to the Pfaffian polynomials discussed above the mere form of
 the result tells us that the generating power series of the normalized
 moment polynomials $\Psi^\circ_n$ is the exponential of the sum on the
 right hand side taken over connected colored trivalent graphs only.
 \qed

 \pfill
 Of course Lemma \ref{efpsi} describes the expansion of the generating
 power series of the normalized moment polynomials $\Psi^\circ_n$ in the
 power series completion of the algebra $\A^\bullet$ of colored trivalent
 graphs. Simplifying the result using IHX--relations we find after some
 not too messy calculations in the power series completion of $\overline
 \A^\bullet$ the following expansion up to degree $3$:
 \begin{eqnarray*}
  \sum_{n\,\geq\,0}\Psi^\circ_n
  &=&
  \exp\Big(\;
  -\,\frac12\,
  \raise-9pt\hbox{\begin{picture}(10,26)
   \put( 5, 2){\circle2}
   \put( 5,24){\circle2}
   \qbezier(4.0, 2.0)( 2.0, 2.0)( 2.0,12.0)
   \qbezier(4.0,24.0)( 2.0,24.0)( 2.0,12.0)
   \qbezier(6.1, 2.0)( 8.0, 2.0)( 8.0,12.0)
   \qbezier(6.1,24.0)( 8.0,24.0)( 8.0,12.0)
   \multiput( 5.0, 4.0)( 0.0,+2.0){10}{\circle*1}
  \end{picture}}
  \,+\,\frac12\,
  \raise-9pt\hbox{\begin{picture}(26,26)
   \put( 2, 2){\circle2}
   \put(24, 2){\circle2}
   \put( 2,24){\circle2}
   \put(24,24){\circle2}
   \put( 2.0, 3.0){\line( 0,+1){20}}
   \put(24.0, 3.0){\line( 0,+1){20}}
   \put( 3.0, 2.0){\line(+1, 0){20}}
   \put( 3.0,24.0){\line(+1, 0){20}}
   \multiput( 3.2, 3.9)( 0.0,+2.0){10}{\circle*1}
   \multiput(22.8, 3.9)( 0.0,+2.0){10}{\circle*1}
  \end{picture}}
  \,-\,\frac14\,
  \raise-9pt\hbox{\begin{picture}(26,26)
   \put( 2, 2){\circle2}
   \put(24, 2){\circle2}
   \put( 2,24){\circle2}
   \put(24,24){\circle2}
   \put( 2.0, 3.0){\line( 0,+1){20}}
   \put(24.0, 3.0){\line( 0,+1){20}}
   \put( 3.0, 2.0){\line(+1, 0){20}}
   \put( 3.0,24.0){\line(+1, 0){20}}
   \multiput( 3.4, 3.4)(+1.2,+1.2){17}{\circle*1}
   \multiput(22.6, 3.4)(-1.2,+1.2){17}{\circle*1}
  \end{picture}}
  \\[-4pt]
  &&
  \qquad\;\;-\,\frac13\,
  \raise-11pt\hbox{\begin{picture}(30,30)
   \put(15, 2){\circle2}
   \put(15,28){\circle2}
   \put( 3, 8){\circle2}
   \put(27, 8){\circle2}
   \put( 3,22){\circle2}
   \put(27,22){\circle2}
   \put( 3.0, 9.0){\line( 0,+1){12}}
   \put(27.0, 9.0){\line( 0,+1){12}}
   \put(15.9, 2.5){\line(+2,+1){10}}
   \put(14.1, 2.5){\line(-2,+1){10}}
   \put(15.9,27.5){\line(+2,-1){10}}
   \put(14.1,27.5){\line(-2,-1){10}}
   \multiput( 4.9,21.9)(+1.5,+0.75){7}{\circle*1}
   \multiput( 4.9, 8.1)(+1.5,-0.75){7}{\circle*1}
   \multiput(26.1,10.1)( 0.0,+2.0){6}{\circle*1}
  \end{picture}}
  \,-\,\frac12\,
  \raise-11pt\hbox{\begin{picture}(30,30)
   \put(15, 2){\circle2}
   \put(15,28){\circle2}
   \put( 3, 8){\circle2}
   \put(27, 8){\circle2}
   \put( 3,22){\circle2}
   \put(27,22){\circle2}
   \put( 3.0, 9.0){\line( 0,+1){12}}
   \put(27.0, 9.0){\line( 0,+1){12}}
   \put(15.9, 2.5){\line(+2,+1){10}}
   \put(14.1, 2.5){\line(-2,+1){10}}
   \put(15.9,27.5){\line(+2,-1){10}}
   \put(14.1,27.5){\line(-2,-1){10}}
   \multiput( 3.9,10.1)( 0.0,+2.0){6}{\circle*1}
   \multiput(15.0, 4.0)( 0.0,+2.0){12}{\circle*1}
   \multiput(26.1,10.1)( 0.0,+2.0){6}{\circle*1}
  \end{picture}}
  \,+\,
  \raise-11pt\hbox{\begin{picture}(30,30)
   \put(15, 2){\circle2}
   \put(15,28){\circle2}
   \put( 3, 8){\circle2}
   \put(27, 8){\circle2}
   \put( 3,22){\circle2}
   \put(27,22){\circle2}
   \put( 3.0, 9.0){\line( 0,+1){12}}
   \put(27.0, 9.0){\line( 0,+1){12}}
   \put(15.9, 2.5){\line(+2,+1){10}}
   \put(14.1, 2.5){\line(-2,+1){10}}
   \put(15.9,27.5){\line(+2,-1){10}}
   \put(14.1,27.5){\line(-2,-1){10}}
   \multiput( 4.9,21.9)(+1.5,+0.75){7}{\circle*1}
   \multiput(16.1, 3.8)(+0.9,+1.5){12}{\circle*1}
   \multiput( 5.0, 8.0)(+2.0, 0.0){11}{\circle*1}
  \end{picture}}
  \,-\,\frac16\,
  \raise-11pt\hbox{\begin{picture}(30,30)
   \put(15, 2){\circle2}
   \put(15,28){\circle2}
   \put( 3, 8){\circle2}
   \put(27, 8){\circle2}
   \put( 3,22){\circle2}
   \put(27,22){\circle2}
   \put( 3.0, 9.0){\line( 0,+1){12}}
   \put(27.0, 9.0){\line( 0,+1){12}}
   \put(15.9, 2.5){\line(+2,+1){10}}
   \put(14.1, 2.5){\line(-2,+1){10}}
   \put(15.9,27.5){\line(+2,-1){10}}
   \put(14.1,27.5){\line(-2,-1){10}}
   \multiput(15.0, 4.0)( 0.0,+2.0){12}{\circle*1}
   \multiput( 4.9, 8.2)(+2.0, 0.0){11}{\circle*1}
   \multiput( 4.9,21.8)(+2.0, 0.0){11}{\circle*1}
  \end{picture}}
  \,-\,\frac16\,
  \raise-11pt\hbox{\begin{picture}(30,30)
   \put(15, 2){\circle2}
   \put(15,28){\circle2}
   \put( 3, 8){\circle2}
   \put(27, 8){\circle2}
   \put( 3,22){\circle2}
   \put(27,22){\circle2}
   \put( 3.0, 9.0){\line( 0,+1){12}}
   \put(27.0, 9.0){\line( 0,+1){12}}
   \put(15.9, 2.5){\line(+2,+1){10}}
   \put(14.1, 2.5){\line(-2,+1){10}}
   \put(15.9,27.5){\line(+2,-1){10}}
   \put(14.1,27.5){\line(-2,-1){10}}
   \multiput(15.0, 4.0)( 0.0,+2.0){12}{\circle*1}
   \multiput( 4.7, 9.0)(+1.9,+1.1){12}{\circle*1}
   \multiput(25.3, 9.0)(-1.9,+1.1){12}{\circle*1}
  \end{picture}}
  \,+\,\ldots\;\Big)\ .
 \end{eqnarray*}
\section{Curvature Identities for Einstein Manifolds}
\label{standard}
 The graphical calculus for stable curvature invariants is well suited to
 derive curvature identities generalizing the Hitchin--Thorpe inequality
 in dimension $m\,=\,4$. In a sense the Hitchin--Thorpe inequality deals
 with the expected value and the variance of sectional curvature, its third
 moment is related to the cubic polynomial $\Theta_3(\,R\,)\,:=\,g_{\LP^2T^*
 \otimes\LP^2T^*}(\,q(R)\star R,\,R\,)$ of central importance to this section.
 Using the expansions of the Pfaffian and normalized moment polynomials we
 find an essentially unique linear relation between $\pf_3$, $\Psi^\circ_3$
 and $\Theta_3$, whose integral over compact Einstein manifolds of dimension
 $m\,\geq\,3$ results in Theorem \ref{ineq}. In contrast to the preceding,
 algebraic sections Riemannian manifolds are tacitly assumed to be endowed
 with positive definite metrics $g\,>\,0$ throughout this section.
 
 \pfill
 It is relatively easy to identify the generators (\ref{gen6}) of the
 reduced graph algebra $\overline\A^\bullet$ up to degree $2$ using
 only the Ricci identity (\ref{ihex}) and the argument leading to equation
 (\ref{eqnorm}) together with the IHX--congruence (\ref{ihq}). In this way
 we find
 $
  \left[\raise-4pt\hbox{\begin{picture}(10,14)
   \put( 5, 0){\circle2}
   \put( 5,14){\circle2}
   \qbezier( 4.0, 0.0)( 1.0, 0.0)( 1.0, 7.0)
   \qbezier( 4.0,14.0)( 1.0,14.0)( 1.0, 7.0)
   \qbezier( 6.0, 0.0)( 9.0, 0.0)( 9.0, 7.0)
   \qbezier( 6.0,14.0)( 9.0,14.0)( 9.0, 7.0)
   \multiput( 5.0, 2.2)( 0.0,+1.6){7}{\circle*1}
  \end{picture}}\right](\,R\,)\,=\,-2\kappa
 $
 as well as:
 \begin{equation}\label{deg2}
  \left[\;\raise-4pt\hbox{\begin{picture}(16,16)
   \put( 0, 0){\circle2}
   \put(16, 0){\circle2}
   \put( 0,16){\circle2}
   \put(16,16){\circle2}
   \put( 0, 1){\line( 0,+1){14}}
   \put(16, 1){\line( 0,+1){14}}
   \put( 1, 0){\line(+1, 0){14}}
   \put( 1,16){\line(+1, 0){14}}
   \multiput( 1.1, 1.7)( 0.0,+1.8){8}{\circle*1}
   \multiput(14.9, 1.7)( 0.0,+1.8){8}{\circle*1}
  \end{picture}}\;\right](\,R\,)
  \;\;=\;\;
  8\,g_{\S^2T^*}(\,\Ric,\,\Ric\,)
  \qquad
  \left[\;\raise-4pt\hbox{\begin{picture}(16,16)
   \put( 0, 0){\circle2}
   \put(16, 0){\circle2}
   \put( 0,16){\circle2}
   \put(16,16){\circle2}
   \put( 0, 1){\line( 0,+1){14}}
   \put(16, 1){\line( 0,+1){14}}
   \put( 1, 0){\line(+1, 0){14}}
   \put( 1,16){\line(+1, 0){14}}
   \multiput( 1.4, 1.4)(+1.1,+1.1){13}{\circle*1}
   \multiput( 1.4,14.6)(+1.1,-1.1){13}{\circle*1}
  \end{picture}}\;\right](\,R\,)
  \;\;=\;\;
  -\,24\,g_{\LP^2T^*\otimes\LP^2T^*}(\,R,\,R\,)\ .
 \end{equation}
 Moreover the expansions of the total Pfaffian and the normalized moment
 polynomials given in Lemma \ref{efp} and Lemma \ref{efpsi} respectively
 duely simplified up to degree $3$ in Section \ref{special} imply:
 $$
  \pf_2
  \;\;=\;\;
  \frac1{32}\;\raise-4pt\hbox{\begin{picture}(10,14)
   \put( 5, 0){\circle2}
   \put( 5,14){\circle2}
   \qbezier( 4.0, 0.0)( 1.0, 0.0)( 1.0, 7.0)
   \qbezier( 4.0,14.0)( 1.0,14.0)( 1.0, 7.0)
   \qbezier( 6.0, 0.0)( 9.0, 0.0)( 9.0, 7.0)
   \qbezier( 6.0,14.0)( 9.0,14.0)( 9.0, 7.0)
   \multiput( 5.0, 2.2)( 0.0,+1.6){7}{\circle*1}
  \end{picture}}^2
  \,-\,
  \frac18\;\raise-5pt\hbox{\begin{picture}(16,16)
   \put( 0, 0){\circle2}
   \put(16, 0){\circle2}
   \put( 0,16){\circle2}
   \put(16,16){\circle2}
   \put( 0, 1){\line( 0,+1){14}}
   \put(16, 1){\line( 0,+1){14}}
   \put( 1, 0){\line(+1, 0){14}}
   \put( 1,16){\line(+1, 0){14}}
   \multiput( 1.1, 1.7)( 0.0,+1.8){8}{\circle*1}
   \multiput(14.9, 1.7)( 0.0,+1.8){8}{\circle*1}
  \end{picture}}
  \,-\,
  \frac1{48}\;\raise-5pt\hbox{\begin{picture}(16,16)
   \put( 0, 0){\circle2}
   \put(16, 0){\circle2}
   \put( 0,16){\circle2}
   \put(16,16){\circle2}
   \put( 0, 1){\line( 0,+1){14}}
   \put(16, 1){\line( 0,+1){14}}
   \put( 1, 0){\line(+1, 0){14}}
   \put( 1,16){\line(+1, 0){14}}
   \multiput( 1.4, 1.4)(+1.1,+1.1){13}{\circle*1}
   \multiput( 1.4,14.6)(+1.1,-1.1){13}{\circle*1}
  \end{picture}}
  \qquad\qquad
  \Psi^\circ_2
  \;\;=\;\;
  \frac18\;\raise-4pt\hbox{\begin{picture}(10,14)
   \put( 5, 0){\circle2}
   \put( 5,14){\circle2}
   \qbezier( 4.0, 0.0)( 1.0, 0.0)( 1.0, 7.0)
   \qbezier( 4.0,14.0)( 1.0,14.0)( 1.0, 7.0)
   \qbezier( 6.0, 0.0)( 9.0, 0.0)( 9.0, 7.0)
   \qbezier( 6.0,14.0)( 9.0,14.0)( 9.0, 7.0)
   \multiput( 5.0, 2.2)( 0.0,+1.6){7}{\circle*1}
  \end{picture}}^2
  \,+\,
  \frac12\;\raise-5pt\hbox{\begin{picture}(16,16)
   \put( 0, 0){\circle2}
   \put(16, 0){\circle2}
   \put( 0,16){\circle2}
   \put(16,16){\circle2}
   \put( 0, 1){\line( 0,+1){14}}
   \put(16, 1){\line( 0,+1){14}}
   \put( 1, 0){\line(+1, 0){14}}
   \put( 1,16){\line(+1, 0){14}}
   \multiput( 1.1, 1.7)( 0.0,+1.8){8}{\circle*1}
   \multiput(14.9, 1.7)( 0.0,+1.8){8}{\circle*1}
  \end{picture}}
  \,-\,
  \frac14\;\raise-5pt\hbox{\begin{picture}(16,16)
   \put( 0, 0){\circle2}
   \put(16, 0){\circle2}
   \put( 0,16){\circle2}
   \put(16,16){\circle2}
   \put( 0, 1){\line( 0,+1){14}}
   \put(16, 1){\line( 0,+1){14}}
   \put( 1, 0){\line(+1, 0){14}}
   \put( 1,16){\line(+1, 0){14}}
   \multiput( 1.4, 1.4)(+1.1,+1.1){13}{\circle*1}
   \multiput( 1.4,14.6)(+1.1,-1.1){13}{\circle*1}
  \end{picture}}\ .
 $$
 Specializing to dimension $m\,=\,4$ we observe that the actual second
 moment polynomial equals $\Psi_2\,=\,\frac1{180}\,\Psi^\circ_2$ or
 equivalently $15\,\Psi_2\,=\,\frac1{12}\,\Psi^\circ_2$, compare Lemma
 \ref{msc}, moreover we can decompose the Ricci tensor into its trace free
 and scalar part $\Ric\,=\,\Ric_\circ\,+\,\frac\kappa4g$ to find:
 $$
  \Big[\,\pf_2\,-\,15\,\Psi_2\,\Big](\,R\,)
  \;\;=\;\;
  \Big[\;\frac1{48}\raise-4pt\hbox{\begin{picture}(10,14)
   \put( 5, 0){\circle2}
   \put( 5,14){\circle2}
   \qbezier( 4.0, 0.0)( 1.0, 0.0)( 1.0, 7.0)
   \qbezier( 4.0,14.0)( 1.0,14.0)( 1.0, 7.0)
   \qbezier( 6.0, 0.0)( 9.0, 0.0)( 9.0, 7.0)
   \qbezier( 6.0,14.0)( 9.0,14.0)( 9.0, 7.0)
   \multiput( 5.0, 2.2)( 0.0,+1.6){7}{\circle*1}
  \end{picture}}^2
  \;-\;
  \frac16\;\raise-5pt\hbox{\begin{picture}(16,16)
   \put( 0, 0){\circle2}
   \put(16, 0){\circle2}
   \put( 0,16){\circle2}
   \put(16,16){\circle2}
   \put( 0, 1){\line( 0,+1){14}}
   \put(16, 1){\line( 0,+1){14}}
   \put( 1, 0){\line(+1, 0){14}}
   \put( 1,16){\line(+1, 0){14}}
   \multiput( 1.1, 1.7)( 0.0,+1.8){8}{\circle*1}
   \multiput(14.9, 1.7)( 0.0,+1.8){8}{\circle*1}
  \end{picture}}\;\Big](\,R\,)
  \;\;=\;\;
  -\;\frac43\,g_{\S^2T^*}(\,\Ric_\circ,\Ric_\circ\,)
  \;-\;\frac{\kappa^2}{12}\ .
 $$
 Integrating this identity over a compact $4$--dimensional Riemannian
 manifold $M$ we obtain
 \begin{equation}\label{htho}
  (\,2\pi\,)^2\,\chi(M)\;+\;\frac43\,|\!|\,\Ric_\circ|\!|^2_{\S^2T^*}
  \;\;=\;\;
  \int_M\frac{\kappa^2}{48}\,|\,\vol_g\,|\;+\;
  15\,\int_M\Big(\;\Psi_2(R)\,-\,\frac{\kappa^2}{144}\;\Big)\,|\,\vol_g\,|
 \end{equation}
 due to the formulation (\ref{cgb}) of the Theorem of Chern--Gauß--Bonnet.
 For a compact Riemannian manifold $M$ the right hand side is strictly positive
 unless $M$ is flat, because the functions $\frac\kappa{12}$ and $\Psi_2(R)
 \,-\,(\frac\kappa{12})^2\,\geq\,0$ are the pointwise expected value and the
 pointwise variance respectively of the sectional curvature considered as a
 random variable. The resulting inequality $4\,\pi^2\,\chi(M)\,+\,\frac43\,
 |\!|\,\Ric_\circ\,|\!|_{\S^2T^*}\,\geq\,0$ is a weak version of the
 Hitchin--Thorpe inequality for oriented manifolds $\cite{h}$, in particular
 every $4$--dimensional Einstein manifold has strictly positive Euler
 characteristic unless it is flat \cite{b}. In passing we recall that
 the signature is not a stable curvature invariant in the sense of this
 article.

 \pfill
 In order to generalize the Hitchin--Thorpe inequality to Einstein manifolds
 of higher dimensions we will make use of the standard curvature term $q(R)$
 defined in \cite{sw}, albeit with a slightly different normalization. The
 letter $q$ in our notation does not stand for quadratic, in contrast to the
 notation $Q(R)$ used in equation (\ref{etc}), but refers instead to the
 so--called quantization map from the symmetric to the universal enveloping
 algebra of a Lie algebra:
 $$
  q:\;\;\S^{\leq\bullet}\mathfrak{g}\;\stackrel\cong\longrightarrow\;
  \U^{\leq\bullet}\mathfrak{g},
  \qquad\X^r\;\longmapsto\;\X^r\ .
 $$
 Specializing to the orthogonal Lie algebra $\so(\,T,\,g\,)$ of a euclidean
 vector space $T$ we can convert every algebraic curvature tensor $R$ over
 $T$ via $q$ linearly into various endomorphisms:

 \begin{Definition}[Standard Curvature Term]
 \hfill\label{sct}\break
  For the orthogonal Lie algebra $\so(\,T,\,g\,)$ of skew symmetric
  endomorphisms of a euclidean vector space $T$ the quantization map
  $q:\,\S^{\leq\bullet}\so(\,T,\,g\,)\longrightarrow\U^{\leq
  \bullet}\so(\,T,\,g\,)$ can be applied to every algebraic curvature
  tensor $R\,\in\,\Curv^-T$ over $T$ considered in analogy to equation
  (\ref{2v})
  $$
   R
   \;\;=\;\;
   \frac12\,\sum_{\mu,\,\nu\,=\,1}^m
   (\,dE^\up_\mu\,\wedge\,dE^\up_\nu\,)\,\otimes\,R_{E_\mu,\,E_\nu}
   \;\;\stackrel!=\;\;
   \frac14\,\sum_{\mu,\,\nu\,=\,1}^m
   (\,dE^\up_\mu\,\wedge\,dE^\up_\nu\,)\,\cdot\,R_{E_\mu,\,E_\nu}
  $$
  as an element of the symmetric square $\S^2\so(\,T,\,g\,)$. Its image
  $q(\,R\,)\,\in\,\U^{\leq2}\so(\,T,\,g\,)$ becomes in turn an endomorphism
  in every representation $\star$ of the Lie algebra $\so(\,T,\,g\,)$:
  $$
   q(\,R\,)\;\star
   \;\;:=\;\;
   \frac14\,\sum_{\mu,\,\nu\,=\,1}^m
   (\,dE^\up_\mu\wedge dE^\up_\nu\,)\;\star\;R_{E_\mu,\,E_\nu}\;\star\ .
  $$
 \end{Definition}

 \noindent
 In this definition of the standard curvature term the scalar factor $\frac14$
 is the proper choice, in difference to the rather unmotivated scalar factor
 $\frac12$ used in \cite{sw}. For this reason the standard curvature term
 equals half the Ricci endomorphism, the symmetric endomorphism of $T$
 corresponding to the Ricci tensor, in the defining representation of the
 Lie algebra $\so(\,T,\,g\,)$:
 $$
  q(R)\,\star\,X
  \;\;=\;\;
  \frac12\,\sum_{\mu,\,\nu\,=\,1}^m
  g(\,dE^\up_\mu,\,R_{E_\mu,\,E_\nu}X\,)\,dE^\up_\nu
  \;\;=\;\;
  \frac12\,\sum_{\nu\,=\,1}^m
  \Ric(E_\nu,X)\,dE^\up_\nu
  \;\;=:\;\;
  \frac12\,\Ric\,X\ .
 $$
 On the bivector representations $\LP^2T$ and $\S^2T$ the standard curvature
 term $q(R)$ acts essentially as the so called curvature operator and its
 symmetric counterpart:
 $$
  \begin{array}{ccccl}
   \K^\mathrm{alt}:&\LP^2T&\longrightarrow&\LP^2T\;,\quad
   & X\wedge Y\;\longmapsto\;-\;R_{X,\,Y}
   \\[4pt]
   \K^\mathrm{sym}:&\S^2T&\longrightarrow&\S^2T,\quad
   & X\,\cdot\,Y\;\longmapsto\;+\,\Sec_{X,\,Y}\ .
  \end{array}
 $$
 In fact we find say for the adjoint or alternating bivector representation
 $\so(\,T,\,g\,)\,\cong\,\LP^2T$
 \begin{eqnarray*}
  \lefteqn{q(\;R\;)\,\star\,(\,X\,\wedge\,Y\,)}
  &&
  \\
  &=&
  \Der_{\frac12\Ric}(\,X \wedge Y\,)\,+\,\frac12\sum_{\mu,\,\nu\,=\,1}^m
  \Big(\,g(dE^\up_\mu,X)dE^\up_\nu \wedge R_{E_\mu,\,E_\nu}Y\,+\,
  R_{E_\mu,\,E_\nu}X\wedge g(dE^\up_\mu,Y)dE^\up_\nu\,\Big)
  \\
  &=&
  \Der_{\frac12\Ric}(\,X\wedge Y\,)\,+\,\frac12\,\sum_{\nu\,=\,1}^m
  dE^\up_\nu\,\wedge\,\Big(\;R_{X,\,E_\nu}Y\;-\;R_{Y,\,E_\nu}X\;\Big)
  \\
  &=&
  \frac12\,\Der_\Ric(\,X\,\wedge\,Y\,)
  \;-\;\K^\mathrm{alt}(\,X\,\wedge\,Y\,)\ ,
 \end{eqnarray*}
 where $\Der_{\frac12\Ric}$ denotes the derivation extension of the Ricci
 endomorphism of $T$ to an endomorphism of $\LP^2T$ via $\Der_{\frac12\Ric}
 (\,X\wedge Y\,)\,:=\,\frac12\,\Ric\,X\wedge Y\,+\,\frac12\,X\wedge\Ric\,Y$.
 Replacing the first Bianchi identity $R_{X,\,E_\nu}Y-R_{Y,\,E_\nu}X\,=\,
 R_{X,\,Y}E_\nu$ used in this argument by the variant $R_{X,\,E_\nu}Y+
 R_{Y,\,E_\nu}X\,=\,-\frac12\,\Sec(X,Y;E_\nu,\,\cdot\,)^\up$ of the definition
 of the sectional curvature tensor $\Sec\,:=\,\Phi^+R$ we find analogously
 $q(\,R\,)\,\star\,=\,\frac12\,\Der_\Ric\,-\,\frac12\,\K^\mathrm{sym}$ on
 the symmetric bivector representation $\S^2T$. Let us now study the action
 of $q(\,R\,)$ on $R\,\in\,\Curv^-T$ itself:
 
 \begin{Lemma}[Standard Curvature Term Polynomial]
 \hfill\label{cubic}\break
  Consider an algebraic curvature tensor $R\,\in\,\Curv^-T$ on a euclidean
  vector space $T$ with scalar product $g$. The scalar product of $q(R)\,
  \star\,R$ with $R$ equals the stable cubic polynomial:
  $$
   \Theta_3(\,R\,)
   \;\;:=\;\;
   g_{\LP^2T^*\,\otimes\,\LP^2T^*}(\;q(\,R\,)\,\star\,R,\;R\;)
   \;\;=\;\;
   \frac1{12}\;\left[
   \frac16\;
   \raise-11pt\hbox{\begin{picture}(30,30)
    \put(15, 2){\circle2}
    \put(15,28){\circle2}
    \put( 3, 8){\circle2}
    \put(27, 8){\circle2}
    \put( 3,22){\circle2}
    \put(27,22){\circle2}
    \put( 3.0, 9.0){\line( 0,+1){12}}
    \put(27.0, 9.0){\line( 0,+1){12}}
    \put(15.9, 2.5){\line(+2,+1){10}}
    \put(14.1, 2.5){\line(-2,+1){10}}
    \put(15.9,27.5){\line(+2,-1){10}}
    \put(14.1,27.5){\line(-2,-1){10}}
    \multiput(15.0, 4.0)( 0.0,+2.0){12}{\circle*1}
    \multiput( 4.7, 9.0)(+1.9,+1.1){12}{\circle*1}
    \multiput(25.3, 9.0)(-1.9,+1.1){12}{\circle*1}
   \end{picture}}
   \;-\;
   \frac23\;
   \raise-11pt\hbox{\begin{picture}(30,30)
    \put(15, 2){\circle2}
    \put(15,28){\circle2}
    \put( 3, 8){\circle2}
    \put(27, 8){\circle2}
    \put( 3,22){\circle2}
    \put(27,22){\circle2}
    \put( 3.0, 9.0){\line( 0,+1){12}}
    \put(27.0, 9.0){\line( 0,+1){12}}
    \put(15.9, 2.5){\line(+2,+1){10}}
    \put(14.1, 2.5){\line(-2,+1){10}}
    \put(15.9,27.5){\line(+2,-1){10}}
    \put(14.1,27.5){\line(-2,-1){10}}
    \multiput(15.0, 4.0)( 0.0,+2.0){12}{\circle*1}
    \multiput( 4.9, 8.2)(+2.0, 0.0){11}{\circle*1}
    \multiput( 4.9,21.8)(+2.0, 0.0){11}{\circle*1}
   \end{picture}}
   \;+\;
   \frac12\;
   \raise-11pt\hbox{\begin{picture}(30,30)
    \put(15, 2){\circle2}
    \put(15,28){\circle2}
    \put( 3, 8){\circle2}
    \put(27, 8){\circle2}
    \put( 3,22){\circle2}
    \put(27,22){\circle2}
    \put( 3.0, 9.0){\line( 0,+1){12}}
    \put(27.0, 9.0){\line( 0,+1){12}}
    \put(15.9, 2.5){\line(+2,+1){10}}
    \put(14.1, 2.5){\line(-2,+1){10}}
    \put(15.9,27.5){\line(+2,-1){10}}
    \put(14.1,27.5){\line(-2,-1){10}}
    \multiput( 4.9,21.9)(+1.5,+0.75){7}{\circle*1}
    \multiput(16.1, 3.8)(+0.9,+1.5){12}{\circle*1}
    \multiput( 5.0, 8.0)(+2.0, 0.0){11}{\circle*1}
   \end{picture}}
   \right](\;R\;)\ .
  $$
 \end{Lemma}

 \proof
 Without discussing the details of the construction of such an algebra we
 consider an extended version of the algebra $\A^\bullet$ of colored trivalent
 graphs similar to the extended graph algebra $\A^\bullet_\ext$ discussed at
 the end of Section \ref{graphs}, in which a pair of vertices may be connected
 by a single or a double red edge with the associated additional Feynman rule:
 \begin{equation}\label{feyn3}
  \raise-19pt\hbox{\begin{picture}(54,45)
   \put(28,13){\circle2}
   \put(28,33){\circle2}
   \put(10, 4){\line(+2,+1){17}}
   \put(10,42){\line(+2,-1){17}}
   \put(46, 4){\line(-2,+1){17}}
   \put(46,42){\line(-2,-1){17}}
   \multiput(27.3,15.0)( 0,+2){9}{\circle*1}
   \multiput(28.7,15.0)( 0,+2){9}{\circle*1}
   \put( 0,39){$X$}
   \put(48,39){$Y$}
   \put( 0, 0){$V$}
   \put(48, 0){$U$}
  \end{picture}}
  \;\;\widehat=\;\;
  \bigg(\,q(R)\,\star\,\Sec\,\bigg)(\;X,\;Y;\;U,\;V\;)\ .
 \end{equation}
 In order to expand elements of such an extended graph algebra into elements
 of the algebra $\A^\bullet$ of colored trivalent graphs we expand the
 additional quartic interaction into the sum 
 \begin{eqnarray*}
  \lefteqn{\bigg(\,q(R)\,\star\,\Sec\,\bigg)(\,X,\,Y;\,U,\,V\,)}
  &&
  \\
  &=&
  \!\!\!\!\sum_{\mu,\,\nu,\,\alpha,\,\beta\,=\,1}^m
  \Sec(E_\mu,E_\nu;E_\alpha,E_\beta)\,
  g\Big(q(R)\star(X\otimes Y\otimes U\otimes V),\,
  dE^\up_\mu\otimes dE^\up_\nu\otimes dE^\up_\alpha\otimes dE^\up_\beta\Big)
 \end{eqnarray*}
 over an arbitrary pair of dual bases $\{\,E_\mu\,\}$ and $\{\,dE_\mu\,\}$.
 Representations of Lie algebras extend in general as derivations to the
 tensor algebra $\bigotimes T$ associated to $T$, hence the quadratic element
 $q(R)\,\in\,\U^{\leq2}\so(\,T,\,g\,)$ can be seen as a second order
 differential operator in the sense that its action on $\bigotimes T$
 is determined by its action on the tensor square $\bigotimes^2T$ with
 \begin{eqnarray*}
  \lefteqn{g(\;q(R)\,\star\,(\,X\,\otimes\,Y\,),\;U\,\otimes\,V\;)}
  &&
  \\
  &=&
  g(\;\Der_{\frac12\,\Ric}(\,X\,\otimes\,Y\,),\;U\,\otimes\,V\;)
  \;+\;
  \sum_{\mu,\,\nu\,=\,1}^mg(\;dE_\mu(X)\,dE^\up_\nu\,\otimes\,R_{E_\mu,\,
   E_\nu}Y,\;U\,\otimes\,V\;)
  \\
  &=&
  g(\;\Der_{\frac12\,\Ric}(\,X\,\otimes\,Y\,),\;U\,\otimes\,V\;)
  \;-\;\frac16\,\Sec(\,X,\,Y;\,U,\,V\,)\;+\;\frac16\,\Sec(\,X,\,V;\,Y,\,U\,)
 \end{eqnarray*}
 for all $X,\,Y,\,U,\,V\,\in\,T$ using $R\,=\,\Phi^-\Sec$. In consequence
 the full expansion of a double red edge (\ref{feyn3}) into ordinary colored
 trivalent graphs results in a sum of 16 terms grouped into two sums of 4
 and 12 terms respectively. In the first of these two sums we replace the
 double by a single red edge and in addition each of the four adjacent
 black flags in turn by
 $$
  \raise-11pt\hbox{\begin{picture}(20,30)
   \put(10,0){\line(0,+1){30}}
  \end{picture}}
  \;\;\rightsquigarrow\;\;
  -\;\frac14\;\raise-11pt\hbox{\begin{picture}(20,30)
   \put( 5, 6){\circle2}
   \put( 5,24){\circle2}
   \put( 5, 0){\line(0,+1){5}}
   \put( 5,30){\line(0,-1){5}}
   \multiput( 5, 8)( 0,+2){8}{\circle*1}
   \qbezier( 6, 6)(12, 6)(12,15)
   \qbezier( 6,24)(12,24)(12,15)
  \end{picture}}
 $$
 in order to implement $\Der_{\frac12\Ric}$ via equation (\ref{ihex}). In the
 second sum we replace similarly the double by a single red edge and
 moreover each of the six pairs of adjacent flags in turn by:
 $$
  \raise-11pt\hbox{\begin{picture}(20,30)
   \put( 3, 0){\line(0,+1){30}}
   \put(17, 0){\line(0,+1){30}}
  \end{picture}}
  \;\;\rightsquigarrow\;\;
  -\;\frac16\;\raise-11pt\hbox{\begin{picture}(20,30)
   \put(10, 6){\circle2}
   \put(10,24){\circle2}
   \multiput(10, 8)( 0,+2){8}{\circle*1}
   \qbezier( 3, 0)( 3, 5)( 9, 6)
   \qbezier(17, 0)(17, 5)(11, 6)
   \qbezier( 3,30)( 3,25)( 9,24)
   \qbezier(17,30)(17,25)(11,24)
  \end{picture}}
  \;+\;\frac16\;\raise-11pt\hbox{\begin{picture}(30,30)
   \put( 5,15){\circle2}
   \put(15,15){\circle2}
   \multiput( 7,15)(+2, 0){4}{\circle*1}
   \qbezier( 3, 0)( 3,14)( 5,14)
   \qbezier(17, 0)(17,14)(15,14)
   \qbezier( 3,30)( 6,22)(15,16)
   \qbezier(17,30)(14,22)( 5,16)
  \end{picture}}\ .
 $$
 Implementing this expansion in the special case relevant to Lemma
 \ref{cubic} we find without effort:
 \begin{eqnarray*}
  \frac14\;\raise-13pt\hbox{\begin{picture}(30,30)
   \put( 4, 4){\circle2}
   \put( 4,26){\circle2}
   \put(26, 4){\circle2}
   \put(26,26){\circle2}
   \multiput( 3.3, 6.0)( 0,+2){10}{\circle*1}
   \multiput( 4.7, 6.0)( 0,+2){10}{\circle*1}
   \multiput(26.0, 6.0)( 0,+2){10}{\circle*1}
   \put(4.9, 3.5){\line(+1, 0){20.2}}
   \put(4.9, 4.5){\line(+1, 0){20.2}}
   \put(4.9,25.5){\line(+1, 0){20.2}}
   \put(4.9,26.5){\line(+1, 0){20.2}}
  \end{picture}}
  &\rightsquigarrow&
  -\,\frac14\cdot\frac44\,
  \raise-11pt\hbox{\begin{picture}(30,30)
   \put( 6, 3){\circle2}
   \put(24, 3){\circle2}
   \put( 2,20){\circle2}
   \put(28,20){\circle2}
   \put( 9,23){\circle2}
   \put(21,23){\circle2}
   \put( 6.9, 3.6){\line(+1, 0){16.5}}
   \put( 6.9, 2.4){\line(+1, 0){16.5}}
   \put( 3.0,20.0){\line(+1, 0){24.0}}
   \qbezier( 2.8,20.7)( 2.8,20.7)( 8.1,22.6)
   \qbezier(27.2,20.7)(27.2,20.7)(21.9,22.6)
   \qbezier( 9.7,23.7)(11.0,27.0)(15.0,27.0)
   \qbezier(20.3,23.7)(19.0,27.0)(15.0,27.0)
   \multiput( 5.5, 5.0)(-0.4,+1.6){9}{\circle*1}
   \multiput(24.5, 5.0)(+0.4,+1.6){9}{\circle*1}
   \multiput(11.2,23.0)(+1.9, 0.0){5}{\circle*1}
  \end{picture}}
  \;-\,\frac16\cdot\frac24\,
  \raise-11pt\hbox{\begin{picture}(30,30)
   \put( 6, 3){\circle2}
   \put(24, 3){\circle2}
   \put( 2,20){\circle2}
   \put(28,20){\circle2}
   \put(10,20){\circle2}
   \put(20,20){\circle2}
   \put( 6.9, 3.6){\line(+1, 0){16.5}}
   \put( 6.9, 2.4){\line(+1, 0){16.5}}
   \put( 3.0,20.6){\line(+1, 0){6.2}}
   \put( 3.0,19.4){\line(+1, 0){6.2}}
   \put(27.1,20.6){\line(-1, 0){6.2}}
   \put(27.1,19.4){\line(-1, 0){6.2}}
   \multiput( 5.5, 5.0)(-0.4,+1.6){9}{\circle*1}
   \multiput(24.5, 5.0)(+0.4,+1.6){9}{\circle*1}
   \multiput(12.0,20.0)(+1.9, 0.0){4}{\circle*1}
  \end{picture}}
  \;+\,\frac16\cdot\frac24\,
  \raise-11pt\hbox{\begin{picture}(30,30)
   \put( 6, 3){\circle2}
   \put(24, 3){\circle2}
   \put( 2,20){\circle2}
   \put(28,20){\circle2}
   \put(15,26){\circle2}
   \put(15,15){\circle2}
   \put( 6.9, 3.6){\line(+1, 0){16.5}}
   \put( 6.9, 2.4){\line(+1, 0){16.5}}
   \put( 2.8,20.5){\line(+2,+1){11.0}}
   \put(27.2,20.5){\line(-2,+1){11.0}}
   \put( 2.7,19.2){\line(+3,-1){11.0}}
   \put(27.3,19.2){\line(-3,-1){11.0}}
   \multiput( 5.5, 5.0)(-0.4,+1.6){9}{\circle*1}
   \multiput(24.5, 5.0)(+0.4,+1.6){9}{\circle*1}
   \multiput(15.0,17.1)( 0.0,+1.7){5}{\circle*1}
  \end{picture}}
  \;-\,\frac16\cdot\frac44\,
  \raise-11pt\hbox{\begin{picture}(30,30)
   \put( 3, 3){\circle2}
   \put(27, 3){\circle2}
   \put( 3,27){\circle2}
   \put(27,27){\circle2}
   \put(10,15){\circle2}
   \put(20,15){\circle2}
   \put( 4.0, 2.4){\line(+1, 0){22.2}}
   \put( 4.0,27.6){\line(+1, 0){22.2}}
   \qbezier( 3.5, 3.9)( 6.3, 8.5)( 9.2,14.2)
   \qbezier( 3.5,26.1)( 6.3,21.5)( 9.2,15.8)
   \qbezier(26.5, 3.9)(23.7, 8.5)(20.8,14.2)
   \qbezier(26.5,26.1)(23.7,21.5)(20.8,15.8)
   \multiput( 3.0, 5.0)( 0,+2){11}{\circle*1}
   \multiput(27.0, 5.0)( 0,+2){11}{\circle*1}
   \multiput(12.0,15.0)(+2, 0){4}{\circle*1}
  \end{picture}}
  \;+\,\frac16\cdot\frac44\,
  \raise-11pt\hbox{\begin{picture}(30,30)
   \put( 3, 3){\circle2}
   \put(27, 3){\circle2}
   \put( 3,27){\circle2}
   \put(27,27){\circle2}
   \put(15,10){\circle2}
   \put(15,20){\circle2}
   \put( 4.0, 2.4){\line(+1, 0){22.2}}
   \put( 4.0,27.6){\line(+1, 0){22.2}}
   \qbezier( 3.5, 3.9)( 6.0,10.5)(14.0,20.0)
   \qbezier( 3.5,26.1)( 6.0,19.5)(14.0,10.0)
   \qbezier(26.2, 3.3)(20.2, 6.1)(15.8, 9.6)
   \qbezier(26.2,26.7)(20.2,23.9)(15.8,20.4)
   \multiput( 3.0, 5.0)( 0,+2){11}{\circle*1}
   \multiput(27.0, 5.0)( 0,+2){11}{\circle*1}
   \multiput(15.0,12.0)( 0,+2){4}{\circle*1}
  \end{picture}}\ .
 \end{eqnarray*}
 Interpreting the left hand side and simplifying the right hand side using
 the congruences
 $$
  \raise-11pt\hbox{\begin{picture}(30,30)
   \put( 6, 3){\circle2}
   \put(24, 3){\circle2}
   \put( 2,20){\circle2}
   \put(28,20){\circle2}
   \put( 9,23){\circle2}
   \put(21,23){\circle2}
   \put( 6.9, 3.6){\line(+1, 0){16.5}}
   \put( 6.9, 2.4){\line(+1, 0){16.5}}
   \put( 3.0,20.0){\line(+1, 0){24.0}}
   \qbezier( 2.8,20.7)( 2.8,20.7)( 8.1,22.6)
   \qbezier(27.2,20.7)(27.2,20.7)(21.9,22.6)
   \qbezier( 9.7,23.7)(11.0,27.0)(15.0,27.0)
   \qbezier(20.3,23.7)(19.0,27.0)(15.0,27.0)
   \multiput( 5.5, 5.0)(-0.4,+1.6){9}{\circle*1}
   \multiput(24.5, 5.0)(+0.4,+1.6){9}{\circle*1}
   \multiput(11.2,23.0)(+1.9, 0.0){5}{\circle*1}
  \end{picture}}
  \;\;\equiv\;\;
  -\;2\;\raise-11pt\hbox{\begin{picture}(30,30)
   \put(15, 2){\circle2}
   \put(15,28){\circle2}
   \put( 3, 8){\circle2}
   \put(27, 8){\circle2}
   \put( 3,22){\circle2}
   \put(27,22){\circle2}
   \put( 3.0, 9.0){\line( 0,+1){12}}
   \put(27.0, 9.0){\line( 0,+1){12}}
   \put(15.9, 2.5){\line(+2,+1){10}}
   \put(14.1, 2.5){\line(-2,+1){10}}
   \put(15.9,27.5){\line(+2,-1){10}}
   \put(14.1,27.5){\line(-2,-1){10}}
   \multiput( 4.9,21.9)(+1.5,+0.75){7}{\circle*1}
   \multiput(16.1, 3.8)(+0.9,+1.5){12}{\circle*1}
   \multiput( 5.0, 8.0)(+2.0, 0.0){11}{\circle*1}
  \end{picture}}
  \qquad\qquad
  \raise-11pt\hbox{\begin{picture}(30,30)
   \put(15, 2){\circle2}
   \put(15,28){\circle2}
   \put( 3, 8){\circle2}
   \put(27, 8){\circle2}
   \put( 3,22){\circle2}
   \put(27,22){\circle2}
   \put( 2.5, 8.9){\line( 0,+1){12.1}}
   \put( 3.5, 8.9){\line( 0,+1){12.1}}
   \put(15.8, 2.9){\line(+2,+1){10.2}}
   \put(16.3, 2.1){\line(+2,+1){10.2}}
   \put(15.8,27.1){\line(+2,-1){10.2}}
   \put(16.3,27.9){\line(+2,-1){10.2}}
   \multiput( 4.7,22.8)(+1.7,+0.85){6}{\circle*1}
   \multiput( 4.7, 7.2)(+1.7,-0.85){6}{\circle*1}
   \multiput(27.0,10.0 )( 0.0,+2.0){6}{\circle*1}
  \end{picture}}
  \;\;\equiv\;\;
  -\;2\;\raise-11pt\hbox{\begin{picture}(30,30)
   \put( 6, 3){\circle2}
   \put(24, 3){\circle2}
   \put( 2,20){\circle2}
   \put(28,20){\circle2}
   \put(15,26){\circle2}
   \put(15,15){\circle2}
   \put( 6.9, 3.6){\line(+1, 0){16.5}}
   \put( 6.9, 2.4){\line(+1, 0){16.5}}
   \put( 2.8,20.5){\line(+2,+1){11.0}}
   \put(27.2,20.5){\line(-2,+1){11.0}}
   \put( 2.7,19.2){\line(+3,-1){11.0}}
   \put(27.3,19.2){\line(-3,-1){11.0}}
   \multiput( 5.5, 5.0)(-0.4,+1.6){9}{\circle*1}
   \multiput(24.5, 5.0)(+0.4,+1.6){9}{\circle*1}
   \multiput(15.0,17.1)( 0.0,+1.7){5}{\circle*1}
  \end{picture}}
  \;\;\equiv\;\;
  4\;\raise-11pt\hbox{\begin{picture}(30,30)
   \put(15, 2){\circle2}
   \put(15,28){\circle2}
   \put( 3, 8){\circle2}
   \put(27, 8){\circle2}
   \put( 3,22){\circle2}
   \put(27,22){\circle2}
   \put( 3.0, 9.0){\line( 0,+1){12}}
   \put(27.0, 9.0){\line( 0,+1){12}}
   \put(15.9, 2.5){\line(+2,+1){10}}
   \put(14.1, 2.5){\line(-2,+1){10}}
   \put(15.9,27.5){\line(+2,-1){10}}
   \put(14.1,27.5){\line(-2,-1){10}}
   \multiput(15.0, 4.0)( 0.0,+2.0){12}{\circle*1}
   \multiput( 4.9, 8.2)(+2.0, 0.0){11}{\circle*1}
   \multiput( 4.9,21.8)(+2.0, 0.0){11}{\circle*1}
  \end{picture}}
 $$
 modulo the ideal of IHX--relations we eventually obtain the formula:
 \begin{eqnarray*}
  \lefteqn{g_{\S^2T^*\otimes\S^2T^*}(\;q(\,R\,)\star\Sec,\;\Sec\;)}
  \qquad
  &&
  \\[4pt]
  &=&
  \left[\;\frac14\;\raise-13pt\hbox{\begin{picture}(30,30)
   \put( 4, 4){\circle2}
   \put( 4,26){\circle2}
   \put(26, 4){\circle2}
   \put(26,26){\circle2}
   \multiput( 3.3, 6.0)( 0,+2){10}{\circle*1}
   \multiput( 4.7, 6.0)( 0,+2){10}{\circle*1}
   \multiput(26.0, 6.0)( 0,+2){10}{\circle*1}
   \put(4.9, 3.5){\line(+1, 0){20.2}}
   \put(4.9, 4.5){\line(+1, 0){20.2}}
   \put(4.9,25.5){\line(+1, 0){20.2}}
   \put(4.9,26.5){\line(+1, 0){20.2}}
  \end{picture}}\;\right](\;R\;)
  \;\;=\;\;
  \left[\;\frac16\;\raise-11pt\hbox{\begin{picture}(30,30)
   \put(15, 2){\circle2}
   \put(15,28){\circle2}
   \put( 3, 8){\circle2}
   \put(27, 8){\circle2}
   \put( 3,22){\circle2}
   \put(27,22){\circle2}
   \put( 3.0, 9.0){\line( 0,+1){12}}
   \put(27.0, 9.0){\line( 0,+1){12}}
   \put(15.9, 2.5){\line(+2,+1){10}}
   \put(14.1, 2.5){\line(-2,+1){10}}
   \put(15.9,27.5){\line(+2,-1){10}}
   \put(14.1,27.5){\line(-2,-1){10}}
   \multiput(15.0, 4.0)( 0.0,+2.0){12}{\circle*1}
   \multiput( 4.7, 9.0)(+1.9,+1.1){12}{\circle*1}
   \multiput(25.3, 9.0)(-1.9,+1.1){12}{\circle*1}
  \end{picture}}
  \;-\;
  \frac23\;\raise-11pt\hbox{\begin{picture}(30,30)
   \put(15, 2){\circle2}
   \put(15,28){\circle2}
   \put( 3, 8){\circle2}
   \put(27, 8){\circle2}
   \put( 3,22){\circle2}
   \put(27,22){\circle2}
   \put( 3.0, 9.0){\line( 0,+1){12}}
   \put(27.0, 9.0){\line( 0,+1){12}}
   \put(15.9, 2.5){\line(+2,+1){10}}
   \put(14.1, 2.5){\line(-2,+1){10}}
   \put(15.9,27.5){\line(+2,-1){10}}
   \put(14.1,27.5){\line(-2,-1){10}}
   \multiput(15.0, 4.0)( 0.0,+2.0){12}{\circle*1}
   \multiput( 4.9, 8.2)(+2.0, 0.0){11}{\circle*1}
   \multiput( 4.9,21.8)(+2.0, 0.0){11}{\circle*1}
  \end{picture}}
  \;+\;
  \frac12\;\raise-11pt\hbox{\begin{picture}(30,30)
   \put(15, 2){\circle2}
   \put(15,28){\circle2}
   \put( 3, 8){\circle2}
   \put(27, 8){\circle2}
   \put( 3,22){\circle2}
   \put(27,22){\circle2}
   \put( 3.0, 9.0){\line( 0,+1){12}}
   \put(27.0, 9.0){\line( 0,+1){12}}
   \put(15.9, 2.5){\line(+2,+1){10}}
   \put(14.1, 2.5){\line(-2,+1){10}}
   \put(15.9,27.5){\line(+2,-1){10}}
   \put(14.1,27.5){\line(-2,-1){10}}
   \multiput( 4.9,21.9)(+1.5,+0.75){7}{\circle*1}
   \multiput(16.1, 3.8)(+0.9,+1.5){12}{\circle*1}
   \multiput( 5.0, 8.0)(+2.0, 0.0){11}{\circle*1}
  \end{picture}}\;\right](\;R\;)\ .  
 \end{eqnarray*}
 In light of the isometry equation (\ref{eqnorm}) this result is
 equivalent to the stipulated formula for the stable cubic polynomial
 $g_{\LP^2T^*\otimes\LP^2T^*}(\,q(R)\star R,\,R\,)$, after all the mutually
 inverse isomorphisms $\Phi^+$ and $\Phi^-$ of Section \ref{tensors} are
 equivariant under the action of the orthogonal group and thus commute
 $q(R)\,\star\,\Phi^-\Sec\,=\,\Phi^-(\,q(R)\,\star\,\Sec\,)$ with the
 standard curvature term.
 \qed

 \pfill
 The stable cubic curvature invariant $\Theta_3$ of Lemma \ref{cubic} arises
 naturally in the study of the curvature tensor $R\,\in\,\Gamma(\,\Curv^-TM\,)$
 of Riemannian manifolds $M$ with parallel Ricci curvature $\Ric\,\in\,
 \Gamma(\,\S^2T^*M\,)$, a class of Riemannian manifolds slightly larger
 than the class of Einstein manifolds of dimension $m\,\geq\,3$ \cite{besse}.
 In order to relate $\Theta_3$ to the Ricci curvature we compose the
 symmetrized second covariant derivative $\nabla^{[2]}_{X,\,Y}\,:=\,\frac12
 (\nabla^2_{X,\,Y}+\nabla^2_{Y,\,X})$ with the Nomizu--Kulkarni product
 $\times$ of equation (\ref{nk}) to obtain a second order differential
 operator:
 $$
  \Cross:\;\;\Gamma(\,\S^2T^*M\,)
  \;\stackrel{\nabla^{[2]}}\longrightarrow\;\Gamma(\,\S^2T^*M\otimes\S^2T^*M\,)
  \;\stackrel\times\longrightarrow\;\Gamma(\,\Curv^-TM\,)\ .
 $$
 This cross operator is a fundamental differential operator in Riemannian
 geometry, because it is in essence the linearization of the second order
 non--linear differential operator, which sends a Riemannian metric $g\,\in\,
 \Gamma(\,\S^2_{\mathrm{reg}}T^*M\,)$ to its curvature tensor $R\,\in\,
 \Gamma(\,\Curv^-TM\,)$, more precisely the curvature tensor varies under
 an arbitrary variation of the metric according to
 \begin{equation}\label{varr}
  \delta R
  \;\;=\;\;
  \frac12\;\Cross\;\delta g\;+\;\frac14\;\Der_{[\,\delta g\,]}R\ ,
 \end{equation} 
 where $[\,\delta g\,]$ denotes the symmetric endomorphism field defined by
 $g([\,\delta g\,]X,Y)\,:=\,\delta g(X,Y)$.
   
 \begin{Lemma}[Laplacian of Curvature Tensor \cite{sw}]
 \hfill\label{pde}\break
  The curvature tensor $R\,\in\,\Gamma(\,\Curv^-TM\,)$ of every Riemannian
  manifold $M$ satisfies:
  $$
   \nabla^*\nabla R\;+\;q(\,R\,)\,\star\,R
   \;\;=\;\;
   \Cross\;\Ric\ .
  $$
  In particular the $L^2$--norm of the covariant derivative $\nabla R\,\in\,
  \Gamma(\,T^*M\otimes\Curv^-TM\,)$ of the curvature tensor of a compact
  Riemannian manifold $M$ with parallel Ricci tensor equals:
  $$
   |\!|\,\nabla R\,|\!|^2_{T^*\otimes\LP^2T^*\otimes\LP^2T^*}
   \;\;=\;\;
   -\;\int_Mg_{\LP^2T^*\otimes\LP^2T^*}
   (\;q(\,R\,)\,\star\,R,\;R\;)\,|\,\vol_g\,|\ .
  $$
 \end{Lemma}

 \proof
 Using the second Bianchi identity six times and keeping all the
 symmetries of algebraic curvature tensors self evident we obtain
 for every local basis $E_1,\,\ldots,\,E_m$:
 \begin{eqnarray*}
  \lefteqn{(\,\nabla^*\nabla R\,)(\;X,\,Y;\,U,\,V\;)}
  \quad
  &&
  \\
  &=&
  -\;\frac12\;\sum_{\mu\,=\,1}^m
  \Big(\;+\;(\,\nabla^2_{E_\mu,\,dE_\mu^\up}R\,)(\,X,\,Y;\,U,\,V\,)
  \;+\;(\,\nabla^2_{E_\mu,\,dE_\mu^\up}R\,)(\,U,\,V;\,X,\,Y\,)\;\Big)
  \\
  &=&
  +\;\frac12\;\sum_{\mu\,=\,1}^m
  \Big(\;+\;(\,\nabla^2_{E_\mu,\,X}R\,)(\,Y,\,dE_\mu^\up;\,U,\,V\,)
  \;-\;(\,\nabla^2_{E_\mu,\,Y}R\,)(\,X,\,dE_\mu^\up;\,U,\,V\,)
  \\[-10pt]
  &&
  \qquad\qquad\quad\;+\;(\,\nabla^2_{E_\mu,\,U}R\,)(\,V,\,dE_\mu^\up;\,X,\,Y\,)
  \;-\;(\,\nabla^2_{E_\mu,\,V}R\,)(\,U,\,dE_\mu^\up;\,X,\,Y\,)\;\Big)
  \\
  &=&
  +\;\frac12\;\sum_{\mu\,=\,1}^m
  \Big(\;+\;(\,R_{E_\mu,\,X}R\,)(\,Y,\,dE_\mu^\up;\,U,\,V\,)
  \;+\;(\,\nabla^2_{X,\,E_\mu}R\,)(\,U,\,V;\,Y,\,dE_\mu^\up\,)
  \\[-10pt]
  &&
  \qquad\qquad\quad\;-\;(\,R_{E_\mu,\,Y}R\,)(\,X,\,dE_\mu^\up;\,U,\,V\,)
  \;-\;(\,\nabla^2_{Y,\,E_\mu}R\,)(\,U,\,V;\,X,\,dE_\mu^\up\,)
  \\[2pt]
  &&
  \qquad\qquad\quad\;+\;(\,R_{E_\mu,\,U}R\,)(\,V,\,dE_\mu^\up,\,X,\,Y\,)
  \;+\;(\,\nabla^2_{U,\,E_\mu}R\,)(\,X,\,Y;\,V,\,dE_\mu^\up\,)
  \\
  &&
  \qquad\qquad\quad\;-\;(\,R_{E_\mu,\,V}R\,)(\,U,\,dE_\mu^\up,\,X,\,Y\,)
  \;+\;(\,\nabla^2_{V,\,E_\mu}R\,)(\,X,\,Y;\,U,\,dE_\mu^\up\,)\;\Big)
  \\
  &=&
  +\;\frac12\;\sum_{\mu\,=\,1}^m
  \Big(\;+\;(\,\nabla^2_{X,\,U}R\,)(\,V,\,E_\mu;\,dE_\mu^\up,\,Y\,)
  \;-\;(\,\nabla^2_{X,\,V}R\,)(\,U,\,E_\mu;\,dE_\mu^\up,\,Y\,)
  \\[-10pt]
  &&
  \qquad\qquad\quad\;-\;(\,\nabla^2_{Y,\,U}R\,)(\,V,\,E_\mu;\,dE_\mu^\up,\,X\,)
  \;+\;(\,\nabla^2_{Y,\,V}R\,)(\,U,\,E_\mu;\,dE_\mu^\up,\,X\,)
  \\[2pt]
  &&
  \qquad\qquad\quad\;+\;(\,\nabla^2_{U,\,X}R\,)(\,Y,\,E_\mu;\,dE_\mu^\up,\,V\,)
  \;-\;(\,\nabla^2_{U,\,Y}R\,)(\,X,\,E_\mu;\,dE_\mu^\up,\,V\,)
  \\
  &&
  \qquad\qquad\quad\;-\;(\,\nabla^2_{V,\,X}R\,)(\,Y,\,E_\mu;\,dE_\mu^\up,\,U\,)
  \;+\;(\,\nabla^2_{V,\,Y}R\,)(\,X,\,E_\mu;\,dE_\mu^\up,\,U\,)\;\Big)
  \\
  &&
  -\;\frac12\;\sum_{\mu\,=\,1}^m
  \Big(\;+\;(\,R_{E_\mu,\,X}R\,)(\,dE_\mu^\up,\,Y;\,U,\,V\,)
  \;+\;(\,R_{E_\mu,\,Y}R\,)(\,X,\,dE_\mu^\up;\,U,\,V\,)
  \\[-10pt]
  &&
  \qquad\qquad\quad\;+\;(\,R_{E_\mu,\,U}R\,)(\,X,\,Y;\,dE_\mu^\up,\,V\,)
  \;+\;(\,R_{E_\mu,\,V}R\,)(\,X,\,Y;\,U,\,dE_\mu^\up\,)\;\Big)\ .
 \end{eqnarray*}
 Evidently the first sum on the right hand side is just $\Cross\,\Ric$,
 note in particular that the symmetrized iterated covariant derivatives
 $\nabla^{[2]}$ appear naturally in this calculation. We leave it to the
 reader to verify that the second sum on the right hand side equals
 $-q(R)\,\star\,R$.
 \qed

 \pfill
 Combining Lemma \ref{pde} with the interpretation (\ref{varr}) of the
 cross operator as the linearization of the second order differential operator
 $\Gamma(\,\S^2_{\mathrm{reg}}T^*M\,)\longrightarrow\Gamma(\,\Curv^-TM\,),
 \,g\longmapsto R,$ we obtain the central formula of Hamilton's theory
 \cite{br} for the Ricci Flow $\delta g\,:=\,\pm\,2\,\Ric$:
 $$
  \delta R
  \;\;=\;\;
  \pm\;\Big(\;\Cross\;\Ric\;+\;\frac12\,\Der_\Ric R\;\Big)
  \;\;=\;\;
  \pm\;\Big(\;\nabla^*\nabla R\;+\;q(\,R\,)\star R
  \;+\;\frac12\,\Der_\Ric R\;\Big)\ .
 $$
 Comparing this formulation with the usual formulation of Hamilton's
 theory \cite{br} we find
 \begin{equation}\label{etc}
  q(\,R\,)\,\star\,R
  \;\;=\;\;
  \frac12\,\Der_\Ric R\;+\;Q(\,R\,)
 \end{equation}
 in terms of Hamilton's quadratic curvature expression $Q(\,R\,)$ defined by:
 \begin{eqnarray*}
  \lefteqn{Q(\,R\,)(\,X,\,Y;\,U,\,V\,)}
  &&
  \\
  &:=&
  2\,g_{\LP^2T}(\,R_{X,\,Y},\,R_{U,\,V}\,)
  \;+\;2\,\sum_{\mu\,=\,1}^m\Big(\;g(\,R_{E_\mu,\,X}U,\,R_{dE_\mu^\up,\,Y}V\,)
  \;-\;g(\,R_{E_\mu,\,X}V,\,R_{dE_\mu^\up,\,Y}U\,)\;\Big)\ .
 \end{eqnarray*}
 In order to illustrate the power of the graphical calculus developed in this
 article we want to discuss an interesting identity of stable curvature
 invariants of Einstein manifolds related to Lemma \ref{pde}. Recall first
 of all that an algebraic curvature tensor of Einstein type is an algebraic
 curvature tensor $R\,\in\,\Curv^-T$ on a euclidean vector space $T$ of
 dimension $m\,\in\,\N$ whose associated Ricci tensor $\Ric\,=\,\frac\kappa m
 \,g$ is a multiple of the metric $g$. On the subspace $\Ein(\,T,\,g\,)\,
 \subseteq\,\Curv^-T$ of curvature tensors of Einstein type the Ricci
 identity (\ref{ihex})
 $$
  \raise-12pt\hbox{\begin{picture}(40,30)
   \put(10,25){\circle2}
   \put(10, 5){\circle2}
   \put(11,25){\line(+1, 0){14}}
   \put(11, 5){\line(+1, 0){14}}
   \put(27,21){$X$}
   \put(27, 1){$Y$}
   \multiput(10, 7)(0,+2){9}{\circle*1}
   \qbezier( 9,25)( 2,24)(2,15)
   \qbezier( 9, 5)( 2, 6)(2,15)
  \end{picture}}
  \;\;\widehat=\;\;
  -\;2\,\frac\kappa m\,g(\,X,\,Y\,)
  \;\;\widehat=\;\;
  -\;2\,\frac\kappa m\,\raise-12pt\hbox{\begin{picture}(16,30)
   \put(10,21){$X$}
   \put(10, 1){$Y$}
   \qbezier( 9,25)( 2,24)(2,15)
   \qbezier( 9, 5)( 2, 6)(2,15)
  \end{picture}}
 $$
 becomes an algebraic simplification on the level of the reduced graph algebra
 $\overline\A^\bullet$ reducing the number of vertices at the expense of
 introducing $\kappa$ as a new independent variable. Generators of the reduced
 graph algebra $\overline\A^\bullet$ with a pair of parallel red and black
 edges are thus redundant when restricted to $\Ein(\,T,\,g\,)$, for the first
 two of the generators (\ref{gen6}) of degree $3$ we find say
 \begin{equation}\label{R1}
  \left[\;\raise-11pt\hbox{\begin{picture}(30,30)
   \put(15, 2){\circle2}
   \put(15,28){\circle2}
   \put( 3, 8){\circle2}
   \put(27, 8){\circle2}
   \put( 3,22){\circle2}
   \put(27,22){\circle2}
   \put( 3.0, 9.0){\line( 0,+1){12}}
   \put(27.0, 9.0){\line( 0,+1){12}}
   \put(15.9, 2.5){\line(+2,+1){10}}
   \put(14.1, 2.5){\line(-2,+1){10}}
   \put(15.9,27.5){\line(+2,-1){10}}
   \put(14.1,27.5){\line(-2,-1){10}}
   \multiput( 4.9,21.9)(+1.5,+0.75){7}{\circle*1}
   \multiput( 4.9, 8.1)(+1.5,-0.75){7}{\circle*1}
   \multiput(26.1,10.1)( 0.0,+2.0){6}{\circle*1}
  \end{picture}}\;\right](\,R\,)
  \;\;\stackrel!=\;\;
  \left[\;\raise-11pt\hbox{\begin{picture}(30,30)
   \put(15, 2){\circle2}
   \put(15,28){\circle2}
   \put( 3, 8){\circle2}
   \put(27, 8){\circle2}
   \put( 3,22){\circle2}
   \put(27,22){\circle2}
   \put( 3.0, 9.0){\line( 0,+1){12}}
   \put(27.0, 9.0){\line( 0,+1){12}}
   \put(15.9, 2.5){\line(+2,+1){10}}
   \put(14.1, 2.5){\line(-2,+1){10}}
   \put(15.9,27.5){\line(+2,-1){10}}
   \put(14.1,27.5){\line(-2,-1){10}}
   \multiput( 3.9,10.1)( 0.0,+2.0){6}{\circle*1}
   \multiput(15.0, 4.0)( 0.0,+2.0){12}{\circle*1}
   \multiput(26.1,10.1)( 0.0,+2.0){6}{\circle*1}
  \end{picture}}\;\right](\,R\,)
  \;\;=\;\;
  -\;2\,\frac\kappa m\;
  \left[\;\raise-9pt\hbox{\begin{picture}(26,26)
   \put( 2, 2){\circle2}
   \put(24, 2){\circle2}
   \put( 2,24){\circle2}
   \put(24,24){\circle2}
   \put( 2.0, 3.0){\line( 0,+1){20}}
   \put(24.0, 3.0){\line( 0,+1){20}}
   \put( 3.0, 2.0){\line(+1, 0){20}}
   \put( 3.0,24.0){\line(+1, 0){20}}
   \multiput( 3.2, 3.9)( 0.0,+2.0){10}{\circle*1}
   \multiput(22.8, 3.9)( 0.0,+2.0){10}{\circle*1}
  \end{picture}}\;\right](\,R\,)
  \;\;=\;\;
  -\;8\,\frac{\kappa^3}{m^2}
 \end{equation}
 for every algebraic curvature tensor $R\,\in\,\Ein(\,T,\,g\,)$ of Einstein
 type in light of equation (\ref{deg2}) in combination with $g_{\S^2T^*}
 (\frac\kappa mg,\frac\kappa mg)\,=\,\frac{\kappa^2}{2m}$. Similarly the
 third generator reduces to:
 \begin{equation}\label{R2}
  \left[\;\raise-11pt\hbox{\begin{picture}(30,30)
   \put(15, 2){\circle2}
   \put(15,28){\circle2}
   \put( 3, 8){\circle2}
   \put(27, 8){\circle2}
   \put( 3,22){\circle2}
   \put(27,22){\circle2}
   \put( 3.0, 9.0){\line( 0,+1){12}}
   \put(27.0, 9.0){\line( 0,+1){12}}
   \put(15.9, 2.5){\line(+2,+1){10}}
   \put(14.1, 2.5){\line(-2,+1){10}}
   \put(15.9,27.5){\line(+2,-1){10}}
   \put(14.1,27.5){\line(-2,-1){10}}
   \multiput( 4.9,21.9)(+1.5,+0.75){7}{\circle*1}
   \multiput(16.1, 3.8)(+0.9,+1.5){12}{\circle*1}
   \multiput( 5.0, 8.0)(+2.0, 0.0){11}{\circle*1}
  \end{picture}}\;\right](\,R\,)
  \;\;=\;\;
  -\;2\,\frac\kappa m\;
  \left[\;\raise-9pt\hbox{\begin{picture}(26,26)
   \put( 2, 2){\circle2}
   \put(24, 2){\circle2}
   \put( 2,24){\circle2}
   \put(24,24){\circle2}
   \put( 2.0, 3.0){\line( 0,+1){20}}
   \put(24.0, 3.0){\line( 0,+1){20}}
   \put( 3.0, 2.0){\line(+1, 0){20}}
   \put( 3.0,24.0){\line(+1, 0){20}}
   \multiput( 3.4, 3.4)(+1.2,+1.2){17}{\circle*1}
   \multiput(22.6, 3.4)(-1.2,+1.2){17}{\circle*1}
  \end{picture}}\;\right](\,R\,)
  \;\;=\;\;
  +\;48\,\frac\kappa m\,g_{\LP^2T^*\otimes\LP^2T^*}(\,R,\,R\,)\ .
 \end{equation}
 Effectively we are thus left with $4$ instead of the $8$ generators
 (\ref{gen6}) of degree up to $3$, namely
 $$
  \raise-9pt\hbox{\begin{picture}(10,26)
   \put( 5, 2){\circle2}
   \put( 5,24){\circle2}
   \qbezier(4.0, 2.0)( 2.0, 2.0)( 2.0,12.0)
   \qbezier(4.0,24.0)( 2.0,24.0)( 2.0,12.0)
   \qbezier(6.1, 2.0)( 8.0, 2.0)( 8.0,12.0)
   \qbezier(6.1,24.0)( 8.0,24.0)( 8.0,12.0)
   \multiput( 5.0, 4.0)( 0.0,+2.0){10}{\circle*1}
  \end{picture}}
  \;\;\widehat=\;\;-\;2\,\kappa
  \qquad\quad
  \raise-9pt\hbox{\begin{picture}(26,26)
   \put( 2, 2){\circle2}
   \put(24, 2){\circle2}
   \put( 2,24){\circle2}
   \put(24,24){\circle2}
   \put( 2.0, 3.0){\line( 0,+1){20}}
   \put(24.0, 3.0){\line( 0,+1){20}}
   \put( 3.0, 2.0){\line(+1, 0){20}}
   \put( 3.0,24.0){\line(+1, 0){20}}
   \multiput( 3.4, 3.4)(+1.2,+1.2){17}{\circle*1}
   \multiput(22.6, 3.4)(-1.2,+1.2){17}{\circle*1}
  \end{picture}}
  \;\;\widehat=\;\;
  -\;24\,g_{\LP^2T^*\otimes\LP^2T^*}(\,R,\,R\,)
  \qquad\quad
  \raise-11pt\hbox{\begin{picture}(30,30)
   \put(15, 2){\circle2}
   \put(15,28){\circle2}
   \put( 3, 8){\circle2}
   \put(27, 8){\circle2}
   \put( 3,22){\circle2}
   \put(27,22){\circle2}
   \put( 3.0, 9.0){\line( 0,+1){12}}
   \put(27.0, 9.0){\line( 0,+1){12}}
   \put(15.9, 2.5){\line(+2,+1){10}}
   \put(14.1, 2.5){\line(-2,+1){10}}
   \put(15.9,27.5){\line(+2,-1){10}}
   \put(14.1,27.5){\line(-2,-1){10}}
   \multiput(15.0, 4.0)( 0.0,+2.0){12}{\circle*1}
   \multiput( 4.9, 8.2)(+2.0, 0.0){11}{\circle*1}
   \multiput( 4.9,21.8)(+2.0, 0.0){11}{\circle*1}
  \end{picture}}
  \qquad\quad
  \raise-11pt\hbox{\begin{picture}(30,30)
   \put(15, 2){\circle2}
   \put(15,28){\circle2}
   \put( 3, 8){\circle2}
   \put(27, 8){\circle2}
   \put( 3,22){\circle2}
   \put(27,22){\circle2}
   \put( 3.0, 9.0){\line( 0,+1){12}}
   \put(27.0, 9.0){\line( 0,+1){12}}
   \put(15.9, 2.5){\line(+2,+1){10}}
   \put(14.1, 2.5){\line(-2,+1){10}}
   \put(15.9,27.5){\line(+2,-1){10}}
   \put(14.1,27.5){\line(-2,-1){10}}
   \multiput(15.0, 4.0)( 0.0,+2.0){12}{\circle*1}
   \multiput( 4.7, 9.0)(+1.9,+1.1){12}{\circle*1}
   \multiput(25.3, 9.0)(-1.9,+1.1){12}{\circle*1}
  \end{picture}}\ ,
 $$
 when discussing stable curvature invariants of algebraic curvature tensors
 of Einstein type. In consequence there exists up to scale a unique linear
 combination of the Pfaffian $\pf_3$, the normalized moment polynomial
 $\Psi^\circ_3$ and the cubic polynomial $\Theta_3$ defined in Lemma
 \ref{cubic}
 \begin{eqnarray*}
  \pf_3(R)
  &=&
  \left[\,-\,\hbox to20pt{\hfill$\displaystyle\frac1{432}$\hfill}\;
  \raise-11pt\hbox{\begin{picture}(30,30)
   \put(15, 2){\circle2}
   \put(15,28){\circle2}
   \put( 3, 8){\circle2}
   \put(27, 8){\circle2}
   \put( 3,22){\circle2}
   \put(27,22){\circle2}
   \put( 3.0, 9.0){\line( 0,+1){12}}
   \put(27.0, 9.0){\line( 0,+1){12}}
   \put(15.9, 2.5){\line(+2,+1){10}}
   \put(14.1, 2.5){\line(-2,+1){10}}
   \put(15.9,27.5){\line(+2,-1){10}}
   \put(14.1,27.5){\line(-2,-1){10}}
   \multiput(15.0, 4.0)( 0.0,+2.0){12}{\circle*1}
   \multiput( 4.7, 9.0)(+1.9,+1.1){12}{\circle*1}
   \multiput(25.3, 9.0)(-1.9,+1.1){12}{\circle*1}
  \end{picture}}
  \,-\,
  \hbox to20pt{\hfill$\displaystyle\frac5{432}$\hfill}\;
  \raise-11pt\hbox{\begin{picture}(30,30)
   \put(15, 2){\circle2}
   \put(15,28){\circle2}
   \put( 3, 8){\circle2}
   \put(27, 8){\circle2}
   \put( 3,22){\circle2}
   \put(27,22){\circle2}
   \put( 3.0, 9.0){\line( 0,+1){12}}
   \put(27.0, 9.0){\line( 0,+1){12}}
   \put(15.9, 2.5){\line(+2,+1){10}}
   \put(14.1, 2.5){\line(-2,+1){10}}
   \put(15.9,27.5){\line(+2,-1){10}}
   \put(14.1,27.5){\line(-2,-1){10}}
   \multiput(15.0, 4.0)( 0.0,+2.0){12}{\circle*1}
   \multiput( 4.9, 8.2)(+2.0, 0.0){11}{\circle*1}
   \multiput( 4.9,21.8)(+2.0, 0.0){11}{\circle*1}
  \end{picture}}\,\right](\,R\,)
  \;+\;\kappa^3\,\frac{m^2-12m+40}{48\,m^2}
  \;+\;\kappa\,\frac{m-8}{4\,m}\,|\,R\,|^2
  \\
  \Psi^\circ_3(R)
  &=&
  \left[\,-\,\hbox to20pt{\hfill$\displaystyle\frac16$\hfill}\;
  \raise-11pt\hbox{\begin{picture}(30,30)
   \put(15, 2){\circle2}
   \put(15,28){\circle2}
   \put( 3, 8){\circle2}
   \put(27, 8){\circle2}
   \put( 3,22){\circle2}
   \put(27,22){\circle2}
   \put( 3.0, 9.0){\line( 0,+1){12}}
   \put(27.0, 9.0){\line( 0,+1){12}}
   \put(15.9, 2.5){\line(+2,+1){10}}
   \put(14.1, 2.5){\line(-2,+1){10}}
   \put(15.9,27.5){\line(+2,-1){10}}
   \put(14.1,27.5){\line(-2,-1){10}}
   \multiput(15.0, 4.0)( 0.0,+2.0){12}{\circle*1}
   \multiput( 4.7, 9.0)(+1.9,+1.1){12}{\circle*1}
   \multiput(25.3, 9.0)(-1.9,+1.1){12}{\circle*1}
  \end{picture}}
  \,-\,
  \hbox to20pt{\hfill$\displaystyle\frac16$\hfill}\;
  \raise-11pt\hbox{\begin{picture}(30,30)
   \put(15, 2){\circle2}
   \put(15,28){\circle2}
   \put( 3, 8){\circle2}
   \put(27, 8){\circle2}
   \put( 3,22){\circle2}
   \put(27,22){\circle2}
   \put( 3.0, 9.0){\line( 0,+1){12}}
   \put(27.0, 9.0){\line( 0,+1){12}}
   \put(15.9, 2.5){\line(+2,+1){10}}
   \put(14.1, 2.5){\line(-2,+1){10}}
   \put(15.9,27.5){\line(+2,-1){10}}
   \put(14.1,27.5){\line(-2,-1){10}}
   \multiput(15.0, 4.0)( 0.0,+2.0){12}{\circle*1}
   \multiput( 4.9, 8.2)(+2.0, 0.0){11}{\circle*1}
   \multiput( 4.9,21.8)(+2.0, 0.0){11}{\circle*1}
  \end{picture}}\,\right](\,R\,)
  \;+\;\kappa^3\,\frac{m^2+12m+40}{6\,m^2}
  \;+\;\kappa\,\frac{6\,(m+8)}m\,|\,R\,|^2
  \\
  \Theta_3(R)
  &=&
  \left[\,+\,\hbox to20pt{\hfill$\displaystyle\frac1{72}$\hfill}\;
  \raise-11pt\hbox{\begin{picture}(30,30)
   \put(15, 2){\circle2}
   \put(15,28){\circle2}
   \put( 3, 8){\circle2}
   \put(27, 8){\circle2}
   \put( 3,22){\circle2}
   \put(27,22){\circle2}
   \put( 3.0, 9.0){\line( 0,+1){12}}
   \put(27.0, 9.0){\line( 0,+1){12}}
   \put(15.9, 2.5){\line(+2,+1){10}}
   \put(14.1, 2.5){\line(-2,+1){10}}
   \put(15.9,27.5){\line(+2,-1){10}}
   \put(14.1,27.5){\line(-2,-1){10}}
   \multiput(15.0, 4.0)( 0.0,+2.0){12}{\circle*1}
   \multiput( 4.7, 9.0)(+1.9,+1.1){12}{\circle*1}
   \multiput(25.3, 9.0)(-1.9,+1.1){12}{\circle*1}
  \end{picture}}
  \,-\,
  \hbox to20pt{\hfill$\displaystyle\frac1{18}$\hfill}\;
  \raise-11pt\hbox{\begin{picture}(30,30)
   \put(15, 2){\circle2}
   \put(15,28){\circle2}
   \put( 3, 8){\circle2}
   \put(27, 8){\circle2}
   \put( 3,22){\circle2}
   \put(27,22){\circle2}
   \put( 3.0, 9.0){\line( 0,+1){12}}
   \put(27.0, 9.0){\line( 0,+1){12}}
   \put(15.9, 2.5){\line(+2,+1){10}}
   \put(14.1, 2.5){\line(-2,+1){10}}
   \put(15.9,27.5){\line(+2,-1){10}}
   \put(14.1,27.5){\line(-2,-1){10}}
   \multiput(15.0, 4.0)( 0.0,+2.0){12}{\circle*1}
   \multiput( 4.9, 8.2)(+2.0, 0.0){11}{\circle*1}
   \multiput( 4.9,21.8)(+2.0, 0.0){11}{\circle*1}
  \end{picture}}\,\right](\,R\,)
  \;+\;
  \kappa\,\frac 2m\,|\,R\,|^2
 \end{eqnarray*}
 in which both additional cubic generators are eliminated, the resulting
 identity reads:
 $$
  \Big(\;\pf_3\;-\;\frac1{40}\,\Psi^\circ_3\;-\;\frac2{15}\,\Theta_3\;\Big)
  (\,R\,)
  \;\;=\;\;
  \kappa^3\,\frac{m^2-18m+40}{60\,m^2}\;+\;\kappa\,\frac{3m-104}{30\,m}
  g_{\LP^2T^*\otimes\LP^2T^*}(\,R,\,R\,)\ .
 $$
 According to Lemma \ref{pde} the cubic curvature invariant $\Theta_3(\,R\,)
 \,:=\,-g_{\LP^2T^*\otimes\LP^2T^*}(\,q(R)\star R,\,R\,)$ integrates for a
 compact Einstein manifold $M$ to the $L^2$--norm $|\!|\,\nabla R\,
 |\!|^2_{T^*\otimes\LP^2T^*\otimes\LP^2T^*}$ of the covariant derivative of
 the curvature tensor, in this way we have proved:
 
 \begin{Theorem}[Cubic Curvature Identity for Einstein Manifolds]
 \hfill\label{ineq}\break
  For every compact connected Einstein manifold $M$ of dimension $m\,\geq\,3$
  with scalar curvature $\kappa\,\in\,\R$ the following identity of integrated
  stable curvature invariants of degree $3$ holds true:
  \begin{eqnarray*}
   \lefteqn{\int_M\pf_3(\,R\,)\,|\,\vol_g\,|
    \;-\;\frac1{40}\,\int_M\Psi^\circ_3(\,R\,)\,|\,\vol_g\,|
    \;+\;\frac2{15}\,|\!|\,\nabla R\,|\!|^2_{T^*\otimes\LP^2T^*\otimes\LP^2T^*}}
   \qquad\qquad
   &&
   \\
   &=&
   \kappa^3\,\frac{m^2-18m+40}{60\,m^2}\,\Vol(\,M,\,g\,)\;+\;
   \kappa\,\frac{3m-104}{30\,m}\,|\!|\,R\,|\!|^2_{\LP^2T^*\otimes\LP^2T^*}
  \end{eqnarray*}
  Specifically in dimension $m\,=\,5$ this identity reduces via
  $\Psi^\circ_3\,=\,10080\,\Psi_3$ to the identity
  $$
   -\,252\int_M\Psi_3(R)\,|\,\vol_g\,|\,+\,
   \frac2{15}\,|\!|\,\nabla R\,|\!|^2_{T^*\otimes\LP^2T^*\otimes\LP^2T^*}
   \;\;=\;\;
   -\,\frac{\kappa^3}{60}\,\Vol(\,M,\,g\,)\,-\,
   \frac{89\,\kappa}{150}\,|\!|\,R\,|\!|^2_{\LP^2T^*\otimes\LP^2T^*}\ ,
  $$
  because $\pf_3(\,R\,)\,=\,0$ due to Remark \ref{vsci}, whereas is becomes
  for $m\,=\,6$ with $\Psi^\circ_3\,=\,25200\,\Psi_3$:
  \begin{eqnarray*}
   \lefteqn{(\,2\pi\,)^3\,\chi(\,M\,)
    \;-\;630\,\int_M\Psi_3(\,R\,)\,|\,\vol_g\,|
    \;+\;\frac2{15}\,|\!|\,\nabla R\,|\!|^2_{T^*\otimes\LP^2T^*\otimes\LP^2T^*}}
   \qquad\qquad
   &&
   \\
   &=&
   -\;\frac{2\,\kappa^3}{135}\,\Vol(\,M,\,g\,)
   \;-\;\frac{43\,\kappa}{90}\,|\!|\,R\,|\!|^2_{\LP^2T^*\otimes\LP^2T^*}\ .
  \end{eqnarray*}
 \end{Theorem}

\begin{thebibliography}{10}
%
\bibitem{b}
 \textsc{Berger, M.~:}\ 
 \textit{Sur quelques variétés d'Einstein compactes},\ 
 \textrm{Annali di Matematica Pura ed Applicata \textbf{53},
  89---95 (1961)}.
%
\bibitem{besse}
 \textsc{Besse, A.~:}\ 
 \textit{Einstein Manifolds},\ 
 \textrm{Results in Mathematics and Related Areas \textbf{10},
  Springer--Verlag, Berlin (1987)}.
%
\bibitem{bgv}
 \textsc{Berline, N.~, Getzler, E.~\& Vergne, M.~:}\ 
 \textit{Heat Kernels and Dirac Operators},\ 
 \textrm{Grundlehren der mathematischen Wissenschaften \textbf{298},
  Springer--Verlag, Berlin (1992)}.
%
\bibitem{br}
 \textsc{Brendle, S.~:}\ 
 \textit{Ricci Flow and the Sphere Theorem},\ 
 \textrm{Graduate Studies in Mathematics \textbf{111},
  American Mathematical Society, Providence (2010)}.
%
\bibitem{ch}
 \textsc{Chern, S.~S.~:}\ 
 \textit{A Simple Intrinsic Proof of the Gauss--Bonnet Formula for Closed
  Riemannian Manifolds},\ 
 \textrm{Annals in Mathematics, Second Series \textbf{45} (IV),
  747---752 (1944)}.
%
\bibitem{fl}
 \textsc{Flanders, H.~:}\ 
 \textit{Development of an Extended Exterior Differential Calculus},\ 
 \textrm{Transactions of the American Mathematical Society \textbf{75} (II),
  311---326 (1953)}.
%
\bibitem{h}
 \textsc{Hitchin, N.~:}\ 
 \textit{Compact Four--Dimensional Einstein Manifolds},\ 
 \textrm{Journal of Differential Geometry \textbf{9} (III),
  435---441 (1974)}.
%
\bibitem{nw}
 \textsc{Nieper--Wißkirchen, M.~:}\ 
 \textit{Chern Numbers and Rozansky--Witten Invariants of Compact
  Hyperkähler Manifolds},\ 
 \textrm{World Scientific Publishing, Singapore (2004)}.
%
\bibitem{sw}
 \textsc{Semmelmann, U.~\& Weingart, G.~:}\ 
 \textit{The Standard Laplace Operator},\ 
 \textrm{manuscripta mathematica \textbf{158}, 273---293 (2019)}.
%
\bibitem{wa}
 \textsc{Walter, R.~:}\ 
 \textit{Differentialgeometrie},\ 
 \textrm{Cornelsen Verlag, Berlin (1989)}.
%
\bibitem{w1}
 \textsc{Weingart, G.~:}\ 
 \textit{Moments of Sectional Curvature},\ 
 \textrm{preprint \texttt{arXiv:1707.06369}, (2017).}
%
\bibitem{w2}
 \textsc{Weingart, G.~:}\ 
 \textit{Curvature Graphs},\ 
 \textrm{implementation of graph algebras for Maxima,
  \texttt{http://www.matcuer.unam.mx/$\sim$gw/CurvGraphs.mac} (2024).}
%
\end{thebibliography}
\end{document}